\documentclass[a4paper,11pt]{article}

\usepackage{amssymb,latexsym}
\usepackage{amsmath,amsfonts}
\usepackage[french]{babel}
\usepackage[latin1]{inputenc}

\newtheorem{theorem}{Th\'{e}or\`{e}me}[section]
\newtheorem{prop}{Proposition}[section]
\newtheorem{lemma}{Lemme}[section]

\newtheorem{defi}{D\'{e}finition}[section]

\newtheorem{rem}{Remarque}[section]

\newenvironment{prooof}
{ \noindent {D\'{e}monstration: }}
{{~} \hfill  $\Box$ \par\medskip}

\newenvironment{prooofe}
{ \noindent {D\'{e}monstration: }}
{}

\newcommand{\real}{\mathbb R}
\newcommand{\korps}{\mathbb K}

\newcommand{\complex}{\mathbb C}
\newcommand{\nat}{\mathbb N}
\newcommand{\entier}{\mathbb Z}

\newcommand{\ra}{\rightarrow}
\newcommand{\lra}{\longrightarrow}

\newcommand{\lla}{\longleftarrow}
\newcommand{\Cinf}{\mathcal{C}^\infty}

\newcommand{\CinfK}[1]{\mathcal{C}^\infty(#1\!,\!\korps)}
\newcommand{\Ginf}{{\Gamma}^\infty}

\newcommand {\id} {{\mathrm{id}}}
\newcommand {\beas}{\begin{eqnarray*}}
\newcommand {\eeas}{\end{eqnarray*}}
\newcommand {\bea}{\begin{eqnarray}}
\newcommand {\eea}{\end{eqnarray}}
\def \ann {{\scriptstyle{\rm ann}}}

\def \rlvt {{\scriptstyle{\rm rlvt}}}

\newcommand{\drm}{\mathrm{d}}
\newcommand{\dk}{\mathbf{d}}
\def \Hom {\mathop{\hbox{\rm Hom}}\nolimits}
\def \Har {\mathop{\hbox{\rm Har}}\nolimits}
\def \Ker {\mathop{\hbox{\rm Ker}}\nolimits}

\def \ot {\otimes}
\def \QM {\mathbb{Q}}
\newcommand{\korpsL}{\korps [[h]]}

\newcommand{\Al}{\mathcal{A}}
\newcommand{\Bl}{\mathcal{B}}
\newcommand{\Il}{\mathcal{I}}

\newcommand{\AlL}{\mathcal{A}[[h]]}

\newcommand{\Gs}{\mathfrak{G}}
\newcommand{\ls}{\mathfrak{l}}

\newcommand{\Hs}{\mathfrak{H}}
\newcommand{\hs}{\mathfrak{h}}

\newcommand{\GsI}{\mathfrak{G}_{\mathcal{I}}}
\newcommand{\GsP}{\widetilde{\mathfrak{G}}}
\newcommand{\gs}{\mathfrak{g}}
\newcommand{\gsI}{\mathfrak{g}_{\mathcal{I}}}
\newcommand{\gsP}{\widetilde{\mathfrak{g}} }
\newcommand{\as}{\mathfrak{a}}

\def \gu{\ve^\cdot\!\underline{\gsI^{\otimes \cdot}}}

\def \gextgsI{\ve^\cdot_{\gs}\,\big({\gs}\ot \underline{\gsI^{\otimes \cdot}}\big)}

\newcommand{\vea}{{\scriptstyle \boldsymbol{\Lambda}}}
\def \ve{{\boldsymbol{\Lambda}}}

\newcommand{\Dop}{\mathbf{D}}

\newcommand{\cois}{co\"{\i}sotrope}

\newcommand{\dH}{\partial_H}
\newcommand{\dK}{\partial_K}
\makeatletter
\@addtoreset{equation}{section}

\makeatother

\begin{document}

\pagestyle{empty}

\begin{center}

  {\Large \bf Formalit\'{e} $G_\infty$ adapt\'{e}e et star-repr\'{e}sentations sur
   des sous-vari\'{e}t\'{e}s co\"{\i}sotropes}

\vspace{0.8cm}

      {\bf Martin Bordemann}\\
{\small Laboratoire de Math\'{e}matiques,
        Universit\'{e} de Haute Alsace, Mulhouse   \\
        e--mail: {\tt Martin.Bordemann@uha.fr}  } \\

\vspace{0.4cm}

      {\bf Gr\'{e}gory Ginot}   \\
{\small Laboratoire Analyse G\'{e}om\'{e}trie et Applications,
     Universit\'{e} Paris~13 \\
       Centre des Math\'{e}matiques et de Leurs Applications, ENS Cachan,\\
     e--mail: {\tt Gregory.Ginot@cmla.ens-cachan.fr} } \\

\vspace{0.4cm}

      {\bf Gilles Halbout}   \\
{\small Institut de Recherche Math\'{e}matique Avanc\'{e}e de l'ULP Strasbourg \\
        e-mail: {\tt halbout@math.u-strasbg.fr}  }\\

\vspace{0.4cm}

      {\bf Hans-Christian Herbig}\\
{\small Fachbereich Mathematik, Johann-Wolfgang-Goethe-Universit\"{a}t, Frankfurt a.M. \\
        e-mail: {\tt herbig@math.uni-frankfurt.de}  } \\

\vspace{0.4cm}

      {\bf Stefan Waldmann}  \\
   {\small Fakult\"{a}t f\"{u}r Mathematik und Physik,
          Albert-Ludwigs Universit\"{a}t Freiburg i.~Br.\\
       e-mail: {\tt Stefan.Waldmann@physik.uni-freiburg.de} }

\markboth{Bordemann-Ginot-Halbout-Herbig-Waldmann}
{Repr\'{e}sentations des star-produits sur des sous-vari\'{e}t\'{e}s co\"{\i}sotropes}


\vspace{0.4cm}
\vspace{0.4cm}

{\bf R\'{e}sum\'{e}}\\[5mm]

\begin{minipage}{13cm}
Soit $X$ une vari\'{e}t\'{e} de Poisson et $C$ une sous-vari\'{e}t\'{e}
co\"{\i}sotrope par rapport au crochet de Poisson. Soient
$\Al=\CinfK{X}$ et
$\Il$ l'id\'{e}al des fonctions nulles sur $C$.
Dans ce travail, nous
construisons des star-produits $\star$ sur $X$ qui sont
{\em adapt\'{e}s \`{a} $C$}, i.e.
pour lesquels $\Il[[h]]$ est un id\'{e}al \`{a} gauche de l'alg\`{e}bre
d\'{e}form\'{e}e $(\Al[[h]],\star)$.
Nous obtenons ainsi une {\em repr\'{e}sentation de $(\Al[[h]],\star)$ sur
$\CinfK{C}[[h]]=:\Bl[[h]]=\Al[[h]]/\Il[[h]]$},
d\'{e}formant la repr\'{e}sentation naturelle
de $\Al$ sur $C^\infty(C)$ pour la situation `locale' $X=\real^n$ et
$C=\real^{n-l}$. Ce r\'{e}sultat se d\'{e}duit d'une
g\'{e}n\'{e}ralisation de la conjecture de formalit\'{e} de Tamarkin appliqu\'{e}e
aux cocha\^{\i}nes adapt\'{e}es \`{a} $C$:
nous d\'{e}montrons d'abord un th\'{e}or\`{e}me \`{a} la Hochschild-Kostant-Rosenberg
entre l'espace $\gsI$ des champs de multivecteurs adapt\'{e}s \`{a} $C$
et  $\GsI$, un espace
naturel d'op\'{e}rateurs multidiff\'{e}rentiels adapt\'{e}s \`{a} $C$.
Ensuite nous mettons en \'{e}vidence l'existence d'une structure $G_\infty$
sur $\GsI$ et montrons enfin que les
obstructions \`{a} l'existence d'une formalit\'{e} sont
contr\^{o}\-l\'{e}es par des groupes de cohomologie. Dans le cas o\`{u}
$X=\real^n$ et
$C=\real^{n-l}$, $l\geq 2$ nous montrons  que ces obstructions sont nulles.
\end{minipage}
\end{center}

\vspace{3mm}

\begin{center}
\begin{minipage}{110mm}\footnotesize{\bf Keywords~:}
Deformation quantization,
star-product, homology \end{minipage}
\end{center}

\begin{center}
\begin{minipage}{110mm}\footnotesize{\bf AMS Classification~:}
Primary 16E40, 53D55, Secondary 18D50, 16S80
\end{minipage}
\end{center}

\newpage
\tableofcontents
\newpage

\pagestyle{plain}

\section*{Introduction}

\addcontentsline{toc}{section}{Introduction}

La th\'{e}orie de quantification par d\'{e}formation, con\c{c}ue en 1978 par Bayen,
Flato, Fr{\o}nsdal, Lichnerowicz et Sternheimer \cite{BFFLS78}, est
maintenant tr\`{e}s bien \'{e}tablie, surtout depuis l'article de Kontsevitch
\cite{Kon03} : pour d\'{e}crire l'alg\`{e}bre des observables quantiques, on
cherche \`{a} construire une d\'{e}formation formelle associative $*$, dite
star-produit, de l'alg\`{e}bre $\Al=\CinfK{X}$ des fonctions de
classe $\Cinf$ \`{a} valeurs dans $\korps=\real$ ou $\korps=\complex$ d\'{e}finie
sur une vari\'{e}t\'{e} de Poisson $(X,P)$. On exige que le commutateur d'ordre $1$
dans le param\`{e}tre de d\'{e}formation $h$ soit proportionnel au crochet de
Poisson d\'{e}fini par $P$. Plus pr\'{e}cis\'{e}ment,
un star-produit $\star$ est
une multiplication associative $\korps [[h]]$-bilin\'{e}aire sur le
$\korpsL$-module $\AlL$ telle que pour tous $f,g \in \Al$,

\begin{center}
 \begin{minipage}{10cm}
\begin{description}
\item[(i)] $f \star g = fg + h \mathsf{C}_1(f,g) + \sum_{k\geq 2} h^k \mathsf{C}_k(f,g)$,
\item[(ii)] $\mathsf{C}_1(f,g)-\mathsf{C}_1(g,f) = P(\drm f,\drm g)$,
\item[(iii)] Les $ (\mathsf{C}_r)_{r \geq 1}$ sont des op\'{e}rateurs bidiff\'{e}rentiels.
\item[(iv)] $\mathsf{C}_r(f,1)=0=\mathsf{C}_r(1,f)$ quel que soit $r\geq 1$.
\end{description}
\end{minipage}
\end{center}

L'existence d'un star-produit s'\'{e}tait av\'{e}r\'{e}e tr\`{e}s difficile: pour les
vari\'{e}t\'{e}s symplectiques le th\'{e}or\`{e}me g\'{e}n\'{e}ral a \'{e}t\'{e} montr\'{e} par DeWilde et
Lecomte en 1983 \cite{DL83} et pour les vari\'{e}t\'{e}s de Poisson par
Kontsevitch \cite{Kon03}. La classification \`{a} \'{e}quivalence pr\`{e}s a \'{e}t\'{e} faite
par Deligne, Nest-Tsygan et Bertelsson-Cahen-Gutt pour le cas symplectique
et par Kontsevitch \cite{Kon03} pour le cas g\'{e}n\'{e}ral.

Une des racines historiques de la quantification par d\'{e}formation \'{e}tait sans doute le
calcul symbolique des op\'{e}rateurs diff\'{e}rentiels sur une vari\'{e}t\'{e} diff\'{e}rentielle $Q$ auquel
correspond le fibr\'{e} cotangent $T^*Q$ comme vari\'{e}t\'{e} symplectique. Dans ce
contexte, l'alg\`{e}bre d\'{e}form\'{e}e $\CinfK{T^*Q}[[h]]$ venait toujours avec un module
\`{a} gauche, \`{a} savoir l'espace des fonctions $\CinfK{Q}[[h]]$ sur lequel
elle agit comme op\'{e}rateurs diff\'{e}rentiels. En g\'{e}n\'{e}ral, la question
alg\'{e}brique d'\'{e}tudier {\em les modules} ou les {\em repr\'{e}sentations}
$\rho:\Al[[h]]\otimes \mathcal{M}[[h]]\ra \mathcal{M}[[h]]$
d'une alg\`{e}bre d\'{e}form\'{e}e est tr\`{e}s
int\'{e}ressante et --en toute g\'{e}n\'{e}ralit\'{e}-- toujours ouverte.
Quelques cas particuliers ont \'{e}t\'{e} \'{e}tudi\'{e}s par Fedosov (modules
projectifs, voir \cite{Fed96}), M.Bordemann et S.Waldmann (d\'{e}finition de la construction de
Gel'fand-Naimark-Segal (GNS) dans le cadre de la quantification par d\'{e}formation, voir
\cite{BW98}), H.Bursztyn et S.Waldmann (\'{e}tude de Morita, voir \cite{bursztyn.waldmann:2002a},
\cite{bursztyn.waldmann:2003a:pre}, \cite{bursztyn.waldmann:2004a},
\cite{bursztyn:2002a}, \cite{bursztyn.waldmann:2000b} et le review \cite{waldmann:2005a}).
Pour une \'{e}tude g\'{e}n\'{e}rale des modules, des morphismes et la r\'{e}duction des
star-produits et quelques r\'{e}sultats dans le cas symplectique, voir \cite{Bor03}.

Une classe de repr\'{e}sentations diff\'{e}rentielles importante est donn\'{e}e par
l'espace $\Bl=\CinfK{C}$ pour une vari\'{e}t\'{e} diff\'{e}rentielle $C$ quelconque:
on peut montrer qu'une condition n\'{e}cessaire \`{a} l'existence d'une repr\'{e}sentation
de $(\CinfK{X}[[h]],*)$ sur $\Bl[[h]]$ en tant qu'op\'{e}rateurs
diff\'{e}rentiels est la donn\'{e}e d'une application de classe $\Cinf$
$i:C\ra X$, {\em co\"{\i}sotrope} par rapport \`{a} la structure de
Poisson $P$ sur $X$: ceci veut dire que {\em l'id\'{e}al annulateur} $\Il:=\Ker i^*$
est une sous-alg\`{e}bre de Poisson de $\Al$. Les sous-vari\'{e}t\'{e}s
co\"{\i}sotropes d'une vari\'{e}t\'{e} de Poisson ont la particularit\'{e} d'\^{e}tre
toujours munies d'un feuilletage (en g\'{e}n\'{e}ral singulier), et l'espace
des feuilles peut \^{e}tre consid\'{e}r\'{e} comme vari\'{e}t\'{e} de Poisson ({\em r\'{e}duction
symplectique}). Nous \'{e}tudions ici le probl\`{e}me r\'{e}ciproque :
la {\em quantification des sous-vari\'{e}t\'{e}s
co\"{\i}sotropes}, c'est-\`{a}-dire la construction d'une structure
de module \`{a} gauche $\mathcal{B}[[h]]$
pour l'alg\`{e}bre d\'{e}form\'{e}e $\Al[[h]]$.
Un star-produit $\star$ sera alors  dit {\em repr\'{e}sentable sur une
sous-vari\'{e}t\'{e} $C$ de $X$} lorsqu'il existe une repr\'{e}sentation de l'alg\`{e}bre
d\'{e}form\'{e}e sur $\Bl[[h]]$.
L'existence de tels star-produits, une fois \'{e}tablie, nous
procure d'un bon candidat
pour {\em l'alg\`{e}bre r\'{e}duite}, \`{a} savoir le {\em commutant} ou l'espace de
tous les homomorphismes de $\Al[[h]]$-modules $\mathcal{B}[[h]]\ra\mathcal{B}[[h]]$
(voir \cite{Bor03} pour des
d\'{e}tails). Ceci correspond \'{e}galement au ``{\em coisotropic creed}'' prononc\'{e}
par Jiang-Hua Lu \cite{Lu93} et d\'{e}crit en partie dans le cadre BRST
dans \cite{BHW00}. Le probl\`{e}me que nous regardons est aussi intimement li\'{e} \`{a} celui de la
{\em quantification des
morphismes de Poisson} (voir \cite{Bor03} pour des d\'{e}tails).

Au lieu de fixer la vari\'{e}t\'{e} de Poisson $(X,P)$ et de chercher des
sous-vari\'{e}t\'{e}s co\"{\i}sotropes par rapport \`{a} $P$, on peut fixer la
sous-vari\'{e}t\'{e} $C$ de $X$ et chercher {\em toutes} les structures de Poisson $P$ telles que
$C$ est co\"{\i}sotrope par rapport \`{a} $P$. Dans ce cas, on parle des
structures de Poisson {\em adapt\'{e}es \`{a} la sous-vari\'{e}t\'{e} $C$}.

Dans le m\^{e}me esprit, on dit qu'un star-produit $\star$ est
{\em adapt\'{e} \`{a} une sous-vari\'{e}t\'{e} $C$ de
$X$} lorsque tous les op\'{e}rateurs bidiff\'{e}rentiels $\mathsf{C}_r$ son tels
que $\mathsf{C}_r(f,g)\in\mathcal{I}$ quel que soit $g\in \Il$. Il
s'ensuit imm\'{e}diatement que dans ce cas l'espace $\Il[[h]]$ est un id\'{e}al \`{a} gauche
de l'alg\`{e}bre d\'{e}form\'{e}e, et par cons\'{e}quent $\Bl[[h]]\cong \Al[[h]]/\Il[[h]]$ est un
module \`{a} gauche pour $\Al[[h]$, donc $\star$ est toujours repr\'{e}sentable.
R\'{e}ciproquement, on peut montrer (voir \cite{Bor03}) que tout star-produit
repr\'{e}sentable est \'{e}quivalent \`{a} un star-produit adapt\'{e}.

Les deux derni\`{e}res d\'{e}finitions se g\'{e}n\'{e}ralisent facilement dans un cadre
utilis\'{e} pour la formalit\'{e}:
\begin{itemize}
\item Dans l'espace $\gs:=\Ginf(X,\Lambda TX)$ des
champs de multivecteurs, on consid\`{e}re le sous-espace $\gsI$ des
{\em champs de multivecteurs adapt\'{e}s \`{a} la sous-vari\'{e}t\'{e} $C$}: $\xi$ de rang
$k$ appartient \`{a} $\gsI$ s'il envoie --vu comme op\'{e}rateur
$k$-diff\'{e}rentiel-- tout choix de $k$ \'{e}l\'{e}ments de $\Il$ dans $\Il$.
\item De la m\^{e}me fa\c{c}on, un op\'{e}rateur $k$-diff\'{e}rentiel dans $\Gs$, (l'espace
des op\'{e}rateurs multi-diff\'{e}rentiels sur $X$) est dit {\em adapt\'{e} \`{a} $C$}
lorsqu'il envoie tout choix $f_1,\ldots,f_{k-1},g\in\Al$ avec $g\in\Il$
dans $\Il$. On note $\GsI$ le sous-espace de tous les op\'{e}rateurs
multidiff\'{e}rentiels adapt\'{e}s \`{a} $C$.
\end{itemize}

\noindent D\'{e}sormais $X=\real^n$ et $C=\real^{n-l}$

Le premier r\'{e}sultat de ce travail est de montrer, dans le paragraphe
\ref{SecHKR}, que $\gsI$ et $\GsI$ h\'{e}ritent de toutes les structures
alg\'{e}briques de $\gs$ et $\Gs$, notamment le crochet de Schouten, le
crochet de Gerstenhaber, et la codiff\'{e}rentielle de Hochschild.
De plus, nous d\'{e}montrons un analogue du th\'{e}or\`{e}me de Hochschild, Kostant,
Rosenberg : $\gsI$ est l'espace de cohomologie de $\GsI$, et nous construisons
une application HKR $\psi^{[1]}:\gsI\ra\GsI$ diff\'{e}rente de l'application uselle
\`{a} cause de l'asym\'{e}trie dans la d\'{e}finition de $\GsI$.

Ensuite, nous d\'{e}montrons l'existence d'un $L_\infty$-morphisme diff\'{e}rentiel entre
$\gsI$ et $\GsI$ en suivant le sch\'{e}ma de la preuve de Tamarkin (\cite{Tam98p}, \cite{GH03})
qui a \'{e}tabli et utilis\'{e} les structures $G_\infty$ pour traiter la conjecture de Deligne
\cite{Del93}. Nous avons
ajout\'{e} un appendice \ref{SecAppendice} pour expliciter,
dans un cadre non-op\'{e}radique, quelques d\'{e}tails et conventions de signe.
Dans le paragraphe \ref{SecFormAdapGinf},
nous d\'{e}montrons que toutes les \'{e}tapes de ce sch\'{e}ma peuvent \^{e}tre reproduites pour
$\gsI$ et $\GsI$ et nous r\'{e}duisons le travail au calcul de la cohomologie
du cocomplexe d'obstructions dans le paragraphe
\ref{SecCalculObstructions}. Nous d\'{e}montrons que les obstructions s'annulent,
et nous obtenons ainsi le r\'{e}sultat principal de cet article d'une {\em formalit\'{e}
$G_\infty$ adapt\'{e}e} --qui a \'{e}t\'{e}
annonc\'{e} sans preuve d\'{e}taill\'{e}e pendant l'Euroconf\'{e}rence PQR en 2003 \`{a}
Bruxelles:

\begin{theorem} \label{TheoremePrincipal}
 Soit $X=\real^n$ et $C=\real^{n-l}$ vue comme sous-vari\'{e}t\'{e} de la mani\`{e}re
 usuelle.
 \begin{enumerate}
  \item $\gsI\cong H\GsI$.
  \item Pour $l\geq 2$: il existe un $G_\infty$-morphisme entre $\gsI$ et $\GsI$.
  \item Pour $l\geq 2$: il existe un $L_\infty$-morphisme entre $\gsI$ et $\GsI$.
  \item Soit $P$ une structure de Poisson adapt\'{e}e \`{a} $C$. Alors il existe
      un star-produit $\star$ adapt\'{e} sur $X$ tel que le commutateur
      d'ordre $1$ en $h$ est \'{e}gal \`{a} $P$.
 \end{enumerate}
\end{theorem}

Les id\'{e}es principales de cet article se sont d\'{e}j\`{a} trouv\'{e}es dans notre
pr\'{e}publication \cite{BGHHW03p}; ici on a ajout\'{e} beaucoup de d\'{e}tails
et notamment la preuve du th\'{e}or\`{e}me \ref{TheoDeuxiemeReduction}.

Nous conjecturons que l'on a la fomalit\'{e} \'{e}galement pour le cas $l=1$~:
une indication forte est le fait qu'il est toujours possible de
repr\'{e}senter un star-produit, ce qui a
\'{e}t\'{e} montr\'{e} par Gl\"{o}{\ss}ner (voir \cite[Lemma 1]{Glo98}, \cite{Bor03}):
\begin{theorem}[Gl\"{o}{\ss}ner 1998]
Si $C$ est une sous-vari\'{e}t\'{e} co\"{\i}sotrope de codimension $1$ dans
$X$ et $*$ un star-produit sur $X$.
Alors on peut construire une star-repr\'{e}sentation.
\end{theorem}
Ce r\'{e}sultat suit d\'{e}j\`{a} d'un calcul des obstructions r\'{e}currentes pour
construire une repr\'{e}sentation ordre par ordre: ces obstructions se
trouvent dans le groupe $HH^2\big(A,\Dop(B,B)\big)$, qui se calcule (voir
le th\'{e}or\`{e}me \ref{TheoCalculCohomol}, \'{e}nonc\'{e} 2) comme
\[
  HH\big(A,\Dop(B,B)\big)\cong \Ginf\big(C,\Lambda(TX/TC)\big).
\]
En codimension $1$ le deuxi\`{e}me groupe de cohomologie de Hochschild
s'annule.

Enfin, ce travail nous donne \'{e}galement un sch\'{e}ma pour d\'{e}crire les obstructions
\`{a} la formalit\'{e} qui peuvent appara\^{\i}tre pour la {\em situation globale} d'une sous-vari\'{e}t\'{e}
$C$ dans une vari\'{e}t\'{e} $X$ quelconque: dans ce cas-l\`{a}, nous conjecturons que tous les
\'{e}nonc\'{e}s du type HKR en paragraphe \ref{SecHKR} ainsi
que le premier th\'{e}or\`{e}me
\ref{propcohomolo} pour r\'{e}duire le cocomplexe des obstructions se
globalisent facilement, tandis que le deuxi\`{e}me th\'{e}or\`{e}me
\ref{TheoDeuxiemeReduction} ne sera plus vrai et
pr\'{e}sentera des v\'{e}ritables obstructions globales. Dans ce cas-l\`{a}, il
faudrait
modifier la codiff\'{e}rentielle $d_L$ (co-induite par le crochet de Schouten sur les champs de
multivecteurs adapt\'{e}s) dans le th\'{e}or\`{e}me de formalit\'{e} par une
codiff\'{e}rentielle $d'_L$ contenant des termes d'ordre sup\'{e}rieur ou \'{e}gal \`{a}
$3$ (voir le th\'{e}or\`{e}me \ref{Theorem 2.1}): on peut donc imaginer que ces op\'{e}rateurs
$k$-diff\'{e}rentiels ($k\geq 3$) pr\'{e}sentent des conditions additionnels \`{a} une structure
de Poisson formelle adapt\'{e}e \`{a} une sous-vari\'{e}t\'{e}. Ceci semble possible,
aux vues des obstructions n\'{e}cessaires jusqu'\`{a} l'ordre $3$ \`{a} la repr\'{e}sentabilit\'{e}
d'un star-produit symplectique li\'{e}es
aux classes d'Atiyah-Molino du feuilletage de $C$ (voir \cite{Bor03}).
Un cadre de globalisation \`{a} la \cite{CFT00} ou \cite{Dol03p} serait dans
doute prometteur.

Signalons une autre approche propos\'{e}e dans \cite{CF05}. Le point de vue est
de quantifier directement le crochet de Poisson sur la vari\'{e}t\'{e} r\'{e}duite
$\underline{C}$, c'est-\`{a}-dire de contruire un star-produit sur
$C^\infty(\underline{C})=N(\Il)/\Il$, o\`{u} $N(\Il)=\{f \in \Al, \{f,\Il\}\subset \Il\}$.
Pour cela, les auteurs montrent d'abord qu'un crochet de Poisson sur $X$ pour
lequel $C$ est co\"{\i}sotrope s'\'{e}tend en une structure $P_\infty$ plate
sur $\gsP=\Gamma\big(C,\Lambda (T_CX/TC)\big)$. Une structure $P_\infty$
veut dire qu'on se donne une famille $(P_i)_{i\geq 0}$:
$A^{\otimes i}\to A$ de multi-d\'{e}rivations faisant de $A$ une alg\`{e}bre
$L_\infty$ et ``plate'' veut dire que $P_0=0$.
Dans le cas ``plat'', $P_1$ est une diff\'{e}rentielle et $P_2$ induit une
structure de Poisson sur la cohomologie $H^0(\Gamma,P_1)\cong
C^\infty(\underline{C})$.

Ils construisent ensuite un morphisme de formalit\'{e} pour $\gsP$, vu comme
alg\`{e}bre de fonctions sur une super-vari\'{e}t\'{e}, en d\'{e}finissant des analogues
pour les champs de tenseurs et op\'{e}rateurs multi-diff\'{e}rentiels.
Ce morphisme est un ``super'' analogue de celui de Kontsevitch. Une structure $P_\infty$
est alors envoy\'{e}e sur une structure $A_\infty$  (c'est-\`{a}-dire la
donn\'{e}e d'une famillie $(a_i)_{i\geq 0}$ de lois v\'{e}rifiant des relations
d'associativit\'{e} \`{a} homotopie pr\`{e}s).
Une telle structure permet d'obtenir une
d\'{e}formation associative de $C^\infty(\underline{C})[[h]]$ \`{a} deux conditions
\begin{itemize}
\item que $a_0=0$. Dans ce cas $a_1$ d\'{e}finit une diff\'{e}rentielle
et $a_2$ induit une structure associative sur $H^0(\Gamma(C,N^*C), a_1)$.
\item que dans ce cas  $H^0(\Gamma(C,N^*C), a_1)\cong C^\infty(\underline{C})[[h]].$
\end{itemize}
Les obstructions pour ces deux conditions sont
respectivement dans les groupes $H^2_\pi(N^*C)$ et $H^1_\pi(N^*C)$ qui sont g\'{e}n\'{e}ralement
non nulles m\^{e}me dans le cas plat lin\'{e}aire. Dans ce dernier cas,
les auteurs retrouvent cependant ``\`{a} la main'' une bonne d\'{e}formation
(\cite{CF04}) (dans ce travail, ils proposaient une construction pour laquelle,
moyennant certaines obstructions, il existe une d\'{e}formation de $X$ adapt\'{e}e
\`{a} la sous-vari\'{e}t\'{e} co\"{\i}sotrope $C$).
La premi\`{e}re condition (envoyer une strcture $P_\infty$ avec $P_0=0$ sur
une structure $A_\infty$ avec $a_0=0$) est similaire \`{a} la condition que nous demandons
\`{a} notre morphisme : envoyer $\gsI$ sur $\GsI$. Dans notre cas, cette
condition est r\'{e}alis\'{e}e dans le cas plat.

\smallskip

\noindent{\bf Notations~:}
Pour deux vari\'{e}t\'{e}s
diff\'{e}rentiables $X$ et $X'$, $\Cinf(X,X')$ d\'{e}signe l'ensemble
de toutes les applications de classe $\Cinf$ de $X$ dans $X'$. Pour un
fibr\'{e} vectoriel $E$ sur une vari\'{e}t\'{e} diff\'{e}rentiable $X$, on \'{e}crira
$\Gamma(X,E)$ pour l'espace de toutes les sections de classe $\Cinf$ du
fibr\'{e} $E$. On \'{e}crira les symboles $\wedge$ et $\Lambda$ pour le produit
ext\'{e}rieur point-par-point dans un cadre de g\'{e}om\'{e}trie diff\'{e}rentielle,
tandis qu'on va utiliser $\scriptstyle \boldsymbol{\Lambda}$ and $\ve$ dans un cadre
alg\`{e}brique, i.e. le produit ext\'{e}rieur (gradu\'{e}) des modules gradu\'{e}s.

\smallskip

\noindent{\bf Remerciements~:}
Nous remercions EUCOR pour des aides financi\`{e}res qui ont rendu possible
ce travail entre Freiburg, Strasbourg et Mulhouse, et Alberto Cattaneo
et Giovanni Felder pour de nombreuses discussions et la concurrence
sympathique.

\section{Op\'{e}rateurs multidiff\'{e}rentiels li\'{e}s \`{a} une sous-va\-ri\'{e}\-t\'{e}}

Soit $X$ une vari\'{e}t\'{e} diff\'{e}rentiable de dimension $n$ et
$i:C\ra X$ une
sous-vari\'{e}t\'{e} ferm\'{e}e de codimension $l\leq n$.
Soit $\Al=\CinfK{X}$, $\Bl=\CinfK{C}$ et
\begin{equation} \label{EqDefIdealAnn}
   \Il=\{f\in\Al~|~f(c)=0~\forall c\in C\}
\end{equation}
l'id\'{e}al annulateur de $C$. En utilisant un voisinage tubulaire autour de $C$
et une partition de l'unit\'{e} on voit que la suite d'alg\`{e}bres commutatives
\begin{equation} \label{EqSuiteIAB}
   \{0\}\lra \Il \lra \Al \stackrel{i^*}{\lra} \Bl \lra \{0\}
\end{equation}
est exacte.
Soit $c\in C$ et
\[
   T_cC^{\ann}=\{\beta\in T_cX^*~|~\beta(v)=0~\forall v\in
    T_cC\}.
\]
Soit $P\in\Gamma(X,\Lambda^2 TX)$ une structure de Poisson.
La sous-vari\'{e}t\'{e} $C$ est dite {\em \cois ~par rapport \`{a} $P$} si
\begin{equation} \label{EqDefPoiCoi}
    P_c(\beta,\gamma)=0~~~\mathrm{quels~que~soient~}c\in C; \beta,\gamma
                              \in T_cC^{\ann}.
\end{equation}
Une condition alg\'{e}brique \'{e}quivalente est
\begin{equation*}
  \Il\hbox{ est une sous-alg\`{e}bre de Poisson de }\Al.
\end{equation*}
De mani\`{e}re g\'{e}n\'{e}rale, on d\'{e}finit
\begin{defi}
\label{compavect}
Un champ de multivecteurs
$P\in\Gamma(X,\Lambda^k TX)$ est dit {\em compatible avec $C$} (ou
{\em adapt\'{e} \`{a} $C$}) si
\begin{equation}
    P_c(\beta_1,\ldots,\beta_k)=0~~~\mathrm{quels~que~soient~}c\in C;
              \beta_1,\ldots,\beta_k
                              \in T_cC^{\ann}.
\end{equation}
\end{defi}
Lorsque $P\in\Al=\Gamma(X,\Lambda^0 TX)$, cette condition \'{e}quivaut
\`{a} $P \in \Il$.

Soit $k$ un entier strictement positif et soient
$\mathcal{M}_1,\ldots,\mathcal{M}_k,\mathcal{M}$ des $\Al$-$\Al$-bimodules.
D\'{e}finissons
\begin{multline*}
   \Dop^k(\mathcal{M}_1,\ldots,\mathcal{M}_k;\mathcal{M})
     =
     \{\phi:\mathcal{M}_1\otimes_{\korps}\cdots\otimes_{\korps}\mathcal{M}_k
            \ra\mathcal{M}~|\cr
\phi\mathrm{~est~}
      k\hbox{-}\mathrm{multidiff\acute{e}rentielle}\}
\end{multline*}
L'expression `$k$-multidiff\'{e}rentiel' correspond \`a la d\'{e}finition alg\'{e}brique de
Grothen\-dieck. Dans la suite, on ne recontrera que des bimodules
plus g\'{e}om\'{e}\-triques comme par exemple $\CinfK{N}$ pour une vari\'{e}t\'{e} $N$ o\`{u}
l'on a l'\'{e}quiva\-lence avec la
d\'{e}finition analytique des op\'{e}rateurs multidiff\'{e}rentiels.
On \'{e}crira $\Dop^k(\mathcal{M}_1;\mathcal{M})$ au lieu de l'expression
ci-dessus lorsque tous les modules $\mathcal{M}_1,$ $\dots,\mathcal{M}_k$
sont \'{e}gaux \`{a} $\mathcal{M}_1$. Pour $k=1$ on \'{e}crit parfois $\Dop$ au lieu de
$\Dop^1$. Pour $k=0$ on pose $\Dop^0(~;\mathcal{M})
=\mathcal{M}$, et on convient que $\Dop^k$ s'annule lorsque $k\leq -1$.
$\Dop^1(\mathcal{M}_1;\mathcal{M})$ est muni d'une structure de
$\Al$-$\Al$-bimodule
avec $(f\phi g)(\eta)=f\phi(g\eta)$ pour $f,g\in\Al,$ $\eta\in \mathcal{M}_1$ et
$\phi\in
\Dop^1(\mathcal{M}_1;\mathcal{M})$.
Soit $\Gs=\oplus_{k\in\entier}\Gs^k$ l'espace des cocha\^{\i}nes
de Gerstenhaber
 o\`{u}
\begin{equation*}
 \Gs^k  =  \Dop^k(\Al;\Al).
\end{equation*}

Par analogie avec les champs de multivecteurs, d\'{e}finissons maintenant
les cocha\^{\i}nes compatibles dont le r\^{o}le sera
tr\`{e}s important pour la suite de notre \'{e}tude~:
\begin{defi}
\label{compacoch}
Le sous-espace $\GsI=\oplus_{k\in\entier}\GsI^k$ de $\Gs$ des
{\sl cocha\^{\i}nes compatibles} est d\'{e}fini comme suit~:
\begin{align*}
 \GsI^k &=
  \{\phi\in \Gs^k|\phi(f_1,\ldots,f_{k-1},g)\in\Il,\cr
&\hskip3cm  \forall f_1,\ldots,f_{k-1}\in\Al,  g\in\Il\}~(\hbox{pour }k\geq 1), \cr
 \GsI^0 & =  \Il,\cr
 \GsI^k & =  \{0\}\hbox{ pour }k\leq -1.
\end{align*}
\end{defi}

Pour la suite de notre travail, nous d\'{e}finissons
$\GsP=\oplus_{k\in\entier}\GsP^k$
o\`{u}
\begin{align*}
\GsP^k & = \Dop^{k-1}\big(\Al;\Dop^1(\Il;\Bl)\big)
               ~(\hbox{pour }k\geq 1), \cr
\GsP^0 & =  \Al/\Il = \Bl, \cr
\GsP^k & =  \{0\}\hbox{ pour }k\leq -1.
\end{align*}
Remarquons que $\Dop^1(\Al;\Bl)$ est l'espace des op\'{e}rateurs diff\'{e}rentiels
le long de l'application $i:C\ra M$, voir \cite{Bor03} pour une d\'{e}finition
analytique, et $\Dop^1(\Il;\Bl)$ est la restriction de $\Dop^1(\Al;\Bl)$
\`{a} $\Il$. En outre, l'espace des op\'{e}rateurs $\Dop^1(\Bl,\Bl)$ s'injecte
dans $\Dop^1(\Al;\Bl)$ par l'application $\mathsf{D}\mapsto \big(f\mapsto
\mathsf{D}(i^*f)\big)$. Un petit calcul en coordonn\'{e}es de sous-vari\'{e}t\'{e}s
montre que la suite de $\Al$-$\Al$-bimodules
\[
  \{0\} \leftarrow \Dop^1(\Il;\Bl) \leftarrow \Dop^1(\Al;\Bl) \leftarrow
         \Dop^1(\Bl,\Bl) \leftarrow \{0\}
\]
est exacte (voir aussi le lemme \ref{LemSuiteExacteDiffOp}).
En utilisant que $\Dop^{k-1}\big(\Al;\Dop^1(\Al,\Bl)\big)$
est \'{e}gal \`{a} $\Dop^k(\Al;\Bl)$ (on travaille avec des morphismes de
$\korps$-modules filtr\'{e}s) on d\'{e}finit
$\Dop^{k-1}\big(\Al;\Dop^1(\Il;\Bl)\big)$ par la restriction du $k$-i\`eme
argument \`{a} $\Il$.

D'autre part, on peut d\'{e}finir des analogues des espaces
$\GsI$ et $\GsP$ dans l'espace des
champs de multivecteurs $\gs=\oplus_{k\in\entier}\gs^k$ o\`{u}
\begin{align*}
   \gs^k & =  \Gamma(X,\Lambda^k TX) ~(\hbox{pour }k\geq 0), \cr
   \gs^k & =  \{0\}\hbox{ pour }k\leq -1.
\end{align*}
L'analogue de $\GsI $ est
$\gsI=\oplus_{k\in\entier}\gsI^k$ le sous-espace de $\gs$
o\`{u}
\begin{align*}
   \gsI^k & =  \{P\in\Gamma(X,\Lambda^k TX)|
                  P_c(\beta_1,\ldots,\beta_k)=0,\cr
&\hskip3cm \forall c\in C; \beta_1,\ldots,\beta_k \in T_cC^{\ann}\}
 ~(\hbox{pour }k\geq 1), \cr
   \gsI^0 & =  \Il \cr
   \gsI^k & =  \{0\}\hbox{ pour }k\leq -1.
\end{align*}

\noindent L'analogue de $\GsP$ est l'espace $\gsP=\oplus_{k\in\entier}\gsP^k$ o\`{u}
\begin{align*}
  \gsP^k  & =  \Gamma\big(C,\Lambda^k (T_CX/TC)\big)
                    ~(\hbox{pour }k\geq 1),\cr
  \gsP^0 & =  \Bl=\Al/\Il \cr
   \gsP^k & = \{0\}\hbox{ pour }k\leq -1.
\end{align*}

\bigskip

\section{Th\'{e}or\`{e}mes de Hochschild-Kostant-Rosenberg}
  \label{SecHKR}



Dans cette section nous allons montrer un analogue du th\'{e}or\`{e}me
de Hochschild-Kostant-Rosenberg dans le cas o\`{u}
$X=\real^n$ et $C=\real^{n-l}$. Dans les deux premi\`{e}res
sous-sections
nous rappellerons les propri\'{e}t\'{e}s alg\'{e}briques de
$\Gs$, $\GsI$, et $\GsP$, puis de $\gs$, $\gsI$ et
$\gsP$. Puis dans la troisi\`{e}me nous donnerons des homotopies
explicites
entres les r\'{e}solutions bar et de Koszul qui nous permettrons de
prouver
notre r\'{e}sultat principal dans la quatri\`{e}me sous-section.

\subsection{Propri\'{e}t\'{e}s alg\'{e}briques de $\Gs$, $\GsI$ et $\GsP$}

\noindent
Rappelons quelques op\'{e}rations d\'{e}finies sur l'espace des cocha\^{\i}nes de
Gerstenhaber~: pour
$k,l\geq 1$ entiers, $\phi\in \Gs^k$, $\psi\in \Gs^l$ et $1\leq i\leq
k$ on pose
\begin{multline} \label{EqDefCirci}
   (\phi\circ_i \psi)\big(f_1,\ldots,f_{k+l-1}\big)\cr
     = \phi\big(f_1,\ldots,f_{i-1},\psi(f_i,\ldots,f_{l+i-1}),f_{l+i},
       \ldots,f_{k+l-1}\big),
\end{multline}
et l'on d\'{e}finit
\begin{equation} \label{EqDefMultGerstDiff}
 \phi\circ_G \psi = \sum_{i=1}^k(-1)^{(i-1)(l-1)}\phi\circ_i \psi,
\end{equation}
(on pose $\phi\circ_G\psi=0$ lorsque $\phi\in\Gs^k$, $k\leq 0$),
et le {\em crochet de Gerstenhaber}~:
\begin{equation} \label{EqDefCrochGerstDiff}
  [\phi,\psi]_G=\phi\circ_G\psi -(-1)^{(k-1)(l-1)}\psi\circ_G\phi.
\end{equation}
Gerstenhaber a montr\'{e} \cite{Ger63} que $(\Gs[1],[~,~]_G)$ est une alg\`{e}bre de
Lie gradu\'{e}e (o\`{u} on a utilis\'{e} le {\em d\'{e}calage}
$(\hs[j])^k:=\hs^{j+k}$ pour $j\in\mathbb Z$ et un $\korps$-espace gradu\'{e}
$\hs:=\oplus_{k\in\mathbb Z}\hs^k$). Soit $\mu$
la multiplication point par point dans $\Al$. On d\'{e}finit alors
l'op\'{e}rateur cobord de Hochschild par~:
\begin{equation} \label{EqDefCobH}
    \mathsf{b}\phi= -[\phi,\mu]_G
\end{equation}
Enfin, on d\'{e}finit la multiplication $\cup$ par~:
\begin{equation} \label{EqDefMultCup}
    (\phi\cup\psi)\big(f_1,\ldots,f_{k+l}\big)=\mu\big(\phi(f_1,\ldots,f_k),
                    \psi(f_{k+1},\ldots,f_{k+l})\big)
\end{equation}
et il est clair que $(\Gs,\cup)$ est une alg\`{e}bre associative gradu\'{e}e.
\begin{prop} \label{PGsI}
 L'espace $\GsI$ a les propri\'{e}t\'{e}s suivantes:
 \begin{enumerate}
  \item $\mu\in\GsI^2$.
  \item $\GsI$ est un id\'{e}al \`{a} gauche de $(\Gs,\cup)$.
  \item $\GsI[1]$
  et une sous-alg\`{e}bre de Lie gradu\'{e}e de $(\Gs[1],[~,~]_G)$.
  \item $(\GsI,\mathsf{b})$ est un sous-complex de $(\Gs,\mathsf{b})$.
 \end{enumerate}
\end{prop}
\begin{prooofe}
\begin{enumerate}
\item C'est \'{e}vident car $\Il$ est un id\'{e}al de $\Al$.
\item Soit $\phi\in\Gs^k$, $\psi\in \GsI^l$ et $f_{k+l}\in\Il$. On a
 $\psi(f_{k+1},\ldots,f_{k+l})\in\Il$ donc
 $(\phi\cup\psi)(f_1,\ldots,f_{k+l})\in\Il$ car $\Il$ est un id\'{e}al
 de $\Al$.
\item Soient $\phi\in\GsI^k$ et $\psi\in\GsI^l$. Regardons
 $\phi\circ_i\psi$ (si $k\leq 0$ il n'y a rien \`{a} montrer) pour $1\leq
 i\leq k$. Soit $f_{k+l-1}\in\Il$. Si $1\leq i\leq k-1$, alors
 $f_{k+l-1}$ est le dernier argument de $\phi$, donc le membre
 droit de (\ref{EqDefCirci}) appartient \`{a} $\Il$ par d\'{e}finition de $\phi$,
 et $\phi\circ_i\psi\in\GsI^{k+l-1}$. Si $i=k$, alors $f_{k+l-1}$
 est le dernier argument de $\psi$ et la valeur de $\psi$ (qui appartient
 \`{a} $\Il$ par d\'{e}finition de $\psi$) est le dernier argument de $\phi$.
 Par d\'{e}finition de $\phi$, il vient que sa valeur appartient \`{a} $\Il$,
 donc $\phi\circ_k \psi\in\GsI^{k+l-1}$.
\item C'est une cons\'{e}quence de 1., 3. et  de la
 d\'{e}finition
de $\mathsf{b}$
 (\ref{EqDefCobH}).\hfill  $\Box$
\end{enumerate}
\end{prooofe}

L'espace $\GsP$ est muni de l'op\'{e}rateur de Hochschild $\widetilde{\mathsf{b}}$
usuel~: soient $\phi\in\GsP^k$, $f_1,\ldots,f_{k}\in\Al$ et $g\in\Il$,
alors
\begin{eqnarray}\label{EqDefCobHP}
 (\widetilde{\mathsf{b}}\phi)(f_1,\ldots,f_k)(g) & = &
      f_1~\phi(f_2,\ldots,f_k)(g)  \cr
          &  & +\sum_{r=1}^{k-1}(-1)^{r}
             \phi(f_1,\ldots,f_rf_{r+1},\ldots,f_{k})(g) \cr
          &&  +(-1)^k\phi(f_1,\ldots,f_{k-1})(f_kg).
\end{eqnarray}
Pour un entier $k$ consid\'{e}rons les projections canoniques
$\Xi^k:\Gs^k\ra\GsP^k$ suivantes o\`{u} $f_1,\ldots,f_{k-1}\in\Al$ et $g\in\Il$~:
\begin{equation*}
  \big(\Xi^k (\phi)\big)\big(f_1,\ldots,f_{k-1}\big)(g) =
             i^*\big(\phi(f_1,\ldots,f_{k-1},g)\big)
\end{equation*}
pour $k\geq 1$, $\Xi^0=i^*$ et $\Xi^k=0$ quel que soit $k\leq -1$.
\begin{prop}
 On a les propri\'{e}t\'{e}s suivantes~:
 \begin{enumerate}
  \item Le diagramme suivant est une suite exacte de complexes:
  \begin{equation*}
  \{0\}\lra (\GsI,\mathsf{b}) \lra (\Gs,\mathsf{b}) \stackrel{\Xi}{\lra}
    (\GsP,\widetilde{\mathsf{b}}) \lra \{0\}.
  \end{equation*}
  En particulier, $\GsP\cong \Gs/\GsI$.
  \item $\GsP$ est un module \`{a} gauche gradu\'{e} de $(\Gs,\cup)$.
  \item $\GsP[1]$ est un module de Lie gradu\'{e} de $(\GsI[1],[~,~]_G)$.
 \end{enumerate}
\end{prop}
\begin{prooofe}
\begin{enumerate}
\item Il est clair que le noyau de $\Xi$ est \'{e}gal \`{a} $\GsI$ et que
$\Xi$ est surjective. Montrons que $\Xi$
est un morphisme de complexes~: soient
 $\phi\in\Gs^k$, $f_1,\dots,f_k \in \Al$, $g \in \Il$, on a
 \begin{align*}
   \big(\Xi^{k+1}(\mathsf{b}\phi)\big)\big(f_1,\ldots,f_{k}\big)(g) = &
     i^*\big((\mathsf{b}\phi)(f_1,\ldots,f_{k},g)\big) \cr
           &\hskip-3cm = i^*\big(f_1\phi(f_2,\ldots,f_k,g)\big)
  +(-1)^{k+1} i^*\big(\phi(f_1,\ldots,f_{k})g\big)\cr
           &\hskip-2.6cm +\sum_{r=1}^{k-1}(-1)^{r}
                i^*\big(\phi(f_1,\ldots,f_rf_{r+1},
                               \ldots,f_{k},g)\big) \cr
  &\hskip-2.6cm+(-1)^ki^*\big(\phi(f_1,\ldots,f_{k-1},f_kg)\big) \cr
          &\hskip-3cm =
     i^*f_1~\big(\Xi^k\phi\big)(f_2,\ldots,f_k)(g)
+ 0 \hbox{ (car }i^*g=0)\cr
          &\hskip-2.6cm +\sum_{r=1}^{k-1}(-1)^{r}
                \big(\Xi^k\phi\big)(f_1,\ldots,f_rf_{r+1},
                               \ldots,f_{k})(g)  \cr
 &\hskip-2.6cm+ (-1)^k\big(\Xi^k\phi\big)(f_1,\ldots,f_{k-1})(f_kg) \cr
          &\hskip-3cm =  (\widetilde{\mathsf{b}}\Xi^k\phi)(f_1,\ldots,f_k)(g).
 \end{align*}
\item C'est une cons\'{e}quence du fait
que $\GsI$ est un id\'{e}al \`{a} gauche de $\Gs$ (Proposition
 \ref{PGsI}, 2.).
\item Puisque $\GsI[1]$ est une sous-alg\`{e}bre de Lie gradu\'{e}e de $\Gs[1]$
 d'apr\`{e}s la Proposition \ref{PGsI} 3., il vient que $\Gs[1]$ et $\GsI[1]$
 sont des modules de Lie gradu\'{e}s de $\GsI[1]$, donc il en est le m\^{e}me
 pour leur quotient $\GsP[1]$.\hfill  $\Box$
\end{enumerate}
\end{prooofe}

\subsection{Propri\'{e}t\'{e}s alg\'{e}briques de
 $\gs$, $\gsI$ et $\gsP$}
   \label{SubSecPropalggs}

\noindent
\'Etudions maintenant les espaces $\gs,\gsI$ et $\gsP$. Rappelons la
d\'{e}finition du
crochet de Schouten $[-,-]_S$ sur $\Gs$~:
soient $f,g\in\Al$, $k,l\in\nat$, $X=X_1\wedge\cdots\wedge X_k\in\gs^k$ et
$Y=Y_1\wedge\cdots\wedge Y_l\in \gs^l$
($X_1,\ldots,X_k,
Y_1,\ldots,Y_l\in\Gamma(X,TX)$), alors le crochet est d\'{e}fini par~:
\begin{eqnarray}
 [fX,gY]_S   \!\!  \!&=&\!\!  fg\sum_{i=1}^k\sum_{j=1}^l(-1)^{i+j}[X_i,Y_j]
       \wedge
        X_1\wedge\cdots\wedge X_{i-1}\wedge X_{i+1}\wedge\cdots\wedge
        X_k \nonumber\\
      & &\!\!\hskip4.5cm\wedge Y_1\wedge\cdots\wedge Y_{j-1}\wedge
         Y_{j+1}\wedge\cdots\wedge Y_l \nonumber\\
      & &\!\!+ f \sum_{i=1}^k (-1)^{k-i}X_i(g)
        X_1\wedge\cdots\wedge X_{i-1}\wedge X_{i+1}\wedge\cdots\wedge
        X_k\wedge Y \nonumber\\
      &  &\!\!+ gX\wedge \sum_{j=1}^l (-1)^{j}Y_j(g)
        Y_1\wedge\cdots\wedge Y_{j-1}\wedge X_{j+1}\wedge\cdots\wedge
        Y_l. \label{EqDefSchouten}
\end{eqnarray}
Il est bien connu que ce crochet ne d\'{e}pend pas de la d\'{e}composition de $X$ et de $Y$
 en produit de champs de vecteurs. En outre $(\gs[1],[-,-]_S)$ est une alg\`{e}bre de
 Lie gradu\'{e}e~; de plus, $(\gs,\wedge)$ est une alg\`{e}bre
 associative commutative gradu\'{e}e. Enfin, pour tous
$X\in\gs^k,Y\in\gs^l,Z\in\gs^l$, on a la r\`{e}gle de Leibniz gradu\'{e}e~:
 \begin{equation*}
    [X,Y\wedge Z]_S = [X,Y]_S\wedge Z + (-1)^{(k-1)l}Y\wedge [X,Z]_S
 \end{equation*}
et donc $(\gs,[-,-]_S,\wedge)$ est une alg\`{e}bre de Gerstenhaber, voir
\'{e}galement (\ref{EqRegleLeibniz}).

\smallskip

Le produit int\'{e}rieur $i$ peut s'\'{e}tendre en une action
(toujours not\'{e}e $i$) de $(\gs,\wedge)$ sur
l'espace $\Gamma(X,\Lambda
 T^*X)$~:
 soit $x=x_1\wedge\cdots\wedge x_k\in \gs^k$,
 $f \in \Al$ et $\alpha\in\Gamma(X,\Lambda
 T^*X)$, alors
 \begin{equation*}
    i(x)\alpha =i(x_1)\cdots i(x_k)\alpha
\hbox{ et }i(f)\alpha = f \alpha.
 \end{equation*}
De m\^{e}me, la d\'{e}riv\'{e}e de Lie des champs de vecteurs s'\'{e}tend en
 une action $L$
de $(\gs[1],[-,-]_S)$ sur
 $\Gamma(X,\Lambda T^*X)$ d\'{e}finie par
 \begin{equation} \label{EqDefLSch}
    L(x)\alpha = [i(x),\drm]\alpha.
 \end{equation}
Par r\'{e}currence sur le degr\'{e} de $y$, on montre ais\'{e}ment que
 \begin{equation} \label{EqLiSch}
     [L(x),i(y)]=i([x,y]),~\forall x,y\in\gs
 \end{equation}
Enfin, puisque $\drm^2=0$ on a
 \begin{equation*}
     [L(x),\drm]=0\hbox{ et }    [L(x),L(y)]=L([x,y]),~\forall x,y \in \gs
 \end{equation*}

\bigskip

Consid\'{e}rons maintenant les applications  $\Psi^0=i^*:\Al\ra
\Al/\Il=\Bl$
et pour $k\geq 1$, $\Psi^k:\gs^k \lra  \gsP^k$
d\'{e}finie par
$$   x_1\wedge\cdots\wedge x_k  \mapsto
                           \big(c\mapsto
                           (x_{\stackrel{1}{c}}\mathrm{~mod~}T_cC)
                           \wedge\cdots\wedge
                           (x_{\stackrel{k}{c}}\mathrm{~mod~}T_cC)\big). $$
\begin{prop}
 La suite d'espaces vectoriels sur $\real$
 \begin{equation}
\{0\}\lra \gsI^k \lra \gs^k \stackrel{\Psi^k}{\lra} \gsP^k \lra \{0\}
 \end{equation}
 est exacte. En particulier $\gsP^k\cong \gs^k/\gsI^k$.
\end{prop}
\begin{prooof}
 Il est clair que $\Psi^k$ est bien d\'{e}finie. Le cas
 $k=0$ est une cons\'{e}quence
 de la suite exacte (\ref{EqSuiteIAB}).\\
 Montrons la surjectivit\'{e} de $\Psi^k$: soit
 $\widetilde{Y}\in\gsP^k$. Soit $E\subset
 T_CX$ un sous-fibr\'{e} tel que $T_CX=TC\oplus E$, soit $C\subset V\subset
 E$ un voisinage ouvert de la section nulle de $E$, soit $C\subset
 U\subset X$ un voisinage ouvert de $C$ dans $X$ et soit $\Phi:V\ra U$ un
 diff\'{e}omorphisme tel que $\Phi|_C$ est l'application identique (le tout
 est dit un voisinage tubulaire de $C$). $E$ est visiblement isomorphe
 \`{a} $T_CX/TC$. Pour $c\in C$ on rappelle le rel\`{e}vement vertical de $Y\in E_c$ \`{a}
 $y\in E_c$: on a $Y^{\rlvt_y}=\frac{\drm}{\drm t}(y+tY)|_{t=0}$,
 donc $Y^{\rlvt_y}\in T_{y}E$. Ceci induit une injection
 $(~)^{\rlvt_{\Phi(y)}}:\Lambda^k (T_cX/T_cC)\ra \Lambda^k T_{\Phi(y)}U$,
 donc un
 rel\`{e}vement vertical des sections $\widetilde{Y}\mapsto
 Y'=(\widetilde{Y})^\rlvt\in \Gamma(X,\Lambda^k TU)$ d\'{e}fini par
 ${Y'}_{\Phi(y)}=(\widetilde{Y}_c)^{\rlvt_{y}}$. Soit
 $(\psi_U,1-\psi_U)$ une partition de l'unit\'{e} subordonn\'{e}e
 au r\'{e}couvrement
 ouvert $(U,X\setminus C)$ de $X$. Alors le
 champ de vecteurs $Y$ d\'{e}fini par $Y=\psi_U Y'$ sur $U$ et
 $Y=0$ sur $X\setminus U$ est un \'{e}l\'{e}ment de $\gs^k$ tel que
 $\Psi^k Y=\widetilde{Y}$.\\
 Finalement, montrons que
$\gsI^k$ est \'{e}gal au noyau de $\Psi^k$:
 pour $c\in C$ soit $\Psi^k_c:\Lambda^kT_cX\ra \Lambda^k(T_cX/T_cC)$.
 Puisque $T_cC^{\ann}\cong (T_cX/T_cC)^*$ alors
$$    \Lambda^k (T_cC^{\ann})\cong \big(\Lambda^k (T_cX/T_cC)\big)^*.
$$
 Soit $Y\in\gs^k$. Alors $Y\in\mathrm{Ker}\Psi^k$ si
et seulement si $\Psi^k_c(Y_c)=0$ quel
 que soit $c\in C$. Il s'ensuit~:
 \[
   \begin{array}{ccll}
    \Psi^k_c(Y_c)=0 &\Longleftrightarrow &
                     \xi_c\big(\Psi^k_c(Y_c)\big)=0 &~~\forall\xi_c\in
                          \Lambda^k (T_cX/T_cC)^{*} \\
                       &\Longleftrightarrow &
                     \xi_c(Y_c)=0 &~~\forall\xi_c\in
                          \Lambda^k (T_cC)^{\ann} \\
                         &\Longleftrightarrow &
                     Y_c(\xi_c)=0 &~~\forall\xi_c\in
                          \Lambda^k (T_cC)^{\ann},
  \end{array}
 \]
 et donc $\gsI^k=\mathrm{Ker}\Psi^k$.
\end{prooof}

\noindent On obtient alors une autre caract\'{e}risation de $\gsI$:
\begin{prop} \label{PgsIIl}
 Soit $k$ un entier strictement positif.
 Un \'{e}l\'{e}ment $X\in\gs^k$ appartient \`{a} $\gsI^k$ si et seulement si
 $$i(X)(\drm g_1\wedge\cdots\wedge \drm g_k)\in\Il\mathrm{~quels~que~soient~}
                   g_1,\cdots,g_k\in\Il
 $$
\end{prop}
\begin{prooof}
 Il est clair que $\drm_c g\in T_cC^{\mathrm{ann}}$ quel que soit $c\in C$ et
 quel que soit $g\in\Il$, donc la condition ci-dessus est n\'{e}cessaire.\\
 D'autre part, soit $\alpha\in T_cC^{\ann}$ et soit
 $\big(U,(x^1,\ldots,x^{n-l},y^1,\ldots,y^l)\big)$ une carte de
 sous-vari\'{e}t\'{e} de $C$ autour de $c$ (c'est-\`{a}-dire
 $U\cap C=\{u\in U~|~y^1(u)=0,\ldots,y^l(u)=0\}$). Alors on trouve
 $\alpha_1,\ldots,\alpha_l\in\real$ tels que
 $\alpha=\sum_{i=1}^l\alpha_i \drm y^i$. Soit $(\psi_U,1-\psi_U)$ une
 partition de l'unit\'{e} subordonn\'{e}e au r\'{e}couvrement ouvert
 $(U,X\setminus \{c\})$
 de $X$. On d\'{e}finit $g= \psi_U\sum_{i=1}^l\alpha_i y^i$ sur $U$ et $g=0$
 sur $X\setminus U$. Visiblement $g\in\Il$ et $\drm_c g=\alpha$. Alors tout
 \'{e}l\'{e}ment de $T_cC^{\ann}$ se repr\'{e}sente comme $\drm_c g$ pour un
 $g\in\Il$, d'o\`{u} la suffisance de la condition.
\end{prooof}

\begin{prop}\label{PropPropgsI}
 L'espace $\gsI$ a les propri\'{e}t\'{e}s suivantes:
 \begin{enumerate}
  \item Toute structure de Poisson $P$ sur $X$ compatible avec $C$ est dans $\gsI^2$.
  \item $(\gsI,\wedge)$ est un id\'{e}al de $(\gs,\wedge)$.
  \item $(\gsI[1]\!,\![-,-]_S)$ est une sous-alg\`{e}bre de Lie gradu\'{e}e
       de $(\gs[1]\!,\![-,-]_S)$.
 \end{enumerate}
\end{prop}
\begin{prooofe}
\begin{enumerate}
\item Ceci est une cons\'{e}quence directe de la d\'{e}finition
 (\ref{EqDefPoiCoi}).
\item $\Psi_c=\sum_{k=0}^n\Psi^k_c:\Lambda T_cX\ra \Lambda (T_cX/T_cC)$ est un
   homomorphisme d'alg\`{e}\-bres de Grassmann, alors $\Psi:\gs\ra\gsP$ aussi,
   donc le noyau de $\Psi$, alors $\gsI$, est un id\'{e}al par rapport \`{a} la
   multiplication ext\'{e}rieure $\wedge$.
\item Soient $x\in\gsI^k$ et $y\in\gsI^l$. Si $k=l=0$ alors
  $[x,y]_S=0\in\gsI^{-1}$. Si $k=1$ et $l=0$, alors $y\in\Il$ et
  $[x,y]_S=x(\drm y)\in\Il=\gsI^{0}$ d'apr\`{e}s la proposition \ref{PgsIIl}.
  On peut supposer que $k+l\geq 2$. Soient $g_1,\ldots,g_{k+l-1}\in\Il$
  et $\gamma=\drm g_1\wedge\cdots\wedge \drm g_{k+l-1}$.
  Puisque $x\wedge y\in \gsI^{k+l}$, alors $i(x)i(y)\gamma=0$. De plus,
  $\drm\gamma=0$. Il s'ensuit, d'apr\`{e}s l'\'{e}quation (\ref{EqLiSch}), que
  $i([x,y])\gamma=[L(x),i(y)]\gamma = i(x)\drm i(y)\gamma
  +(-1)^{kl}i(y)\drm i(x)\gamma$.  La $k-1$ forme $i(y)\gamma$ est une somme
  finie
  d'expressions de la forme~: $\pm(i(y)\gamma')\gamma''$ o\`{u} $\gamma'$
  est une $l$-forme constitu\'{e}e de $l$ \'{e}l\'{e}ments parmi
  $\drm g_1,\ldots,\drm g_{k+l-1}$ et $\gamma''$ est le reste. D'apr\`{e}s la
  proposition \ref{PgsIIl} la fonction $g'=\pm i(y)\gamma'$ est dans
  $\Il$,
alors
  $i(x)\drm i(y)\gamma=i(x)(\drm g'\wedge \gamma'')$ car $\drm\gamma''=0$ et
  la
  $k$-forme $\drm g'\wedge \gamma''$ est le produit ext\'{e}rieur
de $k$ diff\'{e}rentielles
  d'\'{e}l\'{e}ments de l'id\'{e}al, donc $i(x)
  (\drm g'\wedge \gamma'')\in\Il$. Le terme
  $(-1)^{kl}i(y)\drm i(x)\gamma$ appartient \'{e}galement \`{a} $\Il$ par un
  raisonnement enti\`{e}rement analogue. Alors
  $[x,y]_S\in\gsI^{k+l-1}$.\hfill  $\Box$
\end{enumerate}
\end{prooofe}
\begin{prop}
 L'espace $\gsP$ a les propri\'{e}t\'{e}s suivantes~:
 \begin{enumerate}
  \item $(\gsP,\wedge)$ est une alg\`{e}bre commutative associative gradu\'{e}e.
  \item $(\gsP,\wedge)$ est un module \`{a} gauche gradu\'{e} de $(\gs,\wedge)$.
  \item $\gsP[1]$ est un module d'alg\`{e}bre de Lie gradu\'{e}e pour
        $(\gsI[1],[-,-]_S)$. De plus, on a
        \[
         X.(\alpha\wedge\beta)=(X.\alpha)\wedge\beta
        +(-1)^{(k-1)l}\alpha\wedge (X.\beta)
        \]
        quels que soient $X\in\gsI^k,\alpha\in\gsP^l,\beta\in\gsP$.
 \end{enumerate}
\end{prop}
\begin{prooofe}
\begin{enumerate}
\item Ceci est \'{e}vident.
\item Puisque $\gsI$ est un id\'{e}al de $(\gs,\wedge)$ et $\gsP\cong \gs/\gsI$,
   l'\'{e}nonc\'{e} est \'{e}vident.
\item $\gs[1]$ et $\gsI[1]$ sont \'{e}videmment des modules de
    $(\gsI,[~,~]_S)$, et il en est de m\^{e}me pour le quotient
    $\gsP[1]$.\hfill  $\Box$
\end{enumerate}
\end{prooofe}
\begin{prop} \label{PgspEmbed}
 Si $C$ a un voisinage tubulaire global dans $M$ il existe une
 injection $\gsP\ra \gs$ avec
 \begin{enumerate}
  \item $\gs=\gsI\oplus \gsP$ (somme directe d'espaces vectoriels).
  \item $(\gsP,\wedge)$ est une sous-alg\`{e}bre gradu\'{e}e de $(\gs,\wedge)$.
  \item $\gsP[1]$ muni du crochet de Schouten est une sous-alg\`{e}bre
    de Lie gradu\'{e}e ab\'{e}lienne de $(\gs[1],[~,~]_S)$.
 \end{enumerate}
\end{prop}
 \begin{prooofe}
  $M$ est diff\'{e}omorphe \`{a} l'espace total d'un fibr\'{e} vectoriel $\tau:E\ra C$
  et $\gsP\cong \Ginf(C,\Lambda E)$. En utilisant le rel\`{e}vement vertical
  $E_c\ra T_\epsilon E$ pour $\epsilon\in E$ et $\tau(\epsilon)=c$ donn\'{e}
  par $f\mapsto \frac{\partial (f+t\epsilon)}{\partial t}|_{t=0}$ on
  obtient l'injection avec toutes les propri\'{e}t\'{e}s \'{e}nonc\'{e}es.
 \end{prooofe}
\begin{rem}
 Soit $P\in\gsI^2$ une structure de Poisson compatible avec $C$. Alors
 \`{a} l'aide de l'action de $P$ l'espace $\gsP$ devient une alg\`{e}bre
 associative commutative diff\'{e}rentielle gradu\'{e}e. La cohomologie de
 $\gsP$ est dite la {\em cohomologie BRST de $C$ par rapport \`{a} $P$.}
 Cette cohomologie est munie d'une structure d'alg\`{e}bre de Poisson.
\end{rem}

\subsection{Simplification du complexe bar}
\noindent
Dans ce paragraphe $X=\real^n$.

\bigskip

Commen\c{c}ons par rappeler la r\'{e}solution `topologique'
bar de l'alg\`{e}bre $\Al=\CinfK{\real^n}$:
Soit $\Al^e=\CinfK{\real^{2n}}$
et pour tout entier positif $k$
\begin{equation*}
  CH^k= \CinfK{\real^{(k+2)n}}.
\end{equation*}
Nous noterons $(a,x_1,\ldots,x_k,b)$ (o\`{u} $a,x_1,\ldots,x_k,b\in\real^n$)
 pour un point de $\real^{(k+2)n}$.
L'espace $CH^k$ est un $\Al^e$-module~:
\begin{equation*}
 \Al^e\!\times CH^k\! \ra CH^k\!:(f,\Phi)\mapsto\! \big((a,x_1,\ldots,x_k,b)
                                   \mapsto
                                   f(a,b)\Phi(a,x_1,\ldots,x_k,b)\big)
\end{equation*}
Pour $k\geq 1$, rappelons que l'op\'{e}rateur bord de Hochschild
$\dH^k:CH^k\ra CH^{k-1}$ est d\'{e}fini par
\begin{eqnarray}\label{EqDefBordHochBar}
  (\dH^k \Phi)(a,x_1,\ldots,x_{k-1},b)  & = &
                      \Phi(a,a,x_1,\ldots,x_{k-1},b)\nonumber \\
             & &
             +\sum_{r=1}^{k-1}(-1)^r \Phi(a,x_1,\ldots,x_r,x_r,
                                            \ldots,x_{k-1},b)
                          \nonumber \\
             & &
             (-1)^k \Phi(a,x_1,\ldots,x_{k-1},b,b)
\end{eqnarray}
Il est clair que $\dH^k$ est un morphisme de $\Al^e$-modules et
$\dH^k\dH^{k+1}=0$. L'augmentation $\epsilon:CH^0=\Al^e\ra \Al$ est d\'{e}finie
par le morphisme de $\Al^e$-modules suivant
\begin{equation} \label{EqDefEps}
   (\epsilon \Phi)(a)= \Phi(a,a)~~~\forall a\in\real^n.
\end{equation}
Il est bien connu que le complexe bar
\begin{equation} \label{EqCompBar}
  \{0\}\lla \Al \stackrel{\epsilon}{\lla} CH^0 \stackrel{\dH^1}{\lla}
           CH^1 \stackrel{\dH^2}{\lla} CH^2 \stackrel{\dH^3}{\lla}\cdots
              \stackrel{\dH^k}{\lla} CH^k \stackrel{\dH^{k+1}}{\lla}
              \cdots
\end{equation}
est acyclique~: en fait, soit $h_H^{-1}:\Al\ra CH^0$ le prolongement
$\korps$-lin\'{e}aire
\begin{equation*}
   (h_H^{-1}f)(a,b):=f(a)
\end{equation*}
et, pour $k\geq 0$, $h_H^k:CH^k\ra CH^{k+1}$ l'application $\korps$-lin\'{e}aire
\begin{equation*}
   (h_H^k \Phi)(a,x_1,\ldots,x_{k+1},b)=(-1)^{k+1}\Phi(a,x_1,\ldots,x_k,x_{k+1}).
\end{equation*}
En \'{e}crivant $\id_H^{-1}$ pour l'application identique $\Al\ra\Al$
et $\id_H^k$ pour l'application identique $CH^k\ra CH^k$ on
montre que
\begin{eqnarray*}
       \epsilon h_H^{-1} & = & \id_H^{-1}  \\
       h_H^{-1}\epsilon + \dH^1 h_H^0 & = & \id_H^0 \\
       h_H^{k-1}\dH^k + \dH^{k+1} h_H^k & = & \id_H^k~~~\forall k\geq 1
\end{eqnarray*}
ce qui entra\^{\i}ne l'acyclicit\'{e} du complexe bar (\ref{EqCompBar}).

\bigskip

Nous allons maintenant
d\'{e}finir un autre complexe (de Koszul) acyclique pour $\Al$ en tant que
$\Al^e$-module~: soit $E=\korps^n$
\begin{equation*}
   CK^k = \Al^e\otimes_{\korps}\Lambda^k E^*
   ~~\mathrm{quel~que~soit~}k\in\entier.
\end{equation*}
\'Evidemment, chaque $CK^k$ est un $\Al^e$-module libre. Soit
$\xi:\real^{2n}\ra E$ d\'{e}fini par
\begin{equation*}
    \xi(a,b)=a-b
\end{equation*}
Pour tout entier $k$
strictement positif, on d\'{e}finit
l'op\'{e}rateur bord de Koszul $\dK^k:CK^k\ra CK^{k-1}$
par
\begin{equation*}
    \dK^k\omega = i(\xi)\omega~~\mathrm{quel~que~soit~}\omega\in CK^k.
\end{equation*}
Autrement dit, pour $e_1,\ldots,e_k\in E$ et $a,b\in\real^n$ on a
\begin{equation*}
    (\dK^k\omega)(a,b)\big(e_2,\ldots,e_k\big)
             =\omega(a,b)\big(a-b,e_2,\ldots,e_k\big).
\end{equation*}
Il est clair que les $\dK^k$ sont des morphismes de $\Al^e$-modules
et que $\dK^k\dK^{k+1}=0$ quel que soit l'entier strictement positif $k$.
Soit $\epsilon:CK^0=\Al^e\ra\Al$ l'augmentation d\'{e}finie comme dans
(\ref{EqDefEps}). Il en r\'{e}sulte le complexe de Koszul~:
\begin{equation} \label{EqCompKos}
  \{0\}\lla \Al \stackrel{\epsilon}{\lla} CK^0 \stackrel{\dK^1}{\lla}
           CK^1 \stackrel{\dK^2}{\lla} CK^2 \stackrel{\dK^3}{\lla}\cdots
              \stackrel{\dK^k}{\lla} CK^k \stackrel{\dK^{k+1}}{\lla}
              \cdots
\end{equation}
Ce complexe est acyclique~: en effet, soit
$h_K^{-1}=h_H^{-1}:\Al\ra CK^0=\Al^e=CH^0$
le prolongement $\korps$-lin\'{e}aire
\begin{equation*}
   (h_K^{-1}f)(a,b)=f(a).
\end{equation*}
Soit $e_1,\ldots,e_n$ une base de $E$ et $e^1,\ldots,e^n$ la base duale de
$E^*$. Pour $k\geq 0$ soit $h_K^k:CK^k\ra CK^{k+1}$ l'application
$\korps$-lin\'{e}aire
\begin{equation*}
   (h_K^k\omega)(a,b)
     =-\sum_{j=1}^n e^j \wedge
          \int_0^1dt~t^k
                               \frac{\partial\omega}{\partial b^j}(a,tb+(1-t)a).
\end{equation*}
En \'{e}crivant $\id_K^{-1}$ pour l'application identique $\Al\ra\Al$
et $\id_K^k$ pour l'application identique $CK^k\ra CK^k$ on
montre que
\begin{eqnarray}
       \epsilon h_K^{-1} & = & \id_K^{-1} \nonumber \\
       h_K^{-1}\epsilon + \dK^1 h_K^0 & = & \id_K^0 \nonumber\\
       h_K^{k-1}\dK^k + \dK^{k+1} h_K^k & = & \id_K^k~~~\forall k\geq 1
       \label{EqHHomotop}
\end{eqnarray}
ce qui entra\^{\i}ne l'acyclicit\'{e} du complexe de Koszul (\ref{EqCompKos}).

\bigskip

D\'{e}finissons enfin
les applications $F^k:CK^k\ra CH^k$ (voir e.g. \cite[p.~211]{Con94})
par $F^0=\id_H^0=\id_K^0$
et pour tout $\omega\in CK^k$:
\begin{equation*}
    (F^k\omega)(a,x_1,\ldots,x_k,b)=\omega(a,b)\big(x_1-a,\ldots,x_k-a\big)
\end{equation*}
quels que soient $k\in\entier$ avec $k\geq 1$ et
$a,x_1,\ldots,x_k,b\in\real^{n}$. Il est clair que les $F^k$
sont des morphismes de $\Al^e$-modules, et on montre qu'ils sont des morphismes de
complexes, c'est-\`{a}-dire, quel que soit $k\in\entier$,
\begin{equation*}
     F^k\dK^{k+1} = \dH^{k+1}F^{k+1}.
\end{equation*}
Il existe \'{e}galement des applications $G^k:CH^k\ra CK^k$ qui sont des
morphismes de $\Al^e$-modules et des morphismes de complexes~: on d\'{e}finit
$G^0=\id_H^0=\id_K^0$ et pour tout $\Phi\in CH^k$:
\begin{eqnarray}
(G^k \Phi)(a,b)& = &\sum_{i_1,\ldots,i_k=1}^n e^{i_1}\wedge\cdots\wedge e^{i_k}
               \int_{0}^1dt_1 \int_{0}^{t_1}dt_2 \int_{0}^{t_2}dt_3 \cdots
                \int_{0}^{t_{k-1}}dt_k \nonumber \\
           &    & ~~~~\frac{\partial^k \Phi}{\partial x_1^{i_1}\cdots \partial x_k^{i_k}}
                 \big(a,t_1a+(1-t_1)b,\ldots,t_ka+(1-t_k)b,b\big)\nonumber
                 \\
           &    & \label{EqDefGHochKos}
\end{eqnarray}
pour tout entier $k\geq 1$. Il est \'{e}vident que chaque $G^k$ est un
homomorphisme de $\Al^e$-modules, et \`{a} l'aide d'un calcul long mais direct
on montre que quel que soit $k\in\entier$,
\begin{equation*}
    G^k\dH^{k+1} = \dK^{k+1}G^{k+1}.
\end{equation*}
On peut repr\'{e}senter les deux applications $(F^k)_{k\in\entier}$ et
$(G^k)_{k\in\entier}$ dans le diagramme commutatif suivant~:
\begin{equation*}
 \begin{array}{ccccccc}
           \cdots &
              \stackrel{\dH^k}{\lla}& CH^k                       &
              \stackrel{\dH^{k+1}}{\lla}& CH^{k+1}                           &
              \stackrel{\dH^{k+2}}{\lla} &
              \cdots \\
           \cdots &
                                    & F^k\uparrow~\downarrow G^k &
                                        & F^{k+1}\uparrow~\downarrow G^{k+1} &
                                         &
              \cdots  \\
           \cdots &
              \stackrel{\dK^k}{\lla}& CK^k                      &
              \stackrel{\dK^{k+1}}{\lla}& CK^{k+1}                           &
              \stackrel{\dK^{k+2}}{\lla} &
              \cdots
 \end{array}
\end{equation*}

\begin{lemma} \label{LThetaProj}
 Avec les notations ci-dessus on a~:
 \begin{enumerate}
  \item $G^kF^k=\id_K^k$ quel que soit l'entier positif $k$.
  \item L'op\'{e}rateur $\Theta^k=F^kG^k:CH^k\ra CH^k$ est une projection, c'est-\`{a}-dire
       $\Theta^k\Theta^k=\Theta^k$ quel que soit l'entier positif $k$.
 \end{enumerate}
\end{lemma}
\begin{prooof}
  Les deux \'{e}nonc\'{e}s sont montr\'{e}s \`{a} l'aide d'un calcul direct.
\end{prooof}

\bigskip

Nous allons maintenant montrer l'existence d'homotopies $s_H^k:CH^k\ra
CH^{k+1}$: ce sont des homomorphismes de $\Al^e$-modules tels que
\begin{equation} \label{EqVereinfHomtop}
  \id_H^k-\Theta^k = \dH^{k+1}s_H^k + s_H^{k-1}\dH^k
\end{equation}
o\`{u} $s_H^{-1}=0$ et $s_H^0= 0$. Ceci se repr\'{e}sente de la fa\c{c}on
suivante:
\begin{equation*}
 \begin{array}{ccccccccc}
                                           \cdots                   &
    \stackrel{\dH^{k-1}}{\lla}  & CH^{k-1}                           &
    \stackrel{\dH^k}{\lla}     & CH^k                               &
    \stackrel{\dH^{k+1}}{\lla} & CH^{k+1}                           &
    \stackrel{\dH^{k+2}}{\lla} &           \cdots                   \\
                                       \cdots                   &
    \searrow  \hskip0.4cm\hbox{~}
             & \downarrow &
    {\searrow} \hskip0.4cm\hbox{~}
             & \downarrow     &
    \searrow  \hskip0.4cm\hbox{~}
             & \downarrow &
    \searrow \hskip0.4cm\hbox{~}    &           \cdots
    \\
                                       \cdots                   &
    {\scriptstyle s_H^{k-2}}
             & {\scriptstyle \id_H^{k-1}-\Theta^{k-1}} &
    {\scriptstyle s_H^{k-1}}
             & {\scriptstyle \id_H^{k}-\Theta^{k}}     &
    {\scriptstyle s_H^{k}}
             & {\scriptstyle \id_H^{k+1}-\Theta^{k+1}} &
    {\scriptstyle s_H^{k+1}}          &           \cdots                   \\
                                           \cdots                   &
    \hbox{~}\hskip0.4cm\searrow
             & \downarrow &
    \hbox{~}\hskip0.4cm{\searrow}
             & \downarrow     &
    \hbox{~}\hskip0.4cm\searrow
             & \downarrow &
   \hbox{~}\hskip0.4cm \searrow          &           \cdots                   \\
                                           \cdots                   &
    \stackrel{\dH^{k-1}}{\lla}  & CH^{k-1}                           &
    \stackrel{\dH^k}{\lla}     & CH^k                               &
    \stackrel{\dH^{k+1}}{\lla} & CH^{k+1}                           &
    \stackrel{\dH^{k+2}}{\lla} &           \cdots
 \end{array}
\end{equation*}
Puisque $\id_H^0-\Theta^0=0$ et
$(\id_H^k-\Theta^k)\dH^{k+1}=\dH^{k+1}(\id_H^{k+1}-\Theta^{k+1})$ on a
$\dH^1(\id_H^1-\Theta^1)=0$. On construit $s_H^1$ ``sur les g\'{e}n\'{e}rateurs''
de la mani\`{e}re suivante~: soit $F\in CH^1$. On consid\`{e}re $\widetilde{F}$ dans
$\widetilde{C}H^1=\CinfK{\real^{5n}}$ comme $\widetilde{F}(a',a,x,b,b')=F(a',x,b')$. On
prolonge $h_H^1$ et $\id_H^1-\Theta^1$ de $CH^1$ \`{a} tout \'{e}l\'{e}ment $T$ de
$\widetilde{C}H^1$ par $(h_H^1 T)(a',a,x_1,x_2,b,b')= T(a',a,x_1,x_2,b')$
(``les variables $a',b'$ ne sont pas affect\'{e}es'') et on fait de m\^{e}me avec
$\id_H^1-\Theta^1$. Ensuite on d\'{e}finit
\begin{eqnarray*}
    (s_H^1\Phi)(a,x_1,x_2,b)\!\!\!\!\!
      & = &\!\!\!\!
      (h_H^1(id_H^1-\Theta^1)\widetilde{\Phi})(a',a,x_1,x_2,b,b')|_{a'=a,b'=b}\\
      & = &\!\!\!\!\Phi(a,x_1,b) -\sum_{i=1}^n (x_1^i-a^i)\int_0^1
               \frac{\partial \Phi}{\partial x_1^i}(a,ta+(1-t)x_2,b)dt
\end{eqnarray*}
Par construction, $s_H^1$ est un homomorphisme de $\Al^e$-modules. \`A
l'aide de (\ref{EqHHomotop}) on voit que
\begin{equation*}
    \id_H^1-\Theta^1 = \dH^2 s_H^1.
\end{equation*}
On continue par r\'{e}currence~: soient $0=s_H^0,s_H^1,\ldots,s_H^k$ d\'{e}j\`{a}
construits tels que (\ref{EqVereinfHomtop}) soit satisfaite jusqu'\`{a} l'ordre $k$.
 Alors on a
\[
  \dH^{k+1}\big(\id_H^{k+1}-\Theta^{k+1}-s_H^k\dH^{k+1}\big)=0,
\]
et l'expression suivante (pour $\Phi\in CH^{k+1}$)
\begin{multline*}
    \lefteqn{(s_H^{k+1}\Phi)(a,x_1,\ldots,x_{k+2},b)
      =} \cr
       (h_H^{k+1}(id_H^{k+1}-\Theta^{k+1}-s_H^k\dH^{k+1})\widetilde{\Phi})
      \big(a',a,x_1,\ldots,x_{k+2},b,b'\big))|_{a'=a,b'=b}.
      \end{multline*}
est bien d\'{e}finie.
%
%
%
%
%
%
On a alors montr\'{e} le
\begin{lemma}
On peut construire des homotopies $s_H^k:CH^k\ra CH^{k+1}$ quel que soit
 l'entier positif $k$ telles que
 \begin{enumerate}
  \item Toute $s_H^k$ est un homomorphisme de $\Al^e$-modules.
  \item Toute $s_H^k$ est construite par une `suite d'op\'{e}rations qui
  consistent en des int\'{e}grales, des d\'{e}riv\'{e}es et des \'{e}valuations'.
  \item $s_H^0=0$.
  \item $\id_H^k-\Theta^k = \dH^{k+1}s_H^k + s_H^{k-1}\dH^k$ quel que soit
        $k\in \nat$.
 \end{enumerate}
\end{lemma}

\subsection{Calcul de la cohomologie}

\noindent
On consid\`{e}re toujours le cas $X=\real^n$.
Soient
$x=(x^1,\ldots,x^{n})$ les coordonn\'{e}es canoniques de $X$ et on va
\'{e}crire $x'$ pour $(x^1,\ldots,x^{n-l})$ et $x''$ pour
${x''}^1=x^{n-l+1},\ldots,{x''}^l=x^n$. On note $x=(x',x'')$. Il vient que
\begin{equation*}
    C=\{ x\in\real^n~|~x''=0\}.
\end{equation*}
On \'{e}crira parfois $E:=\real^n$, $E':=\real^{n-l}$ et $E'':=\real^l$.

Soit $\mathcal{M}$ un $\Al$-bimodule o\`{u} $\Al = \CinfK{\real^n}$.
On cherche \`{a} calculer la cohomologie de Hochschild \`{a} valeurs
dans $\mathcal{M}$.
Afin d'utiliser le complexe de Koszul du paragraphe pr\'{e}c\'{e}dent
il faut faire quelques hypoth\'{e}ses de r\'{e}gularit\'{e} concernant le
bimodule ainsi que quelques restrictions  pour  les cocha\^{\i}nes.
On suppose que $\mathcal{M}$ est muni de mani\`{e}re canonique de la structure
d'un
$\CinfK{\real^{2n}}$-module \`{a} gauche tel que
$(f \otimes g) \cdot \xi =
f \cdot \xi \cdot g$ pour tous $f, g \in \CinfK{\real^n}$ et $\xi\in\mathcal{M}$.
Dans tous nos exemples ceci est visiblement le cas, et en g\'{e}n\'{e}ral
le proc\'{e}d\'{e} marchera pour un bimodule topologique (par rapport \`{a} la topologie
de Fr\'{e}chet habituelle de $\CinfK{\real^n}$) par continuit\'{e}.

Pour les cocha\^{\i}nes, on se restreint aux cocha\^{\i}nes multidiff\'{e}rentielles
d\'efi\-nies comme suit:
une application $k$-lin\'{e}aire $\phi: \Al \times \cdots \times
\Al \longrightarrow \mathcal{M}$ est dite $k$-diff\'{e}rentielle (par rapport
\`{a} la multiplication de module \`{a} gauche) lorsqu'elle s'exprime comme
\begin{equation}
    \label{eq:MultiDifferential}
    \phi(f_1, \ldots, f_k) = \sum_{I_1, \ldots, I_k}
    \left(
        \frac{\partial^{|I_1|} f_1}{\partial x^{I_1}}
        \cdots
        \frac{\partial^{|I_k|} f_k}{\partial x^{I_k}}
    \right)
    \cdot \xi^{I_1 \cdots I_k}
\end{equation}
en consid\'erant une somme finie sur des multi-indices
$I = (i_1, \ldots, i_n)$ et des \'{e}l\'{e}ments du modules
$\xi^{I_1 \cdots I_k} \in \mathcal{M}$. De plus, on demande
pendant toute cette discussion
que la multiplication de module \`{a} droite
$f \mapsto \xi \cdot f$ (pour $\xi \in \mathcal{M}$ fix\'{e})
soit un op\'{e}rateur diff\'{e}rentiel dans le sens de
\eqref{eq:MultiDifferential}.
\begin{rem}
    \label{rem:ExamplesForBimodules}
    Pour les bimodules $\mathcal{M} = \Al$, $\Bl$ ainsi que
    $\mathcal{M} = \Dop(\Al, \Al)$, $\Dop(\Al, \Bl)$ and $\Dop(\Bl,
    \Bl)$ la d\'{e}finition ci-dessus donne la d\'{e}finition usuelle
    des cocha\^{\i}nes multidiff\'{e}rentielles et la condition concernant la multiplication
    \`{a} droite est \'{e}videmment satisfaite.
\end{rem}

Dans cette situation, on consid\`{e}re les morphismes de
$\CinfK{\real^{2n}}$-modu\-les
$A \in \Hom_{\CinfK{\real^{2n}}} (CH^k, \mathcal{M})$ tels que
l'application induite $a(f_1, \ldots, f_k) = A(1 \otimes f_1 \otimes
\cdots f_k \otimes 1)$
est une cocha\^{\i}ne multidiff\'{e}rentielle dans le sens mentionn\'{e} ci-dessus.
On note $\Hom^{\mathrm{diff}}_{\CinfK{\real^{2n}}}
(CH^\bullet, \mathcal{M})$ le $\korps$-espace engendr\'{e} par ces cocha\^{\i}nes.
Evidemment, la diff\'{e}rentielle du complexe bar pr\'{e}\-serve cet espace des
cocha\^{\i}nes multidiff\'{e}rentielles, et l'on obtient
un isomorphisme de complexes entre
l'espace des cocha\^{\i}nes diff\'{e}rentielles de Hochschild not\'{e} $H\Dop^\bullet (\mathcal{A},
\mathcal{M})$ et l'espace des cocha\^{\i}nes
$\Hom^{\mathrm{diff}}_{\CinfK{\real^{2n}}} (CH^\bullet,
\mathcal{M})$.

On consid\`{e}re \'{e}galement l'espace des morphismes de
$\CinfK{\real^{2n}}$-modu\-les,
$\Hom_{\CinfK{\real^{2n}}} (CK^\bullet, \mathcal{M})$, muni de la
diff\'{e}rentielle
$\dK$ du complexe de Koszul. Alors les applications $F$, $G$ et
les homotopies cha\^{\i}nes $s$ induisent des applications (encore not\'{e}es par $F$, $G$
et $s$) entre les cocha\^{\i}nes correpondantes.  L'\'{e}nonc\'{e} important suivant est
moins \'{e}vident:
\begin{lemma}
    \label{lemma:GsAreNice}
    L'application induite
    \begin{equation}
        \label{eq:GniceMap}
        G:
        \Hom_{\CinfK{\real^{2n}}}
        (CK^\bullet, \mathcal{M})
        \longrightarrow
        \Hom^{\mathrm{diff}}_{\CinfK{\real^{2n}}}
        (CH^\bullet, \mathcal{M})
    \end{equation}
    prend ses valeurs dans l'espace des cocha\^{\i}nes diff\'{e}rentielles.
    De plus, l'application induite
    $s$ pr\'{e}serve les cocha\^{\i}nes diff\'{e}rentielles
    \begin{equation}
        \label{eq:sNiceMap}
        s:\Hom^{\mathrm{diff}}_{\CinfK{\real^{2n}}}
        (CH^\bullet, \mathcal{M})
        \longrightarrow
        \Hom^{\mathrm{diff}}_{\CinfK{\real^{2n}}}
        (CH^{\bullet-1}, \mathcal{M}).
    \end{equation}
\end{lemma}
\begin{prooof}
    Pour $G$, c'est tr\`{e}s facile \`{a} montrer en utilisant la forme explicite
    (\ref{EqDefGHochKos}) car on n'a qu'\`{a} l'\'{e}valuer \`{a} la fin en
    $a = x_1 = \cdots x_k = b$. Pour $s$, c'est l\'{e}g\`{e}rement plus compliqu\'{e}
    puisqu'\`{a} premi\`{e}re vue $s$ est `non-locale'. Cependant, par r\'{e}currence
    on montre qu'apr\`{e}s l'\'{e}valuation $a = x_1 = \cdots = x_k = b$ seuls
    les termes `locaux' survivent et ceux-ci  sont donn\'{e}s par des
    d\'{e}riv\'{e}es.
\end{prooof}
De ce lemme, on d\'{e}duit imm\'{e}diatement l'\'{e}nonc\'{e} g\'{e}n\'{e}ral suivant:
\begin{theorem}
    \label{theorem:ComputeHH}
    Avec les hypoth\`{e}ses sur le bimodule $\mathcal{M}$ faites ci-dessus,
    les groupes de cohomologie de Hochschild diff\'{e}rentiels de
    $\Al$ \`{a} valeurs dans
    $\mathcal{M}$ sont isomorphes aux groupes de cohomologie du complexe
    \begin{equation}
        \label{eq:KoszulHomComputeHH}
        H\Dop^\bullet(\Al, \mathcal{M}) \cong
        H\left(
            \Hom^{\mathrm{diff}}_{\CinfK{\real^{2n}}}
            (CK^\bullet, \mathcal{M}),
            \dK
        \right)
    \end{equation}
    via les applications induites $F$ et $G$, qui sont des r\'{e}ciproques.
\end{theorem}
\begin{rem}
\label{subsec:NotAtAllExt}
   On remarque que le raisonnement usuel pour montrer cet \'{e}nonc\'{e}, \`{a}
   sa voir l'ind\'{e}pendence des $Ext$-groupes vis-\`a-vis du choix de la r\'{e}solution
   projective, ne s'applique pas puisqu'on utilise une classe restreinte
   de cocha\^{\i}nes et puisque le complexe bar n'est pas du tout form\'e de modules
   projectifs sur l'anneau $\CinfK{\real^{2n}}$.
\end{rem}
\begin{rem}
    \label{subsec:ContHH}
    Le m\^{e}me \'{e}nonc\'{e} reste toujours vrai quand on remplace la cohomologie
    de Hochschild diff\'{e}rentielle par la cohomologie de Hochschild
    continue (o\`{u} l'on utilise des cocha\^{\i}nes continues par rapport \`{a} la
    topologie de Fr\'{e}chet de
    $\Al = \CinfK{\real^n}$ et un bimodule topologique
    $\mathcal{M}$). Dans ce cas l\`{a}, la d\'{e}monstration du
    lemme~\ref{lemma:GsAreNice} est beaucoup plus facile comme $G$ et $s$
    sont des applications continues d\'{e}j\`{a} sur le plan des cha\^{\i}nes.
\end{rem}

Par cons\'{e}quent, il faut calculer la cohomologie du complexe de Koszul
$\Hom_{\CinfK{\real^{2n}}} (CK^\bullet, \mathcal{M})$ pour la diff\'{e}rentielle induite
$\dK$. Cette derni\`ere est tr\`{e}s facile \`{a} simplifier car
la diff\'{e}rentielle se calcule de fa\c{c}on explicite:
\begin{lemma}
    On a l'isomorphisme de $\CinfK{\real^{2n}}$-modules
    \begin{equation}
        \label{eq:KoszulMultivektors}
        \Hom_{\CinfK{\real^{2n}}} (CK^\bullet, \mathcal{M})
        \cong
        \mathcal{M} \otimes_{\korps} \Lambda^\bullet (\real^n),
    \end{equation}
    et la diff\'{e}rentielle est donn\'{e}e par
    \begin{equation}
        \label{eq:InducedKoszulDifferential}
        \dK X = e_i \wedge (x^i \cdot X - X \cdot x^i),
    \end{equation}
    o\`{u} $x^1, \ldots, x^n \in C^\infty(\real^n)$ sont des coordonn\'{e}es
    usuelles par rapport \`{a} la base $e_1, \ldots, e_n$ de
    $\real^n$.
\end{lemma}
\begin{prooof}
    L'\'{e}nonc\'{e}~\eqref{eq:KoszulMultivektors} est \'{e}vident car
    $CK^\bullet$ est un $\CinfK{\real^{2n}}$-module libre.
    En dualisant le produit int\'{e}rieur $i(\xi)$ avec le champ de vecteurs
    $\xi(a, b) = a-b$ on obtient la multiplication ext\'{e}rieure $\wedge$ avec $\xi$
    dans le sens de la multiplication de
    $\CinfK{\real^{2n}}$-modules.
    C'est exactement la diff\'{e}rentielle \eqref{eq:InducedKoszulDifferential}.
\end{prooof}

En particulier, si le bimodule
$\mathcal{M}$ est {\em sym\'{e}trique}, i.e.
les multiplications \`{a} gauche et \`{a} droite de
$\Al$ sur $\mathcal{M}$ co\"{\i}ncident, il vient que la
diff\'{e}rentielle ci-dessus s'annule de sorte que
la cohomologie de Hochschild est donn\'{e}e par
\begin{equation}
    \label{eq:HHBimoduleSymmetric}
    H\Dop(\Al, \mathcal{M}) \cong
    \mathcal{M} \otimes_{\korps} \Lambda^\bullet (\real^n),
\end{equation}
o\`{u} les isomorphismes sont donn\'{e}s par
$F$ et $G$. Ceci implique d\'{e}j\`{a} les r\'{e}sultats suivants:
\begin{theorem}
    \label{THKRfacil}
    Avec les notations mentionn\'{e}es ci-dessus:
    \begin{enumerate}
    \item  $H\Dop^k(\Al,\Al)=H\Gs^k\cong
        \Ginf(\Lambda^k TX)=\gs^k$ (Hochschild-Kostant-Rosen\-berg).
    \item $H\Dop^k(\Al,\Bl)\cong \Bl\otimes_{\korps}\Lambda^k \real^n =
        \Ginf(\Lambda^k TX|_C)$.
    \end{enumerate}
\end{theorem}
On aura besoin du r\'{e}sultat technique suivant bien connu concernant {\em l'homologie
de Hochschild continue} de $\Al$, voir \'{e}galement \cite[p.51]{Pfl98} et \cite{Nad99}:
on regarde le
complexe qui est la somme directe des tous les espaces
$C_k(\Al,\Al):=\CinfK{\real^{(k+1)n}}$ dont les \'{e}l\'{e}ments sont des
cha\^{\i}nes de Hochschild `topologiques', i.e. des fonctions lisses
$\phi:(a,x_1,\ldots,x_k)\mapsto \phi(a,x_1,\ldots,x_k)$ \`{a} $k+1$ arguments
dans $\real^n$. Soit
$p_k:CH^k\ra C_k(\Al,\Al)$ la surjection lin\'{e}aire continue donn\'{e}e par
$(p_k\Phi)(a,x_1,\ldots,x_k)=\Phi(a,x_1,\ldots,x_k,a)$.
L'op\'{e}rateur bord de Hochschild $\mathsf{b}'$
de degr\'{e} $-1$ dans $C_k(\Al,\Al)$ est d\'{e}fini comme dans l'\'{e}quation
(\ref{EqDefBordHochBar}) `en enlevant le dernier argument et mettant
$b=a$' de sorte que l'on obtienne une version topologique de la
diff\'{e}rentielle \cite[1.1.1, p.9]{Lo92}. Il s'ensuit que les applications
$p_k$ d\'{e}finissent des morphismes de complexes. Soit $\Omega_\bullet=\Al\otimes
\Lambda^\bullet \real^{n*}$ vu comme complexe avec la diff\'{e}rentielle $0$. L'application
$\tilde{p}_k:CK^k\ra \Omega_k$ d\'{e}finie par
$(\tilde{p}_k\omega)(a):=\omega(a,a)$ est \'{e}galement un morphisme des
complexes. On trouve des uniques applications $F'_k:\Omega_k\ra
C_k(\Al,\Al)$ et $G'_k:C_k(\Al,\Al)\ra \Omega_k$ induites par
les applications $F_k$ et $G_k$ entre les complexes bar et de Koszul
(i.e. $F'_k\circ\tilde{p}_k=p_k\circ F_k$ et
$G'_k\circ p_k=\tilde{p}_k\circ G_k)$ ainsi que des homotopies induites $s^{'\prime k}_H$.
Il vient que
les complexes $\big((C(\Al,\Al),\mathsf{b}'\big)$ et
$\big(\Omega,0\big)$ sont quasi-isomorphes, donc on a
\begin{prop} \label{PropHomolContHKR}
 L'homologie de Hochschild continue de $\Al$ est donn\'{e}e par
 \[
     HH_{k~\mathrm{cont}}(\Al,\Al)\cong \Al\otimes \Lambda \real^{n *}
       ~~~(Hochschild-Kostant-Rosenberg).
 \]
\end{prop}

Pour les autres bimodules dont on a besoin dans cet article, on rappelle
le fait suivant:
\begin{lemma}\label{LemSuiteExacteDiffOp}
 La suite exacte de $\Al^e$-modules
 \[
     \{0\}\lra \Il \lra \Al \lra \Bl \lra \{0\}
 \]
 entra\^{\i}ne la suite exacte de $\Al^e$-modules
 \[
     \{0\}\lla \Dop(\Il,\Bl) \lla \Dop(\Al,\Bl) \lla \Dop(\Bl,\Bl) \lla \{0\}
 \]
 et finalement, la suite exacte de complexes
 \[
     \{0\}\lla \Dop^k\big(\Al,\Dop(\Il,\Bl)\big)
     \lla \Dop^k\big(\Al,\Dop(\Al,\Bl)\big) \lla
     \Dop^k\big(\Al,\Dop(\Bl,\Bl)\big)\lla \{0\}
 \]
\end{lemma}
\begin{prooof}
 Il est clair qu'un op\'{e}rateur diff\'{e}rentiel de $\Al$ dans $\Bl$ s'annule
 sur $\Il$ si et seulement s'il ne contient que des d\'{e}riv\'{e}es partielles
 des variables $x'$, et est donc un prolongement d'un op\'{e}rateur
 diff\'{e}rentiel de $\Bl$ dans $\Bl$.
\end{prooof}
Pour la cohomologie de Hochschild \`{a} valeurs dans des espaces d'op\'{e}ra\-teurs
diff\'{e}rentiels comme ceux mentionn\'{e}s dans le lemme pr\'{e}c\'{e}dent \ref{LemSuiteExacteDiffOp},
on ne peut plus utiliser
\eqref{eq:HHBimoduleSymmetric} puisque la structure de bimodule {\em n'est
pas} sym\'{e}tri\-que. On utilise la d\'{e}composition
$\real^n = E' \oplus E''$ mentionn\'{e}e ci-dessus. On a le
\begin{theorem}\label{TheoCalculCohomol}
    Les groupes de cohomologie des complexes suivants
    se simplifient comme suit~:
    \begin{enumerate}
    \item $H\Dop^k(\Al,\Dop(\Al,\Bl))\cong \{0\}$ si $k\geq 1$
        et
        $H\Dop^0 (\Al, \Dop(\Al, \Bl)) \cong \Bl$.
    \item $H\Dop^k(\Al,\Dop(\Bl,\Bl))\cong \Bl\otimes_\korps \Lambda^k {E''}$ quel
        que soit l'entier $k\geq 0$.
    \item $H\GsP^k=H\Dop^{k-1}(\Al,\Dop(\Il,\Bl))\cong H\Dop^k(\Al,\Dop(\Bl,\Bl))
    \cong \Bl\otimes_\korps \Lambda^k {E''}=\gsP^k$ quel
    que soit l'entier $k\neq 0$.
    \item $H\GsP^0\cong \Bl$.
\end{enumerate}
\end{theorem}
\begin{prooof}
  On esquisse la preuve du deuxi\`{e}me \'{e}nonc\'{e}, celle des autres \'{e}tant
  enti\`{e}rement analogue. Soit
  $\Phi \in \Dop(\Bl, \Bl)$, i.e.
    \[
    (\Phi f)(x') =
    \sum_I \Phi^I(x') \frac{\partial^{|I|}f}{\partial (x')^I}(x')
    \]
    avec des fonctions lisses $\Phi^I \in \CinfK{E'} = \Bl$.
    De mani\`{e}re \'{e}quivalente, on peut d\'{e}crire $\Phi$ par son \emph{symbole}
    $\tilde{\Phi}(x', p') = \sum_I \Phi^I(x') p'_I$ o\`{u} $p'_I =
    (p'_1)^{i_1} \cdots $ $(p'_{n-l})^{i_{n-l}}$, alors $\tilde{\Phi}
    \in \CinfK{T^*E'}$ est une application polyn\^{o}miale
     en les variables d'impulsion
    $p'_{1}, \ldots, p'_{n-l}$. Il vient que l'application $\Phi \leftrightarrow
    \tilde{\Phi}$ est \'{e}videmment une bijection $\korps$-lin\'{e}aire. on
    utilise la d\'{e}composition
    \[
    \Dop(\Bl, \Bl) \otimes_{\korps} \Lambda^\bullet (\real^n)
    \cong
    \Dop(\Bl, \Bl) \otimes_{\korps} \Lambda(E') \otimes_{\korps} \Lambda(E'')
    \]
    pour identifier $\Phi \in \Dop(\Bl, Bl) \otimes_{\korps}
    \Lambda^\bullet (\real^n)$ avec son symbole $\sum\tilde{\Phi} =
    \frac{1}{k!}  \tilde{\Phi}^{i_1 \cdots i_k} e_{i_1}$ $ \wedge \cdots
    \wedge e_{i_k}$ o\`{u} $\tilde{\Phi}^{i_1 \cdots, i_k} \in
    \CinfK{E'}\otimes \mathsf{S}\real^{n}$ comme avant. La diff\'erentielle
    de Koszul \eqref{eq:InducedKoszulDifferential} s'exprime comme suit
    \[
    \dK \tilde{\Phi}
    = \sum_{i=1}^n e_i \wedge \{x^i, \tilde{\Phi}\}
    = \sum_{i=1}^{n-l} e_i \wedge \{x_i, \tilde{\Phi}\}
    = d'_p \tilde{\Phi}
    \]
    o\`{u} $\{~,~\}$ d\'{e}signe le crochet de Poisson usuel dans
    $T^*E'=E'\times E^{\prime *}$,
    $d'_p$ d\'{e}signe la diff\'{e}rentielle de deRham dans la direction des
    impulsions $p'\in E^{\prime *}$ et
    les $e_1, \ldots, e_{n-l}$ s'identifient de mani\`{e}re canonique avec les $1$-formes
    $dp'_1,\ldots, dp'_{n-l}$. Gr\^{a}ce au lemme de Poincar\'{e}, on en d\'{e}duit le
    r\'{e}sultat.
\end{prooof}

En utilisant la suite exacte longue de cohomologie r\'{e}sultant de la suite
exacte courte des complexes $\GsI\ra\Gs\ra\GsP$ et le fait que
l'homomorphisme connectant s'annule, on obtient finalement le r\'{e}sultat
souhait\'{e}~:
\begin{theorem}
\label{theoHKR}
  Quel que soit l'entier $k$:
 \begin{enumerate}
  \item $H\Gs^k\cong \Al\otimes_\korps \Lambda^kE=\gs$ (th\'{e}or\`{e}me HKR classique).
  \item $H\GsP^k\cong \Bl \otimes_\korps \Lambda^k {E''} =\gsP^k$.
  \item $H\GsI^k\cong \gsI^k$.
 \end{enumerate}
\end{theorem}

%
%
%
%
\subsection{(Co)homologie de Hochschild de $\gs$ et de $\gsP$}

On aura besoin de calculer plus tard la cohomologie (diff\'{e}rentielle) et
l'homologie continue gradu\'{e}es des alg\`{e}bres commutatives gradu\'{e}es
$\gs=\Al\otimes \Lambda E$ et $\gsP=\Bl\otimes \Lambda E''$.
On rappelle le cobord de Hochschild gradu\'{e} pour une alg\`{e}bre associative
gradu\'{e}e $G=\oplus_{i\in\mathbb Z}G^i$:
\begin{eqnarray}\label{EqDefCobHgradue}
 (\mathsf{b}\phi)(f_1,\ldots,f_{k+1}) & = &
     (-1)^{|f_1|~|\phi|} f_1~\phi(f_2,\ldots,f_{k+1})  \nonumber \\
          &  & +\sum_{r=1}^{k}(-1)^{r}
             \phi(f_1,\ldots,f_rf_{r+1},\ldots,f_{k+1})\nonumber \\
          & &  +(-1)^{k+1}\phi(f_1,\ldots,f_{k},f_{k+1}).
\end{eqnarray}
pour une cocha\^{\i}ne $\phi\in\Hom(G^{\otimes k},G)$ de degr\'{e} $|\phi|$
et $f_1,\ldots,f_{k+1}\in G$.

\begin{theorem}\label{TheoHKRpourggsP}
 On a les r\'{e}sultats HKR suivants pour la cohomologie diff\'{e}\-rentielle gradu\'{e}e et
 l'homologie continue pour les alg\`{e}bres associatives commutatives gradu\'{e}es
 $\gs$ et $\gsP$ pour tout entier positif $k$:
 \begin{enumerate}
  \item $HH^k(\gs,\gs)\cong
  \gs\otimes\oplus_{p=0}^k \big(\Lambda^p E \otimes \mathsf{S}^{k-p} E^*\big)$.
  \item $HH^k(\gsP,\gsP)\cong
  \gsP\otimes\oplus_{p=0}^k \big(\Lambda^p E' \otimes
               \mathsf{S}^{k-p} E^{\prime\prime *}\big)$.
  \item $HH_k(\gs,\gs)\cong
  \gs\otimes\oplus_{p=0}^k \big(\Lambda^p E^* \otimes \mathsf{S}^{k-p} E\big)$.
  \item $HH_k(\gsP,\gsP)\cong
  \gsP\otimes\oplus_{p=0}^k \big(\Lambda^p E^{\prime *}
    \otimes \mathsf{S}^{k-p} E^{\prime\prime}\big)$.
 \end{enumerate}
\end{theorem}
\begin{prooof}
$\Al$ et $\Bl$ sont des espaces des fonctions $\CinfK{\real^N}$ pour $N=n$
et $N=n-l$. On va calculer les cohomologies pour $\Al\otimes \Lambda W$
pour un $\korps$-espace vectoriel de dimension finie et donc obtenir les \'{e}nonc\'{e}s
comme cas particuliers.\\
Puisque $W$ et donc $G:=\Lambda W$ sont de dimension finie, on peut
utiliser les r\'{e}solutions libres de $G^e$-modules o\`{u} $G^e$ est l'alg\`{e}bre
$G\otimes G^{\mathrm{opp}}$ avec la multiplication
\'{e}vidente $(g_1\otimes g_2)(g'_1\otimes g'_2)
:=(-1)^{|g_2|(|g'_1|+|g'_2|)}g_1 g'_1\otimes g'_2 g_2$ et $G$ est
consid\'{e}r\'{e} comme $G^e$-module gradu\'{e} de fa\c{c}on usuelle, i.e.
$(g_1\otimes g_2)g:=(-1)^{|g_2|~|g|}g_1g g_2$. Dans la cat\'{e}gorie des
$G^e$-modules gradu\'{e}s,
la cohomologie de Hochschild gradu\'{e}e de
$G$ \`{a} valeurs dans $G$ est donn\'{e}e par les groupes $Ext^k_{G^e}(G,G)$
et l'homologie de Hochschild gradu\'{e}e de $G$ \`{a} valeurs dans $G$ est
donn\'{e}e par les groupes $Tor^{G^e}_k(G,G)$,
voir e.g. \cite[p.169]{CE56} pour le cas non gradu\'{e} qui se g\'{e}n\'{e}ralise sans
probl\`{e}me. La premi\`{e}re r\'{e}solution libre de $G$ comme $G^e$-module est donn\'{e}e par
le complexe bar usuel (voir \cite[p.174]{CE56}) d\'{e}fini par la famille
$CH(G)^k:=(G\otimes G^{\otimes k}\otimes G)_{k\in\mathbb N}$ dont l'op\'{e}rateur bord
est
\begin{eqnarray}\label{EqDefCobordHochsGrad}
    \lefteqn{\hat{\partial}^k_H(g\otimes g_1\otimes\cdots\otimes g_k \otimes
    g')
      := gg_1\otimes g_2\otimes\cdots\otimes g_k \otimes g'}
                    \nonumber \\
      & &   +\sum_{r=1}^{k-1}(-1)^r g\otimes g_1\otimes\cdots\otimes g_rg_{r+1}\otimes \cdots
                                \otimes g_k \otimes g'
                    \nonumber \\
      & &   +(-1)^k g\otimes g_1\otimes\cdots\otimes g_{k-1} \otimes g_{k}g'
\end{eqnarray}
sans modification de signe.
La r\'{e}solution de Koszul gradu\'{e}e est d\'{e}finie par la famille de $G^e$-modules gradu\'{e}s
libres $G^e\otimes
\mathsf{S}^kW$, et l'op\'{e}rateur bord $\partial_K^k: G^e\otimes
\mathsf{S}^kW\ra G^e\otimes\mathsf{S}^{k-1}W$ suivant
\begin{eqnarray*}
   \lefteqn{\partial^k_K(g\otimes L\otimes g'):=}\\
      & & (-1)^{|g'|}\sum_{i=1}^{\dim W}\big(g\wedge e_i\otimes i(e^i)(L)\otimes
      g'~-~g \otimes i(e^i)(L)\otimes e_i\wedge g'\big).
\end{eqnarray*}
On voit que --avec cette r\'{e}solution-- les homomorphismes de $G^e$-modules gradu\'{e}s
\`{a} valeurs dans $G$ annulent
$\partial_K^k$, ainsi que les produits tensoriels gradu\'{e}s sur $G^e$ avec
le $G^e$-module $G$, alors
\[
    HH^k(G,G)\cong Ext^k_{G^e}(G,G)\cong \Lambda W\otimes \mathsf{S}W^*
\]
et
\[
    HH_k(G,G)\cong Tor^{G^e}_k(G,G)\cong \Lambda W\otimes \mathsf{S}W.
\]
Pour finalement calculer les (co)homologies de Hochschild pour l'alg\`{e}bre produit tensoriel
$\Al\otimes \Lambda W$ on regarde d'abord la r\'{e}solution `topologique'
donn\'{e}e par le complexe acyclique d\'{e}fini par la famille
$(CH^k\otimes CH(G)^k)_{k\in\mathbb N}$ et les
op\'{e}rateurs bord $\partial^k_H\otimes \hat{\partial}^k_H$. Puisque les deux
complexes $(CH^k)_{k\in\mathbb N}$ et $(CH(G)^k)_{k\in\mathbb N}$ sont
des modules simpliciaux gradu\'{e}s (voir \cite[p.44]{Lo92}), le th\'{e}or\`{e}me
d'Eilenberg-Zilber (voir e.g. \cite[p.238]{Mac75}) nous permet de
remplacer le complexe acyclique $(CH^k\otimes CH(G)^k)_{k\in\mathbb N}$ (dont
l'homologie est \'{e}gale \`{a} $\Al\otimes G$) par le
produit tensoriel des complexes acycliques
$\big((CH^k)_{k\in\mathbb N},(\partial^k_H)_{k\in\mathbb N}\big)$ et
$\big((CH(G)^k)_{k\in\mathbb N},(\hat{\partial}^k_H)_{k\in\mathbb N}\big)$
(dont
l'homologie est aussi \'{e}gale \`{a} $\Al\otimes G$ gr\^{a}ce \`{a} un argument simple de bicomplexes,
voir la proposition \ref{PropBicomplexeElementaire})
comme r\'{e}solution `topologique' de $\Al\otimes \Lambda W$. Le quasi-isomorphisme
 d'Alexander-Whitney (voir e.g. \cite[p.241]{Mac75}) entre ces deux
 complexes est \'{e}galement un morphisme de $\Al^e\otimes G^e$-modules.
En
consid\'{e}rant les homomorphismes diff\'{e}rentiels dans $\Al\otimes G$ et les produits tensoriels
topologiques avec $\Al\otimes G$, respectivement, avec cette derni\`{e}re r\'{e}solution, il r\'{e}sulte
du th\'{e}or\`{e}me de K\"{u}nneth (voir e.g. \cite[p.166]{Mac75}) que la
(co)homologie de $\Al\otimes \Lambda W$ est le produit tensoriel sur
$\korps$ de la (co)homologie de $\Al$ et de celle de $G=\Lambda W$,
d'o\`{u} le r\'{e}sultat d'apr\`{e}s ce qui pr\'{e}c\`{e}de, d'apr\`{e}s la proposition \ref{PropHomolContHKR}
et d'apr\`{e}s le th\'{e}or\`{e}me \ref{THKRfacil}.
\end{prooof}

\subsection{Les applications HKR}

\noindent Dans cette sous-section nous allons construire des
quasi-isomorphismes entre les
alg\`{e}bres de Lie diff\'{e}rentielles gradu\'{e}es $(\GsI,\mathsf{b})$ et
$(\gsI,0)$. Nous verrons que
l'application HKR usuelle correspondant \`{a} l'antisym\'{e}trisation
ne convient pas et doit \^{e}tre modifi\'{e}e.

Pour un entier positif $k$ soit $\alpha^k:\Lambda^k E\ra E^{\otimes k}$
l'application d'antisym\'{e}trisa\-tion usuelle:
$v_1\wedge\cdots\wedge
  v_k\mapsto \frac{1}{k!}\sum_{\sigma\in S_k}\epsilon(\sigma)
     v_{\sigma(1)}\otimes\cdots\otimes v_{(\sigma(k))}$. Soit aussi
$P^1:SE\ra E$ la projection canonique.
En utilisant le fait que $\Gs^k\cong \Al \otimes_\korps (SE)^{\otimes k}$
(\`{a} l'aide du calcul symbolique)
nous d\'{e}finissons
les deux applications suivantes
 \begin{equation}
     \psi^k_{\scriptscriptstyle HKR}:\gs^k\ra \Gs^k: f\otimes T \mapsto
     f\otimes
     \alpha^k(T),
\end{equation}
et
\begin{equation}
   \pi^k_{\scriptscriptstyle HKR}:\Gs^k\ra \gs^k:
   f\otimes L_1\otimes\cdots\otimes L_k\mapsto f\otimes P^1(L_1)\wedge
     \cdots\wedge P^1(L_k).
\end{equation}
Nous \'{e}crirons
$\psi_{\scriptscriptstyle HKR}$ (resp. $\pi_{\scriptscriptstyle HKR}$)
pour la somme de tous
les $\psi^k_{\scriptscriptstyle HKR}$ (resp. $\pi^k_{\scriptscriptstyle
HKR}$) o\`{u} $\psi^k_{\scriptscriptstyle HKR}$ et $\pi^k_{\scriptscriptstyle
HKR}$ s'annulent pour $k\leq -1$ et $k>\dim E$.
Le th\'{e}or\`{e}me de Hochschild, Kostant et Rosenberg (HKR), \ref{theoHKR},
permet de d\'{e}duire le suivant
\begin{theorem}\label{applHKRusuel}
 Quel que soit l'entier positif $k$, les applications
 $\psi^k_{\scriptscriptstyle HKR}$ et
 $\pi^k_{\scriptscriptstyle HKR}$ ont les propri\'{e}t\'{e}s suivantes:
 \begin{enumerate}
  \item $\pi^k_{\scriptscriptstyle HKR}\psi^k_{\scriptscriptstyle HKR}=
          \mathrm{id}_{\gs^k}$.
  \item $\mathsf{b}\psi^k_{\scriptscriptstyle HKR}=0$.
  \item Soit $\phi\in\Gs^k$. Si $\mathsf{b}\phi=0$ et
      $\pi^k_{\scriptscriptstyle HKR}\phi=0$, alors $\phi$ est un cobord,
      i.e. il existe $\phi'\in\Gs^{k-1}$ tel que $\phi=\mathsf{b}\phi'$.
  \item Soient $X,Y\in\gs$, alors il existe $\xi\in \Gs$ tel que
    \[
      [\psi_{\scriptscriptstyle HKR}X,\psi_{\scriptscriptstyle HKR}Y]_G
      -\psi_{\scriptscriptstyle HKR}[X,Y]_S=\mathsf{b}\xi.
    \]
  \item Soient $X\in\gs^k$ et $Y\in\gs^l$, alors il existe $\xi\in \Gs$ tel
      que
      \[
        \psi_{\scriptscriptstyle HKR}X\cup\psi_{\scriptscriptstyle HKR}Y
      -(-1)^{kl}
       \psi_{\scriptscriptstyle HKR}Y\cup\psi_{\scriptscriptstyle HKR}X
      -\psi_{\scriptscriptstyle HKR}(X\wedge Y)=\mathsf{b}\xi.
      \]
 \end{enumerate}
\end{theorem}
Malheureusement, l'application HKR $\psi_{\scriptscriptstyle HKR}$ n'envoie
pas le sous-espace $\gsI$ de $\gs$ dans le sous-espace $\GsI$ de $\Gs$:
en utilisant la d\'{e}composition
\begin{equation}
 \gsI ~~=~~ \Al\otimes_\korps \Lambda {E''}\otimes_\korps \Lambda^+ {E'}~~~\oplus~~~
         \Il \otimes_\korps \Lambda {E''}
\end{equation}
on voit que l'image de $\psi_{\scriptscriptstyle HKR}$ contiendrait des
\'{e}l\'{e}ments
provenant du premier terme de la somme ci-dessus dont le facteur le plus
\`{a} droite est dans ${E''}$ sans que le coefficient soit dans l'id\'{e}al, et un
tel op\'{e}rateur multidiff\'{e}rentiel ne serait plus dans $\GsI$.

 Il faut alors modifier les applications $\alpha^k$: il est bien connu que
 la multiplication ext\'{e}rieure
induit un isomorphisme d'espaces vectoriels
\begin{equation}
    \Phi:\Lambda {E''}\otimes_\korps \Lambda {E'} \ra \Lambda ({E'}\oplus {E''}) =\Lambda E:
      T_2\otimes T_1 \mapsto T_2\wedge T_1.
\end{equation}
Soit $\iota^{(l,k-l)}:{E''}^{\otimes l}\otimes_\korps {E'}^{\otimes (k-l)}\ra
E^{\otimes k}$ l'injection induite par des sous-espaces ${E'},{E''}$.
L'application $\hat{\alpha}^k:=\sum_{l=0}^k\iota^{(l,k-l)}\alpha^l\otimes
\alpha^{k-l}$ envoie $\big(\Lambda {E''}\otimes_\korps \Lambda {E'}\big)_k$ dans
$E^{\otimes k}$. Nous d\'{e}finissons l'application HKR modifi\'{e}e par
\begin{equation}
  {\psi^1}^k_{\scriptscriptstyle HKR}:\gs^k\ra \Gs^k: f\otimes T \mapsto
     f\otimes
     \hat{\alpha}^k\big(\Phi^{-1}(T)\big),
\end{equation}
et $\psi^1_{\scriptscriptstyle HKR}$ comme \'{e}tant la somme des tous les
${\psi^1}^k_{\scriptscriptstyle HKR}$. Ainsi le facteur le plus \`{a} droite
est toujours dans ${E'}$ dans le cas o\`{u} son coefficient est dans $\Al$, et
n'est dans ${E''}$ que si son coefficient est dans $\Il$. Par cons\'{e}quent,
$\psi^1_{\scriptscriptstyle HKR}$ envoie bien $\gsI$ dans $\GsI$.
En \'{e}crivant $\psi^{[1]}$ pour la restriction de
$\psi^1_{\scriptscriptstyle HKR}$
\`{a} $\gsI$ on a l'analogue suivant du
th\'{e}or\`{e}me \ref{applHKRusuel}:
\begin{theorem}\label{applHKR}
 Quel que soit l'entier positif $k$, les applications
 ${\psi^{[1]}}^k$ et
 $\pi^k_{\scriptscriptstyle HKR}$ ont les propri\'{e}t\'{e}s suivantes:
 \begin{enumerate}
  \item $\pi^k_{\scriptscriptstyle HKR}{\psi^{[1]}}^k =
          \mathrm{id}_{\gsI^k}$.
  \item $\mathsf{b}{\psi^{[1]}}^k=0$.
  \item Soit $\phi\in\GsI^k$. Si $\mathsf{b}\phi=0$ et
      $\pi^k_{\scriptscriptstyle HKR}\phi=0$, alors $\phi$ est un cobord,
      i.e. il existe $\phi'\in\GsI^{k-1}$ tel que $\phi=\mathsf{b}\phi'$.
  \item Soient $X,Y\in\gsI$, alors il existe $\phi\in \GsI$ tel que
   \[
      [\psi^{[1]}X,\psi^{[1]}Y]_G
      -\psi^{[1]}[X,Y]_S=\mathsf{b}\phi.
   \]
  \item Soient $X\in\gs^k$ et $Y\in\gs^l$, alors il existe $\xi\in \Gs$ tel
      que
      \[
        \psi^{[1]}X\cup\psi^{[1]}Y
      -(-1)^{kl}
       \psi^{[1]}Y\cup\psi^{[1]}X
      -\psi^{[1]}(X\wedge Y)=\mathsf{b}\xi.
      \]
 \end{enumerate}
\end{theorem}
\begin{prooof}
 1. Le premier \'{e}nonc\'{e} est \'{e}vident. \\
 2. Puisque l'image de $\psi^{[1]}$ consiste en
 des op\'{e}rateurs multidiff\'{e}rentiels \\
 $1$-diff\'{e}rentiels, cet espace
 ne contient que des cocycles, d'o\`{u} le deuxi\`{e}me \'{e}nonc\'{e}.\\
 3. D'apr\`{e}s le troisi\`{e}me \'{e}nonc\'{e} du th\'{e}or\`{e}me \ref{applHKRusuel}, $\phi$
 est un cobord dans $\Gs$, mais puisque le morphisme connectant
 de la suite exacte longue correpondant \`{a} la suite exacte courte
 $\GsI\ra \Gs \ra \GsP$ s'annule, il s'ensuit que $\phi$ est aussi un
 cobord dans $\GsI$.\\
 4. D'apr\`{e}s 1. et le premier \'{e}nonc\'{e} de \ref{applHKRusuel},
 on trouve $\phi_1,\phi_2,\phi'\in\Gs$ tels que
 $\psi^{[1]}X=\psi_{\scriptscriptstyle HKR}X+\mathsf{b}\phi_1$,
 $\psi^{[1]}Y=\psi_{\scriptscriptstyle HKR}Y+\mathsf{b}\phi_2$ et
 $\psi^{[1]}[X,Y]=\psi_{\scriptscriptstyle HKR}[X,Y]+\mathsf{b}\phi'$.
 A l'aide de
 l'\'{e}nonc\'{e} 4. de \ref{applHKRusuel}, on voit que le membre de gauche de 4.
 est un \'{e}l\'{e}ment de $\GsI$ et un cobord dans $\Gs$, donc --\`{a} l'aide du
 m\^{e}me argument que pour 3.-- un cobord dans $\GsI$.\\
 5. Raisonnement analogue \`{a} celui pour 4.
\end{prooof}

\section{Formalit\'{e}, star-produits adapt\'{e}s et une structure $G_\infty$ sur $\GsI$}
  \label{SecFormAdapGinf}

Dans cette section, on rappelle d'abord la construction d'un star-produit
(adapt\'{e}) \`{a} l'aide d'n $L_\infty$-morphisme due \`{a} Kontsevitch. Puis, dans
la deuxi\-\`{e}me sous-section, nous donnons les \'{e}tapes de la construction
d'un tel morphisme en utilisant des structures $G_\infty$.
Nous  nous adapterons, au
cas o\`{u} l'espace est $\GsI$, la preuve de
Tamarkin  de la conjecture de Deligne \cite{Del93}
(il existe une stucture $G_\infty$ sur $\Gs$).

\subsection{Formalit\'{e} et star-produits}

Dans ce paragraphe, $\hs$ est ou bien l'espace gradu\'{e}
$\gs=\Ginf(X,\Lambda TX)$ ou son sous-espace gradu\'{e} $\gsI$ et $\Hs$ est ou
bien l'espace gradu\'{e} des op\'{e}rateurs multidiff\'{e}rentiels $\Gs$ ou son
sous-espace $\GsI$.

D'apr\`{e}s les sections pr\'{e}c\'{e}dentes, l'espace $\hs[1]$ est muni de la structure d'alg\`{e}bre
de Lie gradu\'{e}e
$[~,~]_S$ (le crochet de Schouten (\ref{EqDefSchouten})) et $\Hs[1]$ est muni de la structure
d'alg\`{e}bre de Lie diff\'{e}rentielle gradu\'{e}e $\big([~,~]_G,\mathsf{b}\big)$
(crochet de Gerstenhaber (\ref{EqDefCrochGerstDiff}) et diff\'{e}rentielle de
Hochschild (\ref{EqDefCobH})). Il s'ensuit que $\hs[1]$ est muni d'une structure
$L_\infty$, qu'on note $d_L$ et $\Hs[1]$ est muni d'une structure $L_\infty$ qu'on note
$D_L$, voir l'appendice \ref{AppSubSecFInfini} pour
des d\'{e}tails. Un {\em morphisme de formalit\'{e}} est une somme $\psi=\sum_{k=1}^\infty\psi^{[k]}$
d'applications $\psi^{[k]}:\mathsf{S}^k(\hs[2])\ra\Hs[2]$ qui sont de
degr\'{e} $0$ et des op\'{e}rateurs $k$-diff\'{e}rentiels dans
$\Dop^k\big(\hs;\Hs)$ et qui satisfont \`{a} la condition suivante: le
morphisme de cog\`{e}bres cocommutatives gradu\'{e}es $\overline{\psi}:\mathsf{S}(\hs[2])
\ra \mathsf{S}(\Hs[2])$ co-induit par $\psi$ pr\'{e}serve les
cod\'{e}rivations $\overline{d_L}$ sur $\mathsf{S}(\hs[2])$ et $\overline{D_L}$
sur $\mathsf{S}(\Hs[2])$, i.e. $\overline{\psi}\circ
\overline{d_L}=\overline{D_L}\circ \overline{\psi}$.

\noindent La construction suivante est un r\'{e}sultat c\'{e}l\`{e}bre de Kontsevitch \cite{Kon03}:

\begin{theorem} \label{TheoKontsevitch}
 Soit $\psi$ un morphisme de formalit\'{e} entre les alg\`{e}bres $L_\infty$ $\hs$
 et $\Hs$.
 Il existe un star-produit $\star$ sur $X$ (dans le cas $\hs=\gs$)
 qui est adapt\'{e} \`{a} la sous-vari\'{e}t\'{e} $C$ (dans le cas $\hs=\gs_I$).
\end{theorem}
\begin{prooof}
 Soit $P\in\gs[2][[h]]$ une structure de Poisson formelle, i.e.
 $[P,P]_S=0$. Dans la cog\`{e}bre topologique $\mathsf{S}(\hs[2])[[h]]$ la structure $P$ est de
 degr\'{e} $0$ et la s\'{e}rie $e^{h\bullet P}$ (o\`{u} $\bullet$ d\'{e}signe la multiplication
 shuffle, voir l'appendice \ref{AppSubSecGradues}) est un \'{e}l\'{e}ment de genre groupe. Le
 morphisme de cog\`{e}bres cocommutatives $\overline{\psi}$ envoie $e^{h\bullet P}$ sur un autre
 \'{e}l\'{e}ment de genre groupe dans la cog\`{e}bre topologique
 $\mathsf{S}(\Hs[2])[[h]]$. Gr\^{a}ce \`{a} la colibert\'{e} de $\mathsf{S}(\Hs[2])$
 il vient que l'image $\overline{\psi}\big(e^{h\bullet P}\big)$ est n\'{e}cessairement de la
 forme exponentielle $e^{h\bullet m_\star}$ avec $m_*$ de degr\'{e} $0$ dans
 $\Hs[2]$, donc une s\'{e}rie d'op\'{e}rateurs bidiff\'{e}rentiels. Puisque $P$ est
 une structure de Poisson, il vient $\overline{d_L}\big(e^{h\bullet
 P}\big)=0$. Par cons\'{e}quent, $\overline{D_L}\big(e^{h\bullet m_\star}\big)=0$
 car $\overline{\psi}$ pr\'{e}serve les cod\'{e}rivations. La composante par
 rapport \`{a} $\Hs[2]$ de cette derni\`{e}re \'{e}quation donne
 \[
     0= D_L^{[1]}(h m_\star)+\frac{1}{2}D_L^{[2]}(h^2m_\star\bullet m_\star)
      =\mathsf{b}(h m_\star)+\frac{1}{2}[hm_\star, hm_\star]_G
 \]
 qui est app\'{e}l\'{e}e `l'\'{e}quation Maurer-Cartan' et \'{e}quivaut \`{a} l'associativit\'{e} de
 la s\'{e}rie $\star:=m+hm_\star$ d'op\'{e}rateurs bidiff\'{e}rentiels
 o\`{u} $m$ est la multiplication point-par-point
 dans $\CinfK{M}$: en effet,
 $\star\circ_G\star=[m,hm_\star]_G+\frac{1}{2}[hm_\star,hm_\star]_G$ et
 l'op\'{e}rateur cobord de Hochschild $\mathsf{b}$ est donn\'{e} par le crochet de
 Gerstenhaber avec $m$.
\end{prooof}

\noindent Le but de notre \'{e}tude est donc de trouver un tel morphisme de formalit\'{e} $\psi$
entre les $L_\infty$-alg\`{e}bres $\gsI$ et $\GsI$. Un argument de
perturbation homologique montre que l'on peut toujours trouver un tel
morphisme $\psi$ si on admet une modification $d'_L$ de la structure $L_\infty$,
$d_L$,
sur $\hs$ avec
des termes d'ordre sup\'{e}rieur, voir la proposition \ref{Theorem 2.1}. Mais
l'apparition de $d'_L$ changerait drastiquement l'argument de Kontsevitch
ci-dessus. Pour retransformer $d'_L$ en $d_L$ par un $L_\infty$-morphisme, les
obstructions se trouvent dans le cocomplexe
$\big(\Dop(\mathsf{S}^+(\hs[2]),\hs[2]),[d_L,~]_{NR}\big)$, voir la
proposition \ref{Theorem 3.1}. Malheureusement, cette cohomologie de
Chevalley-Eilenberg de $(\hs[1],[~,~]_S)$ n'est pas concentr\'{e} en $0$, et  il
faut alors  construire $\psi$ par d'autres moyens. Dans le cas $\hs=\gs$,
Kontsevitch a d\'{e}velopp\'{e} une m\'{e}thode g\'{e}om\'{e}trique \`{a} l'aide de graphes, voir
\cite{Kon03}. A cause de l'asym\'{e}trie des op\'{e}rateurs bidiff\'{e}rentiels
adapt\'{e}s dans $\GsI$, la restriction du morphisme de formalit\'{e} de Kontsevitch \`{a} $\gsI$
ne prend pas ses valeurs dans $\GsI$. On suivra la m\'{e}thode de Tamarkin
\cite{Tam98p}, \cite{GH03}: ceci explique le d\'{e}tour par des structures
$G_\infty$.

\subsection{Structures $G_\infty$ sur $\GsI$}

Le lecteur peut trouver les d\'{e}tails et conventions de signe pr\'{e}cises pour les
structures \`{a} homotopie pr\`{e}s dans
l'appendice \ref{AppSubSecFInfini}.

Soit $\as$ un espace vectoriel gradu\'{e}, par exemple $\gs,\gsI,\Gs$ ou $\GsI$.
La {\em cog\`{e}bre de Gerstenhaber colibre}
sur $\as[2]$, est la `cog\`{e}bre sym\'{e}trique colibre en la cog\`{e}bre de Lie colibre
avec un certain d\'{e}calage', plus pr\'{e}cis\'{e}ment, l'espace gradu\'{e}
$\mathsf{S}\big((\underline{\as[1]^\otimes})[1]\big)$ avec deux structures alg\`{e}briques
compatibles: une comultiplication coassociative cocommutative gradu\'{e}e et un
cocrochet de Lie gradu\'{e} sur $\mathsf{S}\big((\underline{\as[1]^\otimes})[1]\big)[-1]$.

\vspace{0.5cm}

{\bf 1. Structures $G_\infty$}: une structure $G_\infty$ sur $\hs:=\as[1]$ est
d\'{e}finie par une cod\'{e}riva\-tion $\overline{d}$ (pour les deux structures alg\`{e}briques de
de $\mathsf{S}\big((\underline{\hs^\otimes})[1]\big)$ de degr\'{e} $1$ et de carr\'{e}
$0$ co-induite par $d\in \Hom\big(\mathsf{S}\big((\underline{\hs^\otimes})[1]\big),\hs[1]\big)$
(voir l'appendice \ref{AppSubSecFInfini} pour plus de d\'{e}tails). Pour simplifier la notation,
on va l'\'{e}crire dans sa version d\'{e}cal\'{e}e $m:=d[-1]=\sum_{r, p_1,\ldots,p_r=1}^\infty
m^{p_1,\ldots,p_r}$ o\`{u} les applications $m^{p_1,\ldots,p_r}$ de degr\'{e} $2-r$ sont des
op\'{e}rateurs $p_1+\cdots +p_r$-diff\'{e}rentiels de
$\underline{\hs^{\otimes p_1}}\vea\cdots\vea\underline{\hs^{\otimes p_n}}$ dans $\hs$.
%
Ici, $m^1$ est une codiff\'{e}rentielle, $m^{1,1}$ est une structure d'alg\'{e}bre
de Lie \`{a} homotopie pr\`{e}s sur $\hs$, et $m^2[-1]$ est une structure d'alg\`{e}bre commutative
associative \`{a} homotopie pr\`{e}s sur $\as=\hs[-1]$.

\vspace{0.5cm}

{\bf 2. Structures $G_\infty$ sur les multivecteurs}:
puisque toute alg\`{e}bre de Gerstenhaber est muni d'une structure $G_\infty$ o\`{u}
$m=m^1+m^{1,1}+m^2$, les espaces des multivecteurs (adapt\'{e}s), $\gs[1]$ et
$\gsI[1]$ sont automatiquement munis d'une structure $G_\infty$
donn\'{e}e par $m^1=0$, $m^{1,1}=[~,~]_S$ et $m^2[-1]=\wedge$.

\vspace{0.5cm}

Il n'est pas aussi facile de trouver une structure $G_\infty$ sur les
espaces d'op\'{e}rateurs multidiff\'{e}rentiels $\Gs[1]$ ou $\GsI[1]$ car les
lois $\mathsf{b}$, $[~,~]_G$  et $\cup[1]$ n'en font pas des alg\`{e}bres de Gerstenhaber:
la r\`{e}gle de Leibniz entre
$[~,~]_G$ et $\cup$ n'est satisfaite qu'\`{a} un cobord de Hochschild pr\`{e}s,
voir \cite[p.285, Theorem 5]{Ger63}.
Tamarkin proc\`{e}de de fa\c{c}on suivante:

\vspace{0.5cm}

{\bf 3. Structure de big\`{e}bre de Lie codiff\'{e}rentielle sur la cog\`{e}bre de Lie
colibre implique structure $G_\infty$}: c'est la proposition
\ref{Lemma 1.1} (montr\'{e} dans \cite{Tam98p} et \cite{GH03}): soit $\hs$ un espace gradu\'{e}
et soit $([~,~],\mathsf{b})$ la structure
d'une big\`{e}bre de Lie gradu\'{e}e codiff\'{e}rentielle sur la cog\`{e}bre de Lie libre
$\underline{\hs^\otimes}$ (voir l'appendice \ref{AppSubSecFInfini} pour
des d\'{e}finitions), co-induite par une somme d'applications
$l^{p_1,p_2}:\underline{\hs^{\otimes p_1}}\otimes \underline{\hs^{\otimes p_1}}\ra\hs$
et $l^p:\underline{\hs^{\otimes p}}\ra \hs$. Alors cette structure se prolonge en une
structure
$G_\infty$ sur $\hs$ avec $m^{p_1,p_2}=l^{p_1,p_2}$ et $m^p=l^p$, les autres
composantes \'{e}tant z\'{e}ro.

\vspace{0.5cm}

{\bf 4. Les accolades d\'{e}finissent la structure d'une big\`{e}bre codiff\'{e}\-rentielle sur les
cog\`{e}bres colibres $\Gs[1]^\otimes$ et $\GsI[1]^\otimes$}: les deux espaces
d'op\'{e}ra\-teurs multidiff\'{e}rentiels $\Gs$ et $\GsI$ sont des sous-espaces
des $\korps$-homo\-mor\-phis\-mes $\Al^\otimes\ra\Al$ avec
$\Al=\Al^0=\CinfK{X}$, alors $\Gs[1],\GsI[1]$ sont des sous-espaces gradu\'{e}s de
$\Hom(\Al[1]^\otimes,\Al[1])$.
 De plus, $\Gs[1]$ et $\GsI[1]$ sont ferm\'{e}s pour les
op\'{e}rations $\circ_i$ (\ref{EqImeComposition}) gr\^{a}ce \`{a} la proposition
\ref{PGsI}.
Dans cette situation
g\'{e}n\'{e}rale, il existe une structure canonique de big\`{e}bre sur la cog\`{e}bre
colibre $\Hs^\otimes$ (o\`{u} par exemple $\Hs=\Gs[1]$ ou $\Hs=\GsI[1]$) d\'{e}finie par des accolades
({\em braces}), voir l'appendice \ref{AppSubSecAccolades} pour des
d\'{e}tails. Ici la comultiplication reste la comultiplication de
d\'{e}concat\'{e}nation, tandis que la multiplication $\bullet_K$ (\ref{EqDefMultKadeiAccol})
se pr\'{e}sente comme une
certaine `d\'{e}formation de la multiplication shuffle $\bullet$' si l'on prend le degr\'{e}
tensoriel comme filtration. La multiplication est co-induite par certaines
applications $m_K=\sum_{p_1,p_2=1}^\infty a^{p_1,p_2}$ avec
$a^{p_1,p_2}:\Hs^{\otimes p_1}\otimes \Hs^{\otimes p_2}\ra\Hs$. Le
commutateur gradu\'{e} (par rapport \`{a} $\bullet_K$) avec la d\'{e}cal\'{e}e $m[1]$ de la multiplication
associative $m$
sur $\Al$ d\'{e}finit une codiff\'{e}rentielle $\mathsf{b}_K$ (\ref{EqDefCobordBKadei})
de la big\`{e}bre qui est co-induite par des applications $a^p:\Hs^{\otimes p}\ra\Hs$)
avec $p=1,2$: $a^1$ est \'{e}gale au cobord de Hochschild $\mathsf{b}$ sur $\Hs$
et $a^2$ est \'{e}gal \`{a} la d\'{e}cal\'{e}e $\cup[1]$ de la multiplication $\cup$.

\vspace{0.5cm}

{\bf 5. Le th\'{e}or\`{e}me de d\'{e}quantification d'Etingof-Kazhdan implique la
structure de big\`{e}bre de Lie sur la cog\`{e}bre de Lie colibre}~: la filtration
par degr\'{e} tensoriel de la big\`{e}bre accolades mentionn\'{e}e
dans {\bf 4.} n'est pas bonne: en passant \`{a} la big\`{e}bre gradu\'{e}e associ\'{e}e, on
obtiendrait
--en divisant par les shuffles $\Hs^{\otimes +}\bullet\Hs^{\otimes +}$, voir l'appendice
\ref{AppSubSecGradues}-- la structure d'une big\`{e}bre de Lie gradu\'{e}e
sur la cog\`{e}bre de Lie colibre, mais la composante $m^2$ serait nulle, donc
on n'aurait pas inclus la multiplication $\cup$. Ce n'est qu'une application du th\'{e}or\`{e}me de
d\'{e}quantification d'Etingof-Kazhdan (voir l'appendice de l'article
\cite{GH03}) qui permet de passer de la big\`{e}bre codiff\'{e}rentielle accolades
(o\`{u} la codiff\'{e}rentielle est le commutateur gradu\'{e} avec la multiplication
point-par-point $m$) \`{a} la structure d'une big\`{e}bre de Lie sur la cog\`{e}bre de
Lie colibre qui inclut la multiplication {\em cup} co-induite par $a^2$:
\begin{prop}\label{Theorem 1.2}
Soit $\Hs$ un espace vectoriel gradu\'{e}.
Supposons donn\'{e}e une structure de big\`{e}bre diff\'{e}rentielle
sur la cog\`{e}bre tensorielle colibre
$\Hs^\otimes=\oplus_{n\geq 0}~\Hs^{\otimes n}$
dont la diff\'{e}rentielle et la
multiplication sont donn\'{e}es respectivement par des applications
$a^n$~: $\Hs^{\otimes n}\rightarrow \Hs$ et
$a^{p_1,p_2}$~: $\Hs^{\otimes p_1}\otimes \Hs^{\otimes p_2}
\rightarrow \Hs$. Alors on a une structure
de big\`{e}bre de Lie diff\'{e}rentielle sur la cog\`{e}bre de Lie
$\underline{\Hs^\otimes}=\oplus_{n\geq 1}~\underline{\Hs^{\otimes n}}$,
dont la diff\'{e}rentielle et le crochet de Lie sont donn\'{e}s
respectivement par les applications $l^n$
et $l^{p_1,p_2}$ o\`{u} $l^1=a^1$, $l^2$ est l'anti-sym\'{e}tris\'{e}e de $a^2$ et
$l^{1,1}$ est l'anti-symmetris\'{e}e de $a^{1,1}$.
\end{prop}

\vspace{0.5cm}

{\bf 6. Bonne structure $G_\infty$ sur $\Gs[1]$ et $\GsI[1]$}~: c'est une cons\'{e}quence
directe des trois \'{e}tapes pr\'{e}c\'{e}dentes {\bf 3.}, {\bf 4.}, {\bf 5.}.
On note la structure $G_\infty$ sur $\Hs=\Gs[1]$ ou $\Hs=\GsI[1]$ par $D$.

\vspace{0.5cm}

{\bf 7. Morphisme de formalit\'{e} (diff\'erentiel) $G_\infty$ entre $(\gsI[1],d')$ et
$(\GsI[1],$ $D)$}~: un raisonnement g\'{e}n\'{e}ral de perturbation homologique (voir
la proposition \ref{Theorem 2.1}) permet de construire une structure $G_\infty$
diff\'{e}rentielle
$d'$ sur $\hs=\gs[1]$ ou $\hs=\gsI[1]$ et un $G_\infty$-morphisme diff\'{e}rentiel de
formalit\'{e} $\overline{\psi}:(\hs,d')\ra (\Hs, D)$ \`{a} partir de l'application
HKR $\psi_{HKR}$ (pour $\Gs$) ou l'application HKR modifi\'{e}e $\psi^{[1]}$
(pour $\GsI$). Les propri\'{e}t\'{e}s de $\psi^{[1]}$ dans
le th\'{e}or\`{e}me \ref{applHKR} permettent de conclure que la structure
$G_\infty$ $d'$ est donn\'{e}e par une modification de la structure usuelle
$m=d[-1]=[~,~]_S+\wedge[1]$ par des termes de rang sup\'{e}rieur
$d^{p_1,\ldots,p_r}$ avec $p_1+\cdots+p_r\geq 3$.

\vspace{0.5cm}

{\bf 8. Le cocomplexe d'obstructions pour transformer $d'$ en $d$}~:
encore une fois, il faut transformer la structure $G_\infty$ $d'$ sur
$\hs$ en la structure usuelle $d$ par un morphisme $G_\infty$ diff\'{e}rentiel $\psi'$
de la cog\`{e}bre de Gerstenhaber colibre de $\hs$ dans elle-m\^{e}me.
La proposition \ref{Theorem 3.1}
d\'{e}crit le cocomplexe d'obstructions \`{a} cette transformation comme l'espace
\[
  \left(\!\Dop({\ve^\cdot\underline{\gsI[1]^{\otimes \cdot}}}, \gsI[1]),
\big[[-,-]_S+\wedge[1],-\big]\right),
\]
pour le cas $\hs=\gsI[1]$.
On va montrer dans la section suivante \ref{SecCalculObstructions}
que la cohomologie de ce cocomplexe
est concentr\'{e} en $0$ donc cette transformation sera toujours possible.

\vspace{0.5cm}

{\bf 9. Morphisme de formalit\'{e} $G_\infty$ entre $(\gsI[1],d)$ et
$(\GsI[1],D)$}~: la composition $\overline{\psi}\circ \overline{\psi'}$ des
deux $G_\infty$-morphismes diff\'{e}rentiels $\overline{\psi}$ ({\bf 7.}) et
$\overline{\psi'}$ ({\bf 8.}) donne
le r\'{e}sultat.

\vspace{0.5cm}

{\bf 10. $L_\infty$-morphisme de formalit\'{e} diff\'{e}rentiel $\gsI[1]\ra \GsI[1]$ \`{a} partir
de la restriction du $G_\infty$-morphisme}~:
afin d'appliquer le th\'{e}or\`{e}me de Kontsevitch \ref{TheoKontsevitch} pour la
construction d'un star-produit adapt\'{e} (dont les op\'{e}rateurs bidiff\'{e}rentiels
r\'{e}sident dans $\GsI$), il faut obtenir le morphisme $L_\infty$ diff\'{e}rentiel \`{a} partir
du morphisme $G_\infty$ diff\'{e}rentiel: c'est toujours possible gr\^{a}ce \`{a} la proposition
\ref{PropGinfiniLinfini}. On obtient les structures $L_\infty$ et le
$L_\infty$-morphisme de formalit\'{e} par restriction aux cog\`{e}bres sym\'{e}triques
vues comme sous-cog\`{e}bres de Gerstenhaber, voir l'appendice
\ref{AppSubSecFInfini}. Ainsi la restriction de $d=d^{1,1}+d^2$ nous donne le crochet
de Schouten $d_L^2=d^{1,1}=[~,~]_S[1]$, et la restriction de
$D=\sum_{p\geq 1}D^p+\sum_{p_1,p_2\geq 1}D^{p_1,p_2}$ donne
$D^1+D^{1,1}=\mathsf{b}+[~,~]_G[1]$. La restriction du
$G_\infty$-morphisme diff\'{e}rentiel de formalit\'{e} d\'{e}finit donc le $L_\infty$-morphisme
diff\'{e}rentiel de formalit\'{e} souhait\'{e}.

\vspace{0.5cm}

{\bf 11. La preuve du th\'{e}or\`{e}me principal \ref{TheoremePrincipal}} est
termin\'{e}e apr\`{e}s le th\'{e}or\`{e}me \ref{theoHKR}, les \'{e}tapes {\bf 1.} \`{a} {\bf 10.} et
l'acyclicit\'{e} du cocomplexe d'obstructions dans le paragraphe
\ref{SecCalculObstructions}.

\section{Calcul des obstructions}
 \label{SecCalculObstructions}

Nous avons vu que les obstructions
\`{a} la construction de star-repr\'{e}senta\-tions r\'{e}sidaient dans le
groupe de cohomologie du {\em cocomplexe d'obstructions}
\begin{equation}\label{EqDefCompObstruction}
   \left(\!\Dop({\ve^\cdot\underline{\gsI[1]^{\otimes \cdot}}}, \gsI[1]),
\big[[-,-]_S+\wedge[1],-\big]\right),
\end{equation}
voir la proposition \ref{Theorem 3.1}.
Dans cette section, nous nous proposons
de calculer ces groupes d'obstructions et de montrer qu'ils s'annulent.

Cette t\^{a}che sera divis\'{e}e en deux parties~: le cocomplexe d'obstructions
(\ref{EqDefCompObstruction}) est un {\em bicomplexe} quand on
munit un \'{e}l\'{e}ment
$x$ de $\underline{\gsI[1]^{\otimes p_1}}\vea \cdots \vea
 \underline{\gsI[1]^{\otimes p_n}}$
du bidegr\'{e} additionnel $(\sum_{i=1}^n p_i-1, n-1)$: la codiff\'{e}rentielle
$D_{CE}:=[d^{1,1},~]_T$ provenant du crochet de Schouten $[~,~]_S$ est de bidegr\'{e}
$(0,1)$ et anticommute avec la codiff\'{e}rentielle $D_{Har}:=[d^2,~]_T$ provenant de la
multiplication $\wedge$
qui est de bidegr\'{e} $(1,0)$ (voir \cite{Tam98p}, \cite{GH03}, \cite{Gin04}).
Le th\'{e}or\`{e}me \'{e}l\'{e}mentaire d'alg\`{e}bre homologique suivant (voir par exemple
\cite[p.25, Cor. 37]{And74})
sera tr\`{e}s utile pour la suite:
\begin{prop} \label{PropBicomplexeElementaire}
La
cohomologie d'un bicomplexe par rapport \`{a} la deuxi\-\`{e}me codiff\'{e}rentielle est
toujours un cocomplexe induit pour la premi\`{e}re codiff\'{e}ren\-tielle. Si
cette cohomologie par rapport \`{a} la deuxi\`{e}me codiff\'{e}rentielle est
concentr\'{e}e en degr\'{e} $0$, alors la cohomologie totale du bicomplexe
est isomorphe \`{a} la cohomologie
du cocomplexe induit.
\end{prop}

\smallskip

Dans la premi\`{e}re sous-section,
nous montrerons que l'homologie de Harrison de $\gsI$ est concentr\'{e}e en degr\'{e} $1$
(apr\`es d\'ecalage), et
que ce r\'{e}sultat implique que le bicomplexe d'obstructions
(\ref{EqDefCompObstruction}) pour
\ la codiff\'{e}rentielle $D_{Har}$ est acyclique. Ceci entra\^{\i}ne une premi\`{e}re r\'{e}duction du
bicomplexe d'obstructions au cocomplexe induit
\begin{equation}\label{EdDefCocomplexeInduit}
  \left({\Hom}_{{\gs}}(\ve_{{\gs}}^\cdot
\big(\gs \ot_{{\gsI}_+}\Omega_{{\gsI}_+}\big),\gsI),\delta^{1,1} \right),
\end{equation}
o\`{u} $\delta^{1,1}$ est la codiff\'{e}rentielle induite de la codiff\'{e}rentielle
$D_{CE}$, l'alg\`{e}bre ${\gsI}_+$ est  $\korps 1\oplus \gsI$ et
$\gs \ot_{{\gsI}_+}\Omega_{{\gsI}_+}$ d\'{e}signe le $\gs$-module
des diff\'{e}rentielles de K\"{a}hler de ${\gsI}_+$. En fait, cette cohomologie
est une version gradu\'{e}e de la cohomologie de Chevalley-Eilenberg de
l'alg\`{e}bre de Lie gradu\'{e}e des diff\'{e}rentielles de K\"{a}hler.

\smallskip

Dans la deuxi\`{e}me sous-section, nous montrerons que
si $X=\real^n$ et $C=\real^{n-l}$ la cohomologie
du cocomplexe induit (\ref{EdDefCocomplexeInduit}) est concentr\'{e} en degr\'{e} $0$.
On dira parfois abusivement qu'un complexe est acyclique pour dire qu'il est
concentr\'{e} en degr\'{e} $0$ (ou bidegr\'{e} $(0,0)$).

\subsection{
    Premi\`{e}re r\'{e}duction du cocomplexe d'obstructions \`{a} l'ai\-de de l'homologie de Harrison}

\noindent  Dans cette partie on s'attache \`{a} d\'{e}montrer
\begin{theorem}
\label{propcohomolo}
 Avec les notations introduites ci-dessus:
 \begin{enumerate}
  \item L'homologie de Harrison $\Har_{\cdot}({\gsI}_+|\korps, \gs)$ est
  isomorphe \`{a} l'espace des diff\'{e}rentielles de K\"{a}hler
  $\gs\ot_{{\gsI}_+}\Omega^1_{{\gsI}_+}$.
  \item Le cocomplexe d'obstructions (\ref{EqDefCompObstruction}) se r\'{e}duit au
  cocomplexe induit (\ref{EdDefCocomplexeInduit}) donn\'{e} par
   $\left({\Hom}_{{\gs}}(\ve_{{\gs}}^\cdot
\big(\gs \ot_{{\gsI}_+}\Omega_{{\gsI}_+}\big),\gsI),\delta^{1,1} \right)$
 \end{enumerate}
\end{theorem}
\begin{prooof}
L'astuce principale est de faire appara\^{\i}tre les $\gs$-modules~:
pour tout espace vectoriel gradu\'{e} $V$, le produit tensoriel $\gs\ot_\korps V$
est un $\gs$-module \`{a} gauche libre engendr\'{e} par $V$
(par multiplication sur le premier facteur).
On a \'{e}videmment un isomorphisme d'espace vectoriels complexes~:
\[
  {\rm Hom}_\korps(\gu,\gsI)
    \cong
    {\rm Hom}_{{\gs}}(\gextgsI,\gsI)
\]
que l'on peut prolonger au produits tensoriels topologiques car on ne
consi\-d\`{e}re que les cocha\^{\i}nes multidiff\'{e}rentiels, donc continues.
La codiff\'{e}rentielle  induite sur le dernier complexe par $[m^2,-]$ est la duale
d'une diff\'{e}rentielle induite par
$\delta^2$ sur $\gextgsI$ qui n'est autre que la {\em diff\'{e}rentielle de
Harrison $\beta$} sur chaque facteur
$\gs\ot \underline{\gsI^{\otimes \cdot}}$. En effet, pour $\varphi:
\gextgsI\to \gsI$ et $x\ot x_1\ot\cdots \ot x_n\in \gextgsI$, on a
\begin{align*}
     [m^2,\varphi](x\ot x_1\ot\cdots \ot x_n)&\\
           &\hskip-2cm
            = \pm m^2(x_1,\varphi(x\ot x_2\cdots))
                \pm m^2(\varphi(x\ot x_1 \cdots),x_n)\\
           &\hskip-1.5cm
              +\sum \pm\varphi(x\ot x_1\cdots m^2(x_i,x_{i+1})\cdots)\\
           &\hskip-2cm
                = \varphi(m^2(x,x_1)\ot x_2\cdots )
                +\sum \pm\varphi(x\ot \cdots m^2(x_i,x_{i+1})\cdots)\\
           &\hskip-2cm
                =\varphi(\beta(x\ot x_1\ot\cdots \ot x_n)).
\end{align*}
Une it\'{e}ration du th\'{e}or\`{e}me des bicomplexes \ref{PropBicomplexeElementaire}
assure que si l'homologie du complexe (tronqu\'{e}, i.e. $k\geq 1$) de Harrison
$\big(\gs\ot \underline{{\gsI}_+^{\ot \cdot}}, \beta\big)$ de $\gsI$
\`{a} coefficients dans $\gs$  est  concentr\'{e}e en degr\'{e} 1 (dans notre cas
\'{e}gale \`{a} $\gs \ot_{{\gsI}_+}\Omega^1_{{\gsI}_+}$), alors l'homologie
du complexe $\big(\ve^k_\gs\big(\gs\otimes\underline{\gsI^{\otimes}}\big),\beta\big)$ est
acyclique $\forall k\geq 1$, et finalement, la cohomologie  du cocomplexe
d'obstructions (\ref{EqDefCompObstruction}) par rapport \`{a} la codiff\'{e}rentielle $D_{Har}$
est acyclique (voir \'{e}galement \cite{Tam98p}, \cite{GH03}).

%
%
%
%
Rappelons que ${\gsI}_+=\korps\oplus \gsI$. Montrons maintenant que
l'homologie de
Harrison est \'{e}gale \`{a} $\gs \ot_{{\gsI}_+}\Omega^1_{{\gsI}_+}$.
On a une suite d'inclusions de sous-alg\`{e}bres gradu\'{e}es commutatives et
unitaires
\[
 \korps\hookrightarrow {\gsI}_+ \hookrightarrow \gs
\]
Comme $\korps$ est de caract\'{e}ristique $0$,
pour tout $\gs$-module $M$, on peut \'{e}crire la version gradu\'{e}e topologique de la
{\em suite exacte de Jacobi-Zariski
associ\'{e}e} (o\`{u} la notation $B|A$ signifie que l'anneau commutatif gradu\'{e}
$B$ est vu comme alg\`{e}bre sur le sous-anneau gradu\'{e} $A$), voir par exemple
\cite[p.61-p.74]{And74} pour des d\'{e}monstrations:
\begin{multline*}
\cdots \Har_{\cdot+1}(\gs|{\gsI}_+, M)\to \Har_{\cdot}({\gsI}_+|\korps, M)\\
\to \Har_{\cdot}(\gs|\korps, M)\to \Har_{\cdot}(\gs|{\gsI}_+, M)\to \cdots .
\end{multline*}

Rappelons que cette propri\'{e}t\'{e} provient du fait que lorsque $A\supset \QM$,
l'homologie de Harrison
$\Har_\cdot (B|A,M)$ d'une $A$-alg\`{e}bre
$B$ \`{a} valeurs dans un $B$-module $M$ est \'{e}gale \`{a} l'homologie
d'Andr\'{e}-Quillen $AQ_\cdot(B|A,M)$ (\`{a} un d\'{e}calage du degr\'{e}
pr\`{e}s) qui admet toujours une telle suite exacte. De plus, elle s'identifie aussi avec
la partie de poids $1$  de la
d\'{e}composition
de Hodge de l'homologie de Hochschild $HH_{\cdot}(B|A,M)$, voir \cite[p.145, Prop.
4.5.13]{Lo92}.
On s'int\'{e}resse au cas $M=\gs$.
Gr\^{a}ce \`{a} la version gradu\'{e}e du th\'{e}or\`{e}me de Hochschild-Kostant-Rosenberg pour l'homologie
continue de l'alg\`{e}bre gradu\'{e}e $\gs=\CinfK{\real^n}\otimes \Lambda E$ (voir le th\'{e}or\`{e}me
\ref{TheoHKRpourggsP}), on a un isomorphisme
\[
  HH_{\cdot}(\gs|\korps,\gs)\cong \Lambda^{.}HH_1(\gs|\korps,\gs)=:\Lambda^{\cdot}\Omega^1_{\gs|\korps}
\]
donc, d'apr\`{e}s la d\'{e}composition de Hodge mentionn\'{e}e ci-dessus, il vient
\[
   \Har_k(\gs|\korps,\gs)\cong \left\{
      \begin{array}{cc}
        \Omega^1_{\gs|\korps} & \mbox{ si }k=1 \\
          \{0\}  & \mbox{ si }k\geq 2
      \end{array}\right.
\]
ce qui ram\`{e}ne le calcul de $\Har_{\cdot}({\gsI}_+|\korps, \gs)$
\`{a} celui de $\Har_{\cdot+1}(\gs|{\gsI}_+, \gs)$ d'apr\`{e}s la suite exacte de Jacobi-Zariski.
Calculons donc $HH_{\cdot}(\gs|{\gsI}_+, \gs)$.

Il est bien connu que  l'homologie de Hochschild $HH_{\cdot}(B|A,M)$ d'une $A$-alg\`{e}bre
unitaire $B$ \`{a} valeurs dans un $B$-module $M$ est \'{e}gale \`{a} l'homologie du complexe de
Hochschild normalis\'{e} $\overline{C_k}(B|A,M)$
d\'{e}fini, en notant $B_{/A}$ le module quotient de $B$ par $A$,
pour tout $k\geq 0$, par $\overline{C_k}(B|A,M)=M \ot_{A} (B_{/A})^{\ot_A k}$ muni de la
diff\'{e}rentielle de Hochschild $\mathsf{b}$. Il nous suffit donc de calculer l'homologie de
$\overline{C_{\cdot}}(\gs_{/{\gsI}_+}, \gs)
=\gs\ot_{{\gsI}_+}(\gs_{/{{\gsI}_+}})\ot_{{\gsI}_+} \cdots \ot_{{\gsI}_+} (\gs_{/{{\gsI}_+}})$.

Comme ${\gsI}_+=\korps\oplus \gsI$ la projection $\gs\to \gsP$ induit un
isomorphisme $\gs_{/{{\gsI}_+}}=\gsP_{/ \korps}$. Puisque $\gsI$ est un id\'{e}al,  si on
munit $\gs_{/{\gsI}_+} \ot_{\korps} \gs_{/{\gsI}_+}$ de la structure de ${\gsI}_+$ module
induit par le morphisme d'anneaux ${\gsI}_+ \to \korps$, on voit que la projection
canonique $\gs_{/{\gsI}_+} \times \gs_{/{\gsI}_+} \to
\gs_{/{\gsI}_+} \ot_{\korps} \gs_{/{\gsI}_+}$ est ${\gsI}_+$-bilin\'{e}aire. Comme toute
application ${\gsI}_+$-bilin\'{e}aire est aussi $\korps$-bilin\'{e}aire, on en d\'{e}duit que toute
application ${\gsI}_+$-bilin\'{e}aire de $\gs_{/{\gsI}_+} \times \gs_{/{\gsI}_+}$ dans un module
$V$ se factorise au travers de $\gs_{/{\gsI}_+} \ot_{\korps} \gs_{/{\gsI}_+}$ et on obtient
un isomorphisme
\[
   \gsP_{/ \korps} \ot_{\korps} \gsP_{/ \korps}
   \cong \gs_{/{\gsI}_+} \ot_{\korps} \gs_{/{\gsI}_+}
   \cong \gs_{/{\gsI}_+} \ot_{{\gsI}_+} \gs_{/{\gsI}_+}.
\]
En it\'{e}rant ce raisonnement et en utilisant que la projection $\gs\to \gsP$ est un morphisme
d'alg\`{e}bres on obtient un isomorphisme de complexes
\begin{align*}
\gs\! \ot_{{\gsI}_+}\! \left(\gs_{/{{\gsI}_+}}\right)
          \ot_{{\gsI}_+}\! \cdots \ot_{{\gsI}_+}\!
             \left(\gs_{/{{\gsI}_+}}\right)
  &\cong
   (\gs/\gsI)\ot_{\korps} \left(\gsP_{/\korps}\right)
      \ot_{\korps}\cdots \ot_{\korps} \left(\gsP_{/\korps}\right)\cr
  &  \cong \gsP\ot_{\korps} \left(\gsP_{/\korps}\right)\ot_{\korps}
         \cdots \ot_{\korps} \left(\gsP_{/\korps}\right).
\end{align*}
Le dernier complexe n'est autre que le complexe de Hochschild normalis\'{e}
$\overline{C_{\cdot}}(\gsP|\korps,\gsP)$ de $\gsP$.
La version gradu\'{e}e du th\'{e}or\`{e}me de Hochschild-Kostant-Rosenberg
pour l'homologie de Hochschild continue gradu\'{e}e
de l'alg\`{e}bre gradu\'{e}e
$\gsP=\CinfK{C}\otimes_\korps \Lambda E''$ (voir le th\'{e}or\`{e}me \ref{TheoHKRpourggsP})
donne un isomorphisme
\[
   H_{\cdot}(\overline{C_\cdot}^{{\korps}}(\gsP,\gsP))\cong
      HH_{\cdot}(\gsP,\gsP)\cong \Lambda^{\cdot}\Omega^1_{\gsP}
\]
dont la composante de poids $1$ est $\Omega^1_{\gsP}$, voir
\cite[p.145, Prop. 4.5.13]{Lo92}, donc les groupes d'homologie
de Harrison s'annulent \`{a} partir de $k\geq 2$.
La suite exacte de Jacobi-Zariski se r\'{e}duit alors \`{a} la
 suite exacte de $\gs$-modules
\[
    \{0\}\ra \gs\ot_{{\gsI}_+}\Omega^1_{{\gsI}_+}\hookrightarrow
 \Omega^1_{\gs}\twoheadrightarrow \Omega^1_{\gsP}\ra \{0\}
\]
(o\`{u} la structure de
 $\gs$-module de $\Omega^1_{\gsP}$ est induite par la
 projection  $\gs\twoheadrightarrow \gsP$) qui donne l'isomorphisme
 cherch\'{e}
 $\Har_{\cdot}({\gsI}_+|\korps, \gs)\cong \gs\ot_{{\gsI}_+}\Omega^1_{{\gsI}_+}$.
\end{prooof}

\subsection{Acyclicit\'{e} du cocomplexe d'obstructions r\'{e}duit}

On rappelle
\[
  \begin{array}{lll}
   x'_\alpha:=x_\alpha & 1\leq \alpha \leq n-l&
                          \mathrm{coordonn\acute{e}es~le~long~de~}C, \\
  x''_\mu := x_{n-l+\mu} & 1\leq \mu \leq l &
   \mathrm{coordonn\acute{e}es~transversales~\grave{a}~}C.
  \end{array}
\]
et l'on pose
\[
  \begin{array}{lll}
   y'_\alpha & 1\leq \alpha \leq n-l&
                          \mathrm{base~de~l'espace~}E', \\
  y''_\mu  & 1\leq \mu \leq l &
   \mathrm{base~de~l'espace~}E''.
  \end{array}
\]
D'apr\`{e}s le th\'{e}or\`{e}me HKR \ref{TheoHKRpourggsP}, il vient que les modules
des diff\'{e}rentielles de K\"{a}hler continues de
$\gs = \CinfK{\mathbb R^{n}}\otimes \Lambda E$ et
de $\gsP  =  \CinfK{\mathbb R^{n-l}}\otimes \Lambda E''$ sont donn\'{e}s
par des modules libres $\Omega^1_\gs:=HH_1(\gs|\korps,\gs)=\gs\otimes(E^*\oplus E)$ (sur $\gs$)
et $\Omega^1_{\gsP}:=HH_1(\gsP|\korps,\gsP)=\gsP\otimes(E^{\prime *}\oplus E'')$ (sur $\gsP$).
On \'{e}crira les bases de $\Omega^1_{\gs}$ et de $\Omega^1_{\gsP}$ comme des diff\'{e}rentielles
des `coordonn\'{e}es' $(x,y)$ comme suit:
\[
   \mbox{Pour }\Omega^1_\gs:~~
   \mathbf{d}x'_\alpha=\drm x'_\alpha,\mathbf{d}x''_\mu=\drm x''_\mu,
   \mathbf{d}y'_\beta, \mathbf{d}y''_\nu
\]
quels que soient $1\leq \alpha,\beta\leq n-l$ et $1\leq \mu,\nu \leq l$,
et
\[
  \mbox{pour }\Omega^1_{\gsP}:~~
   \mathbf{d}x'_\alpha=\drm x'_\alpha,\mathbf{d}y''_\nu
\]
quels que soient $1\leq \alpha,\leq n-l$ et $1\leq \nu \leq l$.
Le morphisme de $\gs$-modules $\Omega^1_{\gs}\twoheadrightarrow
\Omega^1_{\gsP}$ est donc donn\'{e} par
\begin{eqnarray*}
\lefteqn{\sum_{\alpha}\xi_\alpha \dk x'_\alpha+\sum_{\mu}\xi''_\mu \dk x''_\mu
   +\sum_{\beta}\eta'_\beta \dk y'_\beta + \sum_{\nu}\eta''_\nu \dk
   y''_\nu} ~~~~~~~~~~~\\
     & & ~~~~~~~~~~~~~~\mapsto ~~~~~~~\sum_{\alpha}p_{\gsP}(\xi_\alpha) \dk x'_\alpha
             + \sum_{\nu}p_{\gsP}(\eta''_\nu) \dk y''_\nu
\end{eqnarray*}
o\`{u} $p_{\gsP}$ d\'{e}signe le morphisme d'alg\`{e}bres $\gs\ra\gsP$.
Puisqu'on peut identifier le $\gs$-module $\gs\ot_{{\gsI}_+}\Omega^1_{{\gsI}_+}$
avec $\Ker(\Omega^1_{\gs}\twoheadrightarrow \Omega^1_{\gsP})$ et
$\gsI=\Ker(\gs \to \gsP)$ est l'id\'{e}al engendr\'{e} par
$x''_\mu$ et $y'_\alpha$ en tant que $\gs$-module, alors
$\gs\ot_{{\gsI}_+}\Omega^1_{{\gsI}_+}$
est engendr\'{e} par les diff\'{e}rentielles
de K\"{a}hler
\begin{equation}\label{EqBaseDiffKaehlergsI}
  \dk x''_{\mu}, ~\dk y'_\alpha,~ x''_\nu \dk x'_\beta,~
    y'_\gamma \dk x'_{\gamma'}, ~x''_\rho\dk y''_\sigma,~
 y'_{\gamma''}\dk y''_\tau~~~
\end{equation}
quel que soient $1\leq
 \alpha,\beta,\gamma,\gamma',\gamma''\leq n-l$ et $1\leq
 \mu,\nu,\rho,\sigma,\tau\leq l$.
Remarquons que le crochet de Schouten d\'{e}finit la structure d'une alg\`{e}bre de Lie
gradu\'{e}e sur $\Omega^1_{\gs}$ d\'{e}finit par $[\dk \xi,\dk\eta]:=\dk
[\xi,\eta]_S$ et $\dk \xi .\eta :=[\dk \xi,\eta]:=[\xi,\eta]_S$ pour
$\xi,\eta\in\gs$.

%
%
%
%

\smallskip

\noindent Voici le r\'{e}sultat central de ce paragraphe~:

\begin{theorem}\label{TheoDeuxiemeReduction}
Si $X=\real^n$ et $C=\real^{n-l}$ avec $l\geq 2$,
la cohomologie du complexe
$\left({\Hom}_{{\gs}}(\ve_{{\gs}}^\cdot
\big(\gs \ot_{{\gsI}_+}\Omega_{{\gsI}_+}\big)\right.$
$\left.,\gsI),\delta^{1,1} \right)$
est concentr\'{e}e en degr\'{e} $0$.
\end{theorem}

\begin{prooof}
Montrons d'abord que tout morphisme de
${\Hom}_\gs(\ve_\gs^\cdot
\big(\gs \ot_{{\gsI}_+}\Omega_{{\gsI}_+}\big),\gsI)$ peut se relever en un
morphisme
de
$\Hom_\gs(\ve_{{\gs}}^\cdot \Omega_{{\gs}},\gs)$.
Reprenons les notations pr\'{e}d\'{e}dentes~:
le $\gs$-module
$\gs\ot_{{\gsI}_+}\Omega^1_{{\gsI}_+}$ est engendr\'{e} par les diff\'{e}rentielles
de K\"{a}hler (\ref{EqBaseDiffKaehlergsI}).

Soit $\varphi$ dans
${\Hom}_\gs(
\gs \ot_{{\gsI}_+}\Omega_{{\gsI}_+},\gsI)$.
Montrons que pour tout $1\leq \alpha \leq n-l$, il existe un unique $c_\alpha$
dans $\gs$ tel
que pour tous $1\leq \mu \leq l$ et pour tous $1\leq \beta \leq n-l$ on ait~:
\[
   \varphi\big(x''_\mu\dk x'_\alpha\big)=x''_\mu c_\alpha,
          \hbox{ et }
  \varphi\big(y'_\beta\dk x'_\alpha\big)=(-1)^{|\varphi|}y'_\beta c_\alpha.
\]
En effet, par $\gs$-lin\'{e}arit\'{e},
$x''_{l-1}\varphi(x''_l \dk x'_\alpha)=x''_l \varphi(x''_{l-1} \dk x'_\alpha)$
quel que soit $1\leq \alpha\leq n-l$. Par cons\'{e}quent,
$x''_{l-1}\varphi(x''_{l} \dk x'_\alpha)$ et donc
$\varphi(x''_{l} \dk x'_\alpha)$ s'annule sur l'hyperplan
$H_{l}:=\{x\in\mathbb R^n~|~x''_{l}=0\}$. L'id\'{e}al annulateur de $H_{l}$ est
\'{e}gal \`{a} $x''_{l}\CinfK{\mathbb R}$ d'apr\`{e}s le lemme d'Hadamard: si $g\in
\CinfK{\mathbb R^n}$ s'annule sur $H_l$, alors
\[
  g(x)=\int_{0}^1\frac{\partial g}{\partial
  x_{l}}(x_1,\ldots,x_{n-1},tx_n)dt~~x_n
\]
avec $x_n=x''_l$.
Donc il existe un unique $c_\alpha$ dans $\gs$ tel que
$\varphi(x''_l \dk x'_\alpha)=x''_l c_\alpha$
et $\varphi(x''_{l-1} \dk x'_\alpha)=x''_{l-1}c_\alpha$.
Soit $1\leq \mu\leq l-2$, on a encore
$x''_l \varphi(x''_{\mu} \dk x'_\alpha)
=x''_{\mu}\varphi(x''_l \dk x'_\alpha)=x''_\mu x''_l c_\alpha$
et donc $\varphi(x''_{\mu} \dk x'_\alpha)=x''_\mu c_\alpha$. De m\^{e}me, pour
$1\leq \beta \leq n-l$,
$x''_l \varphi(y'_{\beta} \dk x'_\alpha)
=(-1)^{|\varphi|}y'_{\beta}\varphi(x''_l \dk x'_\alpha)
=(-1)^{|\varphi|}y'_\beta x''_l c_\alpha$
et donc $\varphi(y'_{\beta} \dk x'_\alpha)=(-1)^{|\varphi|}y'_\beta c_\alpha$.\\
Il est clair que l'on
peut montrer de la m\^{e}me mani\`{e}re le r\'{e}sultat analogue que pour tout
$1\leq \mu \leq l$, il existe un unique $\hat{c}_\mu$
dans $\gs$ tel
que pour tous $1\leq \nu \leq l$ et pour tous $1\leq \beta \leq n-l$ on ait~:
\[
   \varphi\big(x''_\nu\dk y''_\mu\big)=x''_\nu \hat{c}_\mu,
          \hbox{ et }
  \varphi\big(y'_\beta\dk y''_\mu\big)=(-1)^{|\varphi|}y'_\beta \hat{c}_\mu.
\]
Gr\^{a}ce au fait que $\Omega^1_{{\gs}}$ est un $\gs$-module libre, il s'ensuit que la
prescription $\tilde{\varphi}(\dk x'_\alpha):=c_\alpha$
et $\tilde{\varphi}(\dk y''_\mu):=\hat{c}_\mu$ avec
$\tilde{\varphi}(\dk x''_\mu):=\varphi(\dk x''_\mu)$ et
$\tilde{\varphi}(\dk y'_\beta):=\varphi(\dk y'_\beta)$ d\'{e}finit un
$\gs$-homomorphisme $\tilde{\varphi}:\Omega_{\gs}\ra\gs$ qui co\"{\i}ncide
avec $\varphi$ sur les g\'{e}n\'{e}rateurs de
$\gs\otimes_{{\gsI}_+}\Omega_{{\gsI}_+}$. \\
On montre ensuite, par r\'{e}currence, que
tout morphisme de
${\Hom}_\gs(\ve_\gs^\cdot
\big(\gs \ot_{{\gsI}_+}\Omega_{{\gsI}_+}\big),\gsI)$ peut se prolonger en un
unique morphisme de
$\Hom_\gs(\ve_{{\gs}}^\cdot \Omega_{{\gs}},\gs)$.

Consid\'{e}rons maintenant un morphisme  $\varphi$, $k$-cocycle dans le complexe
$\left({\Hom}_{{\gs}}(\ve_{{\gs}}^\cdot
\big(\gs \ot_{{\gsI}_+}\right.$
$\left.\Omega_{{\gsI}_+}\big),\gsI),\delta^{1,1} \right)$ avec $k\geq 2$.
Relevons $\varphi$ en un morphisme
$\tilde{\varphi}$ de $\Hom_\gs(\ve_{{\gs}}^\cdot \Omega_{{\gs}},\gs)$.
Il est clair que $\tilde{\varphi}$ est aussi un cocycle dans ce dernier complexe.
D'apr\`{e}s \cite{Tam98p} (voir aussi \cite{GH03}), le complexe
$\Hom_\gs(\ve_{{\gs}}^\cdot \Omega_{{\gs}},\gs)$ est acyclique et donc $\tilde{\varphi}$
est un cobord dans ce complexe~: $\tilde{\varphi}= \delta^{1,1} (\xi)$.
Le morphisme $\xi$ est compl\`{e}tement d\'{e}fini sur les g\'{e}n\'{e}rateurs
de $\Omega_{{\gs}}$~: $\dk x_1, \dots, \dk x_n, \dk y_1, \dots ,$ $ \dk y_n$.
D\'{e}finissons le morphisme $\xi_{\gsI}$ sur ces m\^{e}mes g\'{e}n\'{e}rateurs du module
libre
($A_1,\ldots,A_k \in \{\dk x_1, \dots, \dk x_n, \dk y_1, \dots , \dk y_n\}$)~:
\begin{align*}
\xi_{\gsI}(A_1 \vea \cdots \vea A_k)&
        :=p_{\gsI}\left(\xi(A_1 \vea \cdots \vea A_k) \right)
\hbox{ si }\\
&\hskip2cm \forall i,~\big(A_i= \dk y'_\alpha~\hbox{ ou } A_i=\dk x''_\mu \big)\\
\xi_{\gsI}(A_1 \vea \cdots \vea A_k)&
        :=\xi(A_1 \vea \cdots \vea A_k)\hbox{ sinon }
\end{align*}
o\`{u} $p_{\gsI}$ (resp. $p_{\gsP}$)
est la projection $\gs \to \gsI$  (resp. $\gs \to \gsP$)
parall\`{e}lement \`{a} $\gsP$, consid\'{e}r\'{e}e comme sous-alg\`{e}bre ab\'{e}lienne de $(\gs[1],[~,~]_S)$
(resp. $\gsI$). En utilisant les g\'{e}n\'{e}rateurs du sous-module
$\gs\ot_{{\gsI}_+}\Omega^1_{{\gsI}_+}$, on voit que la restriction de
$\xi_{\gsI}$ \`{a} $\ve_\gs^\cdot
\big(\gs \ot_{{\gsI}_+}\Omega_{{\gsI}_+}\big)$ prend ses valeurs dans
$\gsI$.
On a, pour $A_1,...A_k \in \{\dk x_1, \dots, \dk x_n, \dk y_1,
 \dots , \dk y_n\}$~ la formule gradu\'{e}e de Chevalley-Eilenberg pour la
 codiff\'{e}rentielle $\delta^{1,1}$:
\begin{eqnarray}
\label{cobordequation}
 \lefteqn{\tilde{\varphi}(A_0 \vea \cdots \vea A_k) =
         (\delta^{1,1}\xi)(A_0 \vea \cdots \vea A_k)} \nonumber \\
    & &  = \sum_{i=0}^k\pm(-1)^i
             A_i.\big(\xi(A_0 \vea \cdots \vea A_{i-1}\vea A_{i+1}\vea\cdots
             \vea A_k\big)
 \end{eqnarray}
car les crochets de Lie $[A_i,A_j]$ s'annulent.
Montrons que cette \'{e}quation est encore vraie
si l'on remplace $\xi$ par $\xi_{\gsI}$~:

1. Si au moins deux $A_i$ sont \'{e}gaux \`{a} un
$\dk y''_\mu$ ou \`{a} un  $\dk x'_\alpha$, ceci est
imm\'{e}diat par d\'{e}finition de $\xi_{\gsI}$.

2. Si tous les $A_i$ sont \'{e}gaux \`{a}
un $\dk y'_\alpha$ ou \`{a} un  $\dk x''_\mu$, ceci est
encore vrai car
$\tilde{\varphi}(A_1 \vea \cdots \vea A_n)=
\varphi(A_1 \vea \cdots \vea A_n)$ est dans $\gsI$, le crochet
$[y'_\alpha,\gsP]\subset \gsP $
et $[x''_\mu,\gsP]\subset\gsP $.

3. Si maintenant un seul des $A_i$, notons le $B$, est \'{e}gal \`{a} un
$\dk y''_\mu$ ou \`{a} un  $\dk x'_\alpha$, on peut encore
une fois
remplacer $\xi$ par $\xi_{\gsI}$ car
$B.\big(\xi(\cdots)\big)=B.\big(\xi_{\gsI}(\cdots)\big)$
car $B\in \dk \gsP$ et $\gsP$ est une sous-alg\'{e}bre ab\'{e}lienne de
$\gs$ par rapport \`{a} $[~,~]_S$, voir la proposition \ref{PgspEmbed}.

En conclusion, on a trouv\'{e} un $\xi_{\gsI}$ qui se restreint bien en un
morphisme de ${\Hom}_\gs(\ve_\gs^\cdot
\big(\gs \ot_{{\gsI}_+}\Omega_{{\gsI}_+}\big),\gsI)$ tel que $\varphi =
\delta^{1,1}(\xi_{\gsI})$.
\end{prooof}

Dans le cas o\`{u} $X=\real^1$ et $C=\real^{0}$, un
calcul simple nous donne~:
\begin{prop}
La cohomologie du complexe
$\left({\Hom}_{{\gs}}(\ve_{{\gs}}^\cdot
\big(\gs \ot_{{\gsI}_+}\Omega_{{\gsI}_+}\big),\right.$ $\left.\gsI),\delta^{1,1} \right)$
est concentr\'{e} en degr\'{e} $0$.
\end{prop}

\medskip

\begin{appendix}
\section{Appendice}
 \label{SecAppendice}

Dans ce paragraphe, on rappellera --sans utiliser les op\'{e}rades--
quelques notions autour des structures
$A_\infty$, $C_\infty$, $L_\infty$ et $G_\infty$ --qui se trouvent
dans la litt\'{e}rature de fa\c{c}on dispers\'{e}e-- pour fixer nos notations et
conventions de signe, voir \'{e}galement \cite{Mac75}, \cite{LS93}, et pour
un cadre op\'{e}radique \cite{GH03} et \cite{MSS02}. On fixe un corps $\korps$
de caract\'{e}ristique $0$.

\subsection{Quelques cog\`{e}bres colibres pour des espaces gradu\'{e}s}
 \label{AppSubSecGradues}

On rappelle d'abord la {\em cat\'{e}gorie des $\korps$-espaces vectoriels gradu\'{e}s}:
les objets sont des $\korps$-espaces vectoriels gradu\'{e}s
$\hs=\oplus_{k\in\mathbb Z}\hs^k$ (o\`{u} `gradu\'{e}' veut toujours dire $\mathbb Z$-gradu\'{e}).
On \'{e}crira
$|x|\in\mathbb Z$
pour le degr\'{e} d'un \'{e}l\'{e}ment homog\`{e}ne $x\in\hs$. Par la suite, on n'utilisera que des \'{e}l\'{e}ments
homog\`{e}nes si rien d'autre n'est dit.
Pour un autre
$\korps$-espace gradu\'{e} $\hat{\hs}=\oplus_{k\in\mathbb Z}\hat{\hs}^k$ et un
entier $j$
on note $\Hom(\hs,\hat{\hs})^j$
l'espace vectoriel de toutes les applications lin\'{e}aires $\hs\ra
\hat{\hs}$ de degr\'{e} $j$
(i.e $f(\hs^k)\subset \hat{\hs}^{k+l}$). L'espace $\Hom(\hs,\hat{\hs})$ est
d\'{e}fini
par l'espace gradu\'{e} $\oplus_{j\in\mathbb Z}\Hom(\hs,\hat{\hs})^j$.
De plus, le produit tensoriel $\hs\otimes \hat{\hs}$ est gradu\'{e} par
$\big(\hs\otimes \hat{\hs}\big)^k=\oplus_{j\in\mathbb Z}\hs^j\otimes \hat{\hs}^{k-j}$,
voir e.g. \cite[p.175]{Mac75}. Le produit tensoriel de deux morphismes
$\phi:\hs\ra\hs'$ et $\psi:\hat{\hs}\ra\hat{\hs}'$ est d\'{e}fini en
appliquant la {\em r\`{e}gle de signe de Koszul}
\begin{equation}\label{EqRegleKoszul}
   \big(\phi\otimes \psi\big)(a\otimes b)=(-1)^{|\psi|~|a|}\phi(a)\otimes \psi(b)
\end{equation}
pour $a\in \hs^{|a|}$ et $\psi$ de degr\'{e} $|\psi|$, voir e.g.
\cite[p.176]{Mac75}. On rappelle \'{e}galement la {\em transposition gradu\'{e}e}
$\tau:\hs\otimes\hat{\hs}\ra\hat{\hs}\otimes\hs$ par
$\tau(x\otimes y):= (-1)^{|x|~|y|}y\otimes x$. Les deux r\`{e}gles
pr\'{e}c\'{e}dentes d\'{e}termineront tous les signes qui appara\^{\i}tront.\\
 Pour un entier $j$, l'espace
gradu\'{e} {\em d\'{e}cal\'{e}} $\hs[j]$ et d\'{e}fini par $\hs[j]^k:=\hs^{k+j}$. ,
L'application identique
$\hs\ra\hs$ d\'{e}finit une application $s^n:\hs[j]\ra\hs[j-n]$ pour tout $n\in\mathbb Z$
qui est de degr\'{e}
$n$ car $s^n(\hs[j]^k)=\hs^{j+k}=\hs[j-n]^{k+n}$. On regarde $s^n$ comme la $n$\`{e}me
it\'{e}ration de la {\em suspension} $s:=s^1:\hs[j]\ra\hs[j-1]$. La suspension
sera `visible' pour des {\em applications multilin\'{e}aires d\'{e}cal\'{e}es}: soit
$\phi:\hs^{\otimes k}\ra \hat{\hs}^{\otimes l}$ une application lin\'{e}aire
de degr\'{e} $i$. On d\'{e}finit sa d\'{e}cal\'{e}e
$\phi[j]:\hs[j]^{\otimes k}\ra \hat{\hs}[j]^{\otimes l}$ par
$\phi[j]:= (s^{\otimes l})^{-j}\circ \phi
\circ (s^{\otimes k})^{j}$. Son degr\'{e} est \'{e}gal \`{a} $j(k-l)+i$ et on a \'{e}videmment
$(\phi[j])[j']=\phi[j+j']$. On remarque
que $(s^{\otimes k})^{j}
=(-1)^{\frac{k(k-1)}{2}\frac{j(j-1)}{2}}(s^j)^{\otimes k}$.
Pour calculer la d\'{e}cal\'{e}e, on \'{e}crit d'abord pour
$\xi:=x_1\otimes\cdots\otimes x_k\in\hs^{\otimes k}$ la valeur $\phi(\xi)$
avec la convention de Sweedler comme
$\sum\phi_{(1)}(\xi)\otimes \cdots \otimes\phi_{(l)}(\xi)$ avec
$\phi_{(i)}(\xi)\in\hat{\hs}$. D'apr\`{e}s la r\`{e}gle de signe (\ref{EqRegleKoszul}),
la valeur de la d\'{e}cal\'{e}e $\phi[j]$ sur
$\eta:=y_1\otimes\cdots\otimes y_k\in\hs[j]^{\otimes k}$ se calcule
de la mani\`{e}re suivante avec $\tilde{\eta}:=s^j(y_1)\otimes\cdots\otimes
s^j(y_k)$:
\bea\label{EqApplMultDecalees}
   \lefteqn{\phi[j](y_1\otimes\cdots\otimes y_k)=} \nonumber \\
     & & (-1)^{\frac{k(k-1)}{2}\frac{j(j-1)}{2}+\frac{l(l-1)}{2}\frac{-j(-j-1)}{2}}
         (-1)^{j\big((k-1)|y_1|+(k-2)|y_2|+\cdots+(k-(k-1))|y_{k-1}|\big)}
         \nonumber \\
     & & ~~~\sum (-1)^{j\big((l-1)|\phi_{(1)}(\tilde{\eta})|+
                           (l-2)|\phi_{(2)}(\tilde{\eta})|+\cdots+
                           (l-(l-1))|\phi_{(l-1)}(\tilde{\eta})|\big)}
                           \nonumber \\
     & & ~~~~~~~~~~\phi_{(1)}(y_1\otimes\cdots \otimes y_k)\otimes\cdots\otimes
         \phi_{(l)}(y_1\otimes\cdots \otimes y_k).
\eea

{\em L'alg\`{e}bre libre} sur $\hs$,
$\hs_{\mathrm{alg}}^\otimes:=\oplus_{k\in\mathbb N}\hs^{\otimes k}$, est
une $\korps$-alg\`{e}bre associative gradu\'{e}e avec \'{e}l\'{e}ment neutre $1$, voir e.g. \cite[p.179]{Mac75}.
Pour \'{e}viter des confusions,
on n'utilise pas le symbole $\otimes$ pour la multiplication
$\mu=\mu_{\hs^\otimes_{\mathrm{alg}}}$
dans l'alg\`{e}bre
libre. De plus, c'est
une big\`{e}bre gradu\'{e}e: soit $\hs^{\otimes +}:=\oplus_{k>0}\hs^{\otimes k}$ l'id\'{e}al
d'augmentation. La co-unit\'{e}
$\epsilon=\epsilon_{\hs^\otimes}:\hs^\otimes_{\mathrm{alg}}\ra\korps$ est d\'{e}finie
par la condition que $Ker \epsilon:=\hs^{\otimes +}$ et $\epsilon(1):=1$;
et la {\em comultiplication shuffle} gradu\'{e}e
$\Delta_{sh}$ est l'homomorphisme d'alg\`{e}bres associatives $\hs_{\mathrm{alg}}^\otimes\ra
\hs_{\mathrm{alg}}^\otimes \otimes \hs_{\mathrm{alg}}^\otimes$
induit par sa valeur
$\Delta_{sh}(x):=x\otimes 1
+ 1\otimes x$ sur les g\'{e}n\'{e}rateurs $x\in\hs$. Puisque la multiplication $\mu^{[2]}$
sur $\hs_{\mathrm{alg}}^\otimes \otimes \hs_{\mathrm{alg}}^\otimes$ est donn\'{e}e par
$(\mu\otimes \mu)\circ (\id \otimes \tau \otimes \id)$ avec la
transposition gradu\'{e}e, il y a donc des signes dans les
formules pour
$\Delta_{sh}$, par exemple $\Delta_{sh}(xy)=xy\otimes 1 +x\otimes y
+(-1)^{|x|~|y|}y\otimes x +1\otimes xy$ pour $x\in \hs^{|x|}$ et
$y\in\hs^{|y|}$. Cette comultiplication est cocommutative gradu\'{e}e (i.e.
$\tau\Delta_{sh}=\Delta_{sh}$) et de
degr\'{e} $0$.\\
On peut dualiser cette structure de big\`{e}bre
gradu\'{e}e pour obtenir sur l'espace gradu\'{e} $\hs_{\mathrm{alg}}^\otimes$ une
deuxi\`{e}me structure de big\`{e}bre
gradu\'{e}e, plus importante pour nos fins,
donn\'{e}e par la comultiplication (non cocommutative gradu\'{e}e) de
{\em d\'{e}concat\'{e}nation}, $\Delta=\Delta_{\hs^\otimes}$ (qui dualise la multiplication libre et
est d\'{e}finie par $\Delta(x_1\cdots x_k)
=1\otimes x_1\cdots x_k +\sum_{r=2}^{k}x_1\cdots x_{r-1}\otimes x_{r}\cdots x_{k}+
x_1\cdots x_k\otimes 1$) et la
{\em multiplication shuffle} $\mu_{sh}$ gradu\'{e}e commutative, parfois
\'{e}crite $\bullet$,
(qui dualise la comultiplication
$\Delta_{sh}$). Pour une formule explicite de $\mu_{sh}$,
voir (\ref{EqDefMultShuffleAvecSign}).
Les deux op\'{e}rations $\Delta$ et $\mu_{sh}$ sont de degr\'{e}
$0$. Ici $1$ reste l'\'{e}l\'{e}ment neutre et $\epsilon$ la co-unit\'{e}. On note cette big\`{e}bre gradu\'{e}e
par $\hs^\otimes$ et la projection canonique
sur $\hs^{\otimes 1}=\hs$ par $\mathrm{pr}_\hs$.\\
En g\'{e}n\'{e}ral, une cog\`{e}bre gradu\'{e}e $(C,\epsilon_C,\Delta_C)$ est dite
{\em augment\'{e}e} lorsqu'il existe un unique \'{e}l\'{e}ment de genre groupe, not\'{e}
$1$, alors $C=\korps 1 \oplus \mathrm{Ker}\epsilon_C$. Le sous-espace
$C^+:=\mathrm{Ker}\epsilon_C$ est isomorphe \`{a} la cog\`{e}bre gradu\'{e}e quotient
$C/\korps 1$ sans co-unit\'{e} ($\korps 1$ est une sous-cog\`{e}bre, donc un co-id\'{e}al de
$\hs^\otimes$). Une cog\`{e}bre gradu\'{e}e sans co-unit\'{e} est dite
{\em nilpotente} (voir \cite[166-168]{MSS02}) lorsque pour tout \'{e}l\'{e}ment $x$ il existe un
entier positif
$N$ tel que la $N$\`{e}me it\'{e}ration de la comultipliciation sur $x$ s'annule.
Les cog\`{e}bres gradu\'{e}es augment\'{e}es telles que $C^+$ (vu comme cog\`{e}bre quotient)
est nilpotente forment une
sous-cat\'{e}gorie $\mathcal{C}_{AN}$ de la cat\'{e}gorie des cog\`{e}bres gradu\'{e}es.
La cat\'{e}gorie $\mathcal{C}_{AN}$ est ferm\'{e}e pour les produits tensoriels et contient
$\hs^\otimes$ et toutes les cog\`{e}bres consid\'{e}r\'{e}es dans cet appendice.
La cog\`{e}bre $\hs^\otimes$ est
{\em colibre} dans $\mathcal{C}_{AN}$, i.e. pour toute cog\`{e}bre $C$
dans $\mathcal{C}_{AN}$ et toute application lin\'{e}aire $\phi:C\ra \hs$ de degr\'{e}
$0$ s'annulant sur $1$
il existe un unique morphisme de cog\`{e}bres
gradu\'{e}es
$\overline{\phi}:C\ra \hs^\otimes$ (i.e.
$\overline{\phi}\otimes\overline{\phi}\circ \Delta_C=\Delta\circ
\overline{\phi}$) tel que $\mathrm{pr}_\hs\circ \overline{\phi}=\phi$:
\begin{equation}\label{EqDiagrammeColibre}
    \begin{array}{ccccc}
      \hs^\otimes &  & \stackrel{\overline{\phi}}{\longleftarrow} & &  C
      \\
                  &\stackrel{\mathrm{pr}_\hs}{\searrow} &  & \stackrel{\phi}{\swarrow}  &  \\
                  &              &  \hs        &                    &
    \end{array}
\end{equation}
On dit que $\overline{\phi}$ est co-induit par $\phi$. Evidemment, deux
morphismes de cog\`{e}bres $\Phi,\Psi:\hs^\otimes\leftarrow C$ co\"{\i}ncident
si et seulement si leurs composantes $\mathrm{pr}_\hs\circ \Phi$ et
$\mathrm{pr}_\hs\circ
\Psi$ co\"{\i}ncident. On calcule par exemple que la multiplication $\mu_{\mathrm{sh}}$
est co-induite par l'application $\hs^\otimes\otimes\hs^\otimes\ra \hs$ donn\'{e}e par
$\mathrm{pr}_\hs\otimes \epsilon +\epsilon\otimes \mathrm{pr}_\hs$.\\
L'alg\`{e}bre colibre dans la cat\'{e}gorie de toutes les cog\`{e}bres coassociatives
gradu\'{e}es {\em n'est pas} donn\'{e}e par $\hs^\otimes$, mais par un espace
`plus grand' entre $\hs^\otimes$ et son compl\'{e}t\'{e} par rapport \`{a} la
filtration donn\'{e}e par degr\'{e} tensoriel, voir \cite[pp.124-135]{Swe69} et
\cite[166-168]{MSS02}.\\
Pour deux applications lin\'{e}aires $\psi_1,\psi_2$ d'une cog\`{e}bre coassociative gradu\'{e}e
$(C,\Delta_C)$ dans une alg\`{e}bre associative gradu\'{e}e $(A,\mu_A)$
on rappelle
la {\em convolution
$\psi_1\star \psi_2$ de $\psi_1$ et de $\psi_2$ par rapport \`{a} $\mu_A$ et $\Delta_C$} donn\'{e}e
par $\psi_1\star \psi_2:=\mu_A\circ (\psi_1\otimes \psi_2 )\circ \Delta_C$, qui
est une multiplication associative gradu\'{e}e dans
$\mathrm{Hom}(C,A)$. Le morphisme de cog\`{e}bres
$\overline{\phi}$ co\"{\i}nduit par $\phi:C^+\ra \hs$ se calcule comme
s\'{e}rie g\'{e}om\'{e}trique $\overline{\phi}=\sum_{r=0}^\infty
\phi^{\star r}$ avec $\phi^{\star 0}:=1\epsilon_C$ \`{a} l'aide de la
convolution $\star$ par rapport \`{a} la multiplication libre $\mu$ et
$\Delta_C$.\\
Soit $d:C^+\ra \hs$ une application lin\'{e}aire
de degr\'{e} $j\in\mathbb Z$. On peut montrer qu'il existe une unique {\em cod\'{e}rivation
gradu\'{e}e de degr\'{e} $j$ le long de $\overline{\phi}$}, $\overline{d}:C\ra \hs^\otimes$,
(i.e. $\Delta\circ \overline{d}
=(\overline{d}\otimes \overline{\phi}+\overline{\phi}\otimes
\overline{d})\circ \Delta_C$) avec $\mathrm{pr}_1\circ \overline{d}=d$.
Cette cod\'{e}rivation
$\overline{d}$, dite co-induite par $d$, se calcule comme
$\overline{\phi}\star d\star\overline{\phi}$.
En particulier, pour $d_1,d_2:\hs^{\otimes +}\ra\hs$ et $\overline{\phi}=\id_{\hs^{\otimes}}$
on note
\begin{equation}\label{EqMultGerstenhaber}
   d_1\circ_G d_2 := d_1\circ \overline{d_2}
   =\sum_{i,j=0}^{\infty} d_1\circ\big(\id_\hs^{\otimes i}\otimes d_2 \otimes
               \id_\hs^{\otimes j}\big)~:~~
   \hs^{\otimes +}\ra\hs
\end{equation}
la {\em multiplication de Gerstenhaber} gradu\'{e}e. Il s'en suit que le {\em commutateur
gradu\'{e}} $[\overline{d_1},\overline{d_2}]=\overline{d_1}\circ\overline{d_2}-
(-1)^{|d_1|~|d_2|}\overline{d_2}\circ\overline{d_1}$ des deux cod\'{e}rivations le long de
$\id_{\hs^\otimes}$ est une cod\'{e}rivation de $\hs^\otimes$ le long de
$\id_{\hs^\otimes}$,
alors $[\overline{d_1},\overline{d_2}]$ est co-induit par sa composante
$\mathrm{pr}_\hs[\overline{d_1},\overline{d_2}]
=d_1\circ_G d_2-(-1)^{|d_1|~|d_2|}d_2\circ_G d_1=:[d_1,d_2]_G$,
le {\em crochet de Gerstenhaber gradu\'{e}} de $d_1$ et
$d_2$. On en d\'{e}duit l'identit\'{e} de Gerstenhaber suivante
\begin{equation}\label{EqIdGerstenhaber}
   (d\circ_G d_1)\circ_G d_2
     -(-1)^{|d_1|~|d_2|}(d\circ_G d_2)\circ_G d_1
     =
       d\circ_G[d_1,d_2]_G.
\end{equation}
Ainsi $\Hom(\hs^{\otimes +},\hs)$ muni du crochet de Gerstenhaber est une alg\`{e}bre de Lie
gradu\'{e}e. Les structures
$\circ_G$ et $[~,~]_G$ \'{e}taient d\'{e}finies dans \cite{Ger63} pour la situation
$\hs=V[1]$ o\`{u} l'espace gradu\'{e} \'{e}tait $V=V^0$. Soit $m:\hs\otimes\hs\ra\hs$
une multiplication associative gradu\'{e}e (de degr\'{e} $0$). Alors
l'application d\'{e}cal\'{e}e $d:=m[1]:\hs[1]\otimes \hs[1]\ra \hs[1]$ est de degr\'{e} $1$,
 et l'associativit\'{e} de $m$ est \'{e}quivalente
\`{a} $d\circ_G d=0$.
L'application $\mathsf{b}:\Hom(\hs[1]^{\otimes +},\hs[1])\ra \Hom(\hs[1]^{\otimes +},\hs[1])$
d\'{e}finie par $\mathsf{b}\phi:=[\phi,d]_G$ co\"{\i}ncide avec {\em l'op\'{e}rateur
cobord de Hochschild gradu\'{e}}. Les formules (\ref{EqDefMultGerstDiff})
et (\ref{EqDefCrochGerstDiff}) de la premi\`{e}re partie sont des cas
particuliers $\hs=\hs^0=\CinfK{\real^n}$. On observe que $m[1]=-m$ dans ce
cas. Alternativement, on peut d\'{e}finir le cobord de Hochschild directement sur
$\Hom(\hs^{\otimes +},\hs)$ sans d\'{e}calage par la formule \'{e}quivalente et plus simple
$\mathsf{b}\phi=\big[\phi[1],m[1]\big]_G[-1]$, voir par exemple l'\'{e}quation
(\ref{EqDefCobordHochsGrad}).

La {\em big\`{e}bre sym\'{e}trique gradu\'{e}e} sur $\hs$,
$\mathsf{S}\hs=\oplus_{r\in\mathbb N}\mathsf{S}^r\hs$, est
d\'{e}finie par le quotient de l'alg\`{e}bre libre
$\hs_{\mathrm{alg}}^\otimes$
modulo l'id\'{e}al
bilat\`{e}re gradu\'{e} $I_+(\hs)$ engendr\'{e} par tous les \'{e}l\'{e}ments de la forme
$xy-(-1)^{|x|~|y|}yx$
quels que soient $x\in\hs^{|x|}$ et $y\in\hs^{|y|}$. Le quotient $\mathsf{S}\hs$ est une
alg\`{e}bre associative commutative gradu\'{e}e. De plus, la premi\`{e}re
comultiplication $\Delta_{sh}$ passe au quotient et d\'{e}finit
sur $\mathsf{S}\hs$ la structure d'une comultiplication cocommutative gradu\'{e}e,
not\'{e}e \'{e}galement $\Delta_{sh}$. L'espace $\mathsf{S} \hs$ est donc une big\`{e}bre
commutative cocommutative gradu\'{e}e. Il est connu que $\mathsf{S}\hs$ est {\em l'alg\`{e}bre
commutative gradu\'{e}e libre engendr\'{e}e par $\hs$}.

Pour un entier $n$, une permutation $\sigma$ du groupe sym\'{e}trique $S_n$ et
$\xi=x_1\cdots x_n\in\hs^{\otimes n}$, on note $\xi^\sigma :=
x_{\sigma(1)}\cdots x_{\sigma(n)}$ l'action \`{a} droite usuelle de $S_n$ sur
$\hs^{\otimes n}$. On d\'{e}finit la {\em signature gradu\'{e}e de $\sigma$ par
rapport \`{a} $\xi$} par
\begin{eqnarray}\label{EqDefSignatureGraduee}
   e(x_1\cdots x_n, \sigma)
        &:= & \prod_{1\leq i<j\leq n}
             \frac{\sigma(i)+(-1)^{ |x_{\sigma(i)}|~|x_{\sigma(j)}| }~\sigma(j) }
                                    {i+(-1)^{|x_i|~|x_j|}~j} \nonumber \\
        & = & \prod_{i<j~\mathrm{et}~\sigma(i)>\sigma(j)}
                 (-1)^{ |x_{\sigma(i)}|~|x_{\sigma(j)}| }
\end{eqnarray}
et on en d\'{e}duit l'action \`{a} droite gradu\'{e}e $\xi^{.\sigma}:=e(\xi,\sigma)\xi^\sigma$ de
$S_n$ sur
$\hs^{\otimes n}$ car $e(\xi,\sigma\tau)=e(\xi,\sigma)e(\xi^\sigma,\tau)$.
A l'aide de cette action, on peut donner une formule plus explicite de la
multiplication shuffle: on note $Sh(k,n-k)$ l'ensemble des {\em permutations de
battage ou shuffle}, i.e. pour lesquelles $\sigma(1)<\cdots <\sigma(k)$ et
$\sigma(k+1)<\ldots < \sigma(n)$:
\begin{equation}\label{EqDefMultShuffleAvecSign}
  (x_1\cdots x_k)\bullet(x_{k+1}\cdots x_n)=
    \sum_{\sigma^{-1}\in Sh(k,n-k)}
            e(x_1\cdots x_n,\sigma) x_{\sigma(1)}\cdots x_{\sigma(n)}
\end{equation}
Pour tout entier positif $n$ soit
$\mathcal{S}^{(n)}:\hs^{\otimes n}\ra\hs^{\otimes n}$ la
sym\'{e}trisation $\xi\mapsto \sum_{\sigma\in S_n}e(\xi,\sigma)\xi^\sigma$. On d\'{e}finit
$\mathcal{S}:\hs_{\mathrm{alg}}^\otimes\ra \hs^\otimes$ par
$\mathcal{S}:=\sum_{n=0}^\infty \mathcal{S}^{(n)}$. On voit que
$\mathcal{S}$ est un morphisme de la big\`{e}bre $\hs_{\mathrm{alg}}^\otimes$
dans la big\`{e}bre (duale) $\hs^\otimes$ et l'image de $\mathcal{S}$ est
la sous-big\`{e}bre commutative et cocommutative $\mathsf{S}\hs$. L'application
$\mathcal{S}$ est
co-induite par l'application de la cog\`{e}bre cocommutative gradu\'{e}e
$(\hs_{\mathrm{alg}}^\otimes,\Delta_{\mathrm{sh}})$ dans $\hs$ donn\'{e}e
par $\mathrm{pr}_\hs$.
Encore une fois, la cog\`{e}bre cocommutative gradu\'{e}e $\mathsf{S}\hs$ est colibre
dans la cat\'{e}gorie $\mathcal{CS}_{AN}$ des cog\`{e}bres cocommutatives augment\'{e}es gradu\'{e}es \`{a}
$C^+$ nilpotent, et la co-induction des morphismes et cod\'{e}rivations
dans le sens du diagramme (\ref{EqDiagrammeColibre}) se fait comme pour
$\hs^\otimes$, o\`{u} les applications co-induites pour la cog\`{e}bre $\hs^\otimes$,
$\overline{\phi}$ et $\overline{d}$, prennent leurs valeurs
automatiquement dans la sous-cog\`{e}bre $\mathsf{S}\hs$. On peut \'{e}galement utiliser la
multiplication shuffle $\mu_{\mathrm{sh}}$ de $\mathsf{S}\hs$ au lieu de $\mu$ de $\hs^\otimes$
pour une deuxi\`{e}me convolution
$\tilde{\star}$ pour obtenir $\overline{\phi}=e^{\tilde{\star}\phi}$ au
lieu de la s\'{e}rie g\'{e}om\'{e}trique, et $\overline{d}=e^{\tilde{\star}\phi}\tilde{\star}d$
au lieu de $\overline{\phi}\star d\star\overline{\phi}$. Exactement de la
m\^{e}me fa\c{c}on qu'en (\ref{EqMultGerstenhaber}) pour $\hs^\otimes$ on d\'{e}finit la
{\em multiplication de
Nijenhuis-Richardson} $\circ_{NR}$ sur l'espace $\Hom(\mathsf{S}^+\hs,\hs)$
qui satisfait la m\^{e}me identit\'{e} (\ref{EqIdGerstenhaber}) avec $\circ_G$ remplac\'{e} par
$\circ_{NR}$. L'espace
$\Hom(\mathsf{S}^+\hs,\hs)$
est une alg\`{e}bre de Lie gradu\'{e}e par le crochet de Nijenhuis-Richardson
$[d_1,d_2]_{NR}:=d_1\circ_{NR} d_2-(-1)^{|d_1|~|d_2|}d_2\circ_{NR}d_1$.
Les structures
$\circ_{NR}$ et $[~,~]_{NR}$ \'{e}taient d\'{e}finies dans \cite{NR67} pour la situation
$\hs=V[1]$ o\`{u} l'espace gradu\'{e} \'{e}tait $V=V^0$.\\
Si l'on regarde dans $\hs^\otimes_{\mathrm{alg}}$ l'id\'{e}al bilat\`{e}re gradu\'{e}
$I_-(\hs)$ engendr\'{e} par tous les \'{e}l\'{e}ments de la forme
$xy+(-1)^{|x|~|y|}yx$
pour $x\in\hs^{|x|}$ et $y\in\hs^{|y|}$, on obtient comme
quotient {\em l'alg\`{e}bre de Grassmann} $\ve \hs :=\oplus_{r\in\mathbb N}
\ve^r\hs$: ceci est une alg\`{e}bre associative $\mathbb Z \times \mathbb Z$-gradu\'{e}e par
le degr\'{e} usuel et le degr\'{e} tensoriel qui est donc bigradu\'{e}e commutative (et non pas
gradu\'{e}e commutative). Soit $j$ un entier, $\tau_{\hs\otimes\hs}$ la transposition
dans $\hs\otimes\hs$ et $\tau_{\hs[j]\otimes\hs[j]}$ la transposition
dans $\hs[j]\otimes\hs[j]$. Gr\^{a}ce \`{a} l'identit\'{e}
\begin{equation}
   (\tau_{\hs\otimes\hs})[j]= (-1)^j\tau_{\hs[j]\otimes\hs[j]}
\end{equation}
ou $s^j\otimes s^j\circ\tau_{\hs[j]\otimes\hs[j]}
=(-1)^j\tau_{\hs\otimes\hs}\circ s^j\otimes s^j$,
on voit sans peine que si $j$ est impair, la bijection lin\'{e}aire
$\sum_{r=0}^\infty (s^{\otimes r})^j:
\hs[j]_{\mathrm{alg}}^\otimes \ra\hs_{\mathrm{alg}}^\otimes$
envoie l'id\'{e}al $I_+(\hs[j])$ sur l'id\'{e}al $I_-(\hs)$, et l'on obtient l'isomorphisme
lin\'{e}aire (pas d'alg\`{e}bres associatives)
\[
    \forall~r\in\mathbb N:~~  S^r\big(\hs[j]\big)\cong \big(\ve^r \hs\big)[jr]
         ~~~\mbox{pour $j$ impair}.
\]
Pour $\hs=\hs^0=V$ on constate que $\mathsf{S}\big(V[-1]\big)=\ve V$ comme alg\`{e}bres
commutatives gradu\'{e}es. En g\'{e}n\'{e}ral, une application lin\'{e}aire
$\phi:\hs^{\otimes k}\ra \hat{\hs}$ est dite {\em antisym\'{e}trique} lorsque
$\phi\circ(\id^{i}\otimes \tau\circ \otimes \id^{k-2-i})=-\phi$ quel que
soit $0\leq i\leq k-2$. D'apr\`{e}s ce qui pr\'{e}c\`{e}de, sa d\'{e}cal\'{e}e $\phi[j]$ est
sym\'{e}trique pour $j$ impair:
\[
  \phi\in \Hom(\ve \hs,\hat{\hs})~~ \Leftrightarrow~~
     \phi[j]\in \Hom (\mathsf{S}(\hs[j]),\hat{\hs}[j])
     \mbox{ pour }j\mbox{ impair}.
\]
En particulier, soit $(\mathfrak{l},[~,~])$ une {\em alg\`{e}bre de Lie gradu\'{e}e}, i.e. pour le
crochet $[~,~]\in \Hom(\ve^2\mathfrak{l},\mathfrak{l})$ de degr\'{e} $0$ on a {\em l'identit\'{e} de
Jacobi gradu\'{e}e}
$[~,~]([~,~]\otimes \id)(\id^{\otimes3}+\tau_{12}\tau_{23}+\tau_{23}\tau_{12})
=0$. Sur l'espace d\'{e}cal\'{e} $\mathfrak{l}[j]$ ($j$ impair),
le crochet d\'{e}cal\'{e} $[~,~][j]$ est de degr\'{e} $j$, se trouve dans
$\Hom(\mathsf{S}^2(\mathfrak{l}[j],\mathfrak{l[j]})$, est donc sym\'{e}trique
gradu\'{e} et satisfait
$[~,~][j]\circ([~,~][j]\otimes\id)\circ
(\id^{\otimes3}+\tau'_{12}\tau'_{23}+\tau'_{23}\tau'_{12})$ avec
$\tau':=\tau_{\mathfrak{l}[j]\otimes\mathfrak{l}[j]}$.
De plus, la d\'{e}cal\'{e}e $d=[~,~][1]$ du crochet de Lie
satisfait $d\circ_{NR} d=0$ et d\'{e}finit par $\phi\mapsto [\phi,d]_{NR}$ {\em l'op\'{e}rateur
cobord de Chevalley-Eilenberg} sur $\Hom(\mathsf{S}^+(\hs[1]),\hs[1])$.
Comme dans le cas d'une alg\`{e}bre associative, on peut alternativement
d\'{e}finir le cobord de Chevalley-Eilenberg directement sur
$\Hom(\ve\mathfrak{l},\mathfrak{l})$ par la formule isomorphe $\psi\mapsto
\big[\psi[1],[~,~][1]\big]_{NR}[-1]$.

L'id\'{e}al d'augmentation $\hs^{\otimes +}=\oplus_{k\geq 1}\hs^{\otimes k}$
est \'{e}galement
un id\'{e}al bilat\`{e}re de la multiplication shuffle, ainsi que son carr\'{e}
$\hs^{\otimes +}\bullet \hs^{\otimes +}\subset \hs^{\otimes +}$. L'espace gradu\'{e}
\[
     \underline{\hs^{\otimes}}:=\bigoplus_{k\geq 1}\underline{\hs^{\otimes k}}
        := \hs^{\otimes +}/(\hs^{\otimes +}\bullet \hs^{\otimes +})
\]
est muni de la structure d'une {\em cog\`{e}bre de Lie gradu\'{e}e}
$\delta:\underline{\hs^{\otimes}}\ra \underline{\hs^{\otimes}}\otimes
\underline{\hs^{\otimes}}$ de degr\'{e} $0$ (i.e. $\tau\delta=-\delta$ et
l'identit\'{e} de Jacobi gradu\'{e}e $(\id^{\otimes 3}+\tau_{12}\tau_{23}+\tau_{23}\tau_{12})
(\delta\otimes \id)~\delta=0$ avec $\tau_{12}:=\tau\otimes \id$ et
$\tau_{23}:=\id\otimes \tau$) provenant de
l'antisym\'{e}trisation gradu\'{e}e $\Delta-\tau\Delta$ de la comultiplication
$\Delta$ qui passe au quotient. Une cog\`{e}bre de Lie gradu\'{e}e
$(\mathfrak{c},\delta_{\mathfrak{c}})$ est dite {\em nilpotente} lorsque
pour tout $x\in\mathfrak{c}$ il existe un entier $N$ tel que tous les
$N$\`{e}mes it\'{e}r\'{e}s de $\delta_{\mathfrak{c}}$ s'annulent sur $x$. Dans la
cat\'{e}gorie $\mathcal{CL}_N$ de ces cog\`{e}bres de Lie, l'espace
$(\underline{\hs^{\otimes}},\delta)$ est
une cog\`{e}bre de Lie gradu\'{e}e {\em colibre} dans le sens du diagramme
\ref{EqDiagrammeColibre},
en rempla\c{c}ant cog\`{e}bres par cog\`{e}bres de Lie. De la m\^{e}me fa\c{c}on, on peut co-induire des
cod\'{e}rivations des cog\`{e}bres de Lie de degr\'{e} $j$ le long d'un morphisme de cog\`{e}bres de Lie
gradu\'{e}es o\`{u} l'on remplace $\Delta$ par $\delta$ dans les d\'{e}finitions.
Les applications lin\'{e}aires $d$ de $\Hom(\underline{\hs^\otimes},\hs)$ s'identifient
avec les applications lin\'{e}aires $d'$ de $\Hom(\hs^\otimes,\hs)$ s'annulant sur
$\hs^{\otimes +}\bullet \hs^{\otimes +}$
et sont connues comme des {\em cocha\^{\i}nes de Harrison gradu\'{e}es}. On montre que la
cod\'{e}rivation gradu\'{e}e $\overline{d'}$ est \'{e}galement une d\'{e}rivation gradu\'{e}e de la multiplication
shuffle $\bullet$, donc la multiplication de Gerstenhaber $d'_1\circ_G d'_2$ s'annule toujours
sur $\hs^{\otimes +}\bullet \hs^{\otimes +}$, ce qui prouve que $\circ_G$
pr\'{e}serve le sous-espace
et d\'{e}finit sur
$\Hom(\underline{\hs^\otimes},\hs)$ une multiplication $\circ_H$ satisfaisant
(\ref{EqIdGerstenhaber}) avec $\circ_G$ remplac\'{e} par $\circ_H$. L'espace
$\Hom(\underline{\hs^\otimes},\hs)$ muni
du crochet $[~,~]_H$ (l'antisym\'{e}tris\'{e} gradu\'{e} de $\circ_H$) est donc une alg\`{e}bre de Lie
gradu\'{e}e. Soit $m:\hs\otimes \hs\ra\hs$ une multiplication associative commutative gradu\'{e}e
(de degr\'{e} $0$) sur $\hs$. Alors sa d\'{e}cal\'{e}e $m[1]$ s'annule sur
$\hs[1]\bullet \hs[1]$, et l'application induite
$d\in \Hom(\underline{\hs[1]^{\otimes 2}},\hs[1])$
satisfait $d\circ_H d=0$. La formule
$\beta:\phi\mapsto [\phi,d]_H$ d\'{e}finit le {\em cobord de Harrison gradu\'{e}} sur l'espace
$\Hom(\underline{\hs[1]^\otimes},\hs[1])$ des cocha\^{\i}nes de Harrison. On peut regarder
le cocomplexe de Harrison comme un sous-cocomplexe du cocomplexe de Hochschild.

La structure d'une {\em alg\`{e}bre de Gerstenhaber} (voir aussi \cite[p.104]{Gin04} sur un
espace vectoriel
gradu\'{e} $\gs$ est la donn\'{e}e de la structure d'une alg\`{e}bre associative
commutative gradu\'{e}e $m_\gs:\gs\otimes \gs\ra \gs$ avec \'{e}l\'{e}ment neutre
$1$ et de la structure $[~,~]$ d'une alg\`{e}bre de Lie gradu\'{e}e sur l'espace d\'{e}cal\'{e}
$\gs[1]$ satisfaisant la {\em r\`{e}gle de Leibniz gradu\'{e}e} suivante:
soit $\lambda_\gs:\gs\otimes\gs\ra\gs$ le crochet d\'{e}cal\'{e}
$\lambda_\gs:=[~,~][-1]$ (qui est donc de degr\'{e} $-1$, commutative gradu\'{e}e
($\lambda_\gs\circ\tau=\lambda_\gs$) et satisfait l'identit\'{e} de Jacobi gradu\'{e}e), alors
\begin{equation}\label{EqRegleLeibniz}
     \lambda_\gs\circ (\id\otimes m_\gs)=m_\gs \circ (\lambda_\gs\otimes \id) +
     m_{\gs}\circ (\id\otimes \lambda_\gs)\circ \tau_{12}.
\end{equation}
Dans la suite, on note une alg\`{e}bre de Gerstenhaber soit par
$(\gs,m_{\gs},1,\lambda_\gs)$ avec le crochet d\'{e}cal\'{e} $\lambda_\gs$, soit
par $(\gs[1],m_{\gs}[1],1,[~,~])$ avec la multiplication d\'{e}cal\'{e}e
$m_\gs[1]$.\\
L'exemple principal est
construit de mani\`{e}re suivante: soit $(\ls,[~,~]_\ls)$ une alg\`{e}bre de Lie
gradu\'{e}e, soit $\hs=\ls[-1]$ et soit $\lambda_\ls$ le crochet d\'{e}cal\'{e}
$[~,~][-1]$ donn\'{e} par
$s\circ[~,~]_\ls\circ(s\otimes s)^{-1}$ qui est de degr\'{e} $-1$. On pose
$\gs:=\mathsf{S}\hs=\mathsf{S}\big(\ls[-1]\big)$. L'application
$\hs\leftarrow \mathsf{S}\hs\otimes\mathsf{S}\hs$ donn\'{e}e par $\lambda_\ls\circ
(\mathrm{pr}_\hs\otimes\mathrm{pr}_\hs)$ co-induit une cod\'{e}rivation
$\overline{\lambda_\ls}=:\lambda_\gs:\mathsf{S}\hs\otimes\mathsf{S}\hs\ra
\mathsf{S}\hs$ le long de la multiplication shuffle $\mu_{\mathrm{sh}}$, et
on voit sans peine que
$(\mathsf{S}\big(\ls[-1]\big),\mu_{\mathrm{sh}},1,\lambda_\gs)$ est une alg\`{e}bre
de Gerstenhaber gradu\'{e}e. \\
Un exemple g\'{e}om\'{e}trique est donn\'{e} par {\em l'espace
de tous les champs de multivecteurs sur une vari\'{e}t\'{e} diff\'{e}rentielle $M$},
$\gs=\Gamma^\infty(M,\Lambda TM)$, muni de la multiplication ext\'{e}rieure
point-par-point $\wedge=:m$ et du {\em crochet de Schouten} $[~,~]_S$ sur
$\gs[1]$, voir le paragraphe \ref{SubSecPropalggs}. \\
Plus g\'{e}n\'{e}ralement, soit
$V=V^0$ un espace vectoriel, $m:V\otimes V\ra V$ une multiplication
associative, $\Gs:=\Hom(V[1]^\otimes,V[1])$ muni du cobord de Hochschild
$\mathsf{b}$ et $\hs$ la cohomologie du co-complexe $(\Gs,\mathsf{b})$.
D'apr\`{e}s un r\'{e}sultat classique de Gerstenhaber \cite{Ger63}, $\hs[-1]$ est
une alg\`{e}bre de Gerstenhaber o\`{u} le crochet de Lie sur $\hs$ est induit par
le crochet de Gerstenhaber sur $\Gs$ et la multiplication sur $\hs[-1]$
provient de la multiplication `{\em cup}' $\cup$ sur
$\Hom(V^\otimes,V)$ d\'{e}finie par la convolution de deux \'{e}l\'{e}ments de
$\Hom(V^\otimes,V)$ par rapport \`{a} l'alg\`{e}bre
associative $(V,m)$ et la cog\`{e}bre $(V^\otimes,\Delta)$.

De mani\`{e}re analogue, la structure d'une {\em cog\`{e}bre de Gerstenhaber gradu\'{e}e}
sur $\gs$ est la donn\'{e}e de la structure d'une cog\`{e}bre coassociative
cocommutative gradu\'{e}e $\Delta_\gs:\gs\ra \gs\otimes \gs$ avec co-unit\'{e}
$\epsilon_\gs$
et de la structure $\delta_\gs$ d'une cog\`{e}bre de Lie gradu\'{e}e sur l'espace d\'{e}cal\'{e}
$\gs[-1]$ satisfaisant la {\em r\`{e}gle de co-Leibniz gradu\'{e}e} suivante:
soit $\kappa_\gs:\gs\ra\gs\otimes \gs$ le cocrochet d\'{e}cal\'{e}
$\kappa_\gs:=\delta_\gs[1]$ (qui est donc de degr\'{e} $-1$,
cocommutative gradu\'{e}e
($\tau\circ \kappa_\gs=\kappa_\gs$) et satisfait l'identit\'{e} de co-Jacobi gradu\'{e}e), alors
\begin{equation}\label{EqRegleCoLeibniz}
  (\id\otimes \Delta_\gs)\circ \kappa_\gs= (\kappa_\gs\otimes \id)\circ \Delta_\gs +
     \tau_{12}\circ (\id\otimes \kappa_\gs)\circ\Delta_\gs.
\end{equation}
ou $(\Delta_\gs\otimes \id)\circ \kappa_\gs
= \tau_{23}\circ(\kappa_\gs\otimes \id)\circ \Delta_\gs +
(\id\otimes \kappa_\gs)\circ\Delta_\gs$ comme condition \'{e}quivalente.
Dans la suite, on notera une cog\`{e}bre de Gerstenhaber
$(\gs,\Delta_{\gs},\epsilon_\gs,\kappa_\gs)$
avec le cocrochet d\'{e}cal\'{e} $\kappa_\gs$. On retrouve $\delta_\gs$ par
$\kappa_\gs[-1]$. Un {\em morphisme de cog\`{e}bres de Gerstenhaber}
$\phi:(\gs,\Delta_{\gs},\epsilon_\gs,\kappa_\gs)\ra
(\gs',\Delta_{\gs'},\epsilon_{\gs'},\kappa_{\gs'})$
est un morphisme de cog\`{e}bres
qui est en m\^{e}me temps un morphisme de cog\`{e}bres de Lie, i.e.
$(\phi\otimes \phi)\circ \Delta_\gs=\Delta_{\gs'}\circ\phi$ et
$(\phi\otimes \phi)\circ \kappa_\gs=\kappa_{\gs'}\circ\phi$. De plus, une {\em
cod\'{e}rivation de cog\`{e}bres de Gerstenhaber}
$d:(\gs,\Delta_{\gs},\epsilon_\gs,\kappa_\gs)\ra
(\gs',\Delta_{\gs'},\epsilon_{\gs'},\kappa_{\gs'})$
de degr\'{e} $j$ le long d'un morphisme de cog\`{e}bres de Gerstenhaber $\phi$ est une
application lin\'{e}aire de degr\'{e} $j$ telle que
$(\phi\otimes d + d\otimes \phi)\circ \Delta_\gs=\Delta_{\gs'}\circ d$ et
$(\phi\otimes d + d\otimes \phi)\circ \kappa_\gs=(-1)^j \kappa_{\gs'}\circ d$.\\
Soit $(\mathfrak{c},\delta_\mathfrak{c})$ une cog\`{e}bre de Lie gradu\'{e}e. Soit
$\hs:=\mathfrak{c}[1]$ donc $\mathfrak{c}=\hs[-1]$, et
$\delta_\mathfrak{c}[1]:\hs\ra\hs\otimes\hs$ le cocrochet d\'{e}cal\'{e}. L'application
$\hs\stackrel{\delta_\mathfrak{c}[1]}{\ra}
\hs\otimes \hs \ra \mathsf{S}\hs\otimes \mathsf{S}\hs$ se
prolonge de mani\`{e}re unique
comme une d\'{e}rivation gradu\'{e}e $\kappa_\delta$ de degr\'{e} $-1$ de l'alg\`{e}bre commutative gradu\'{e}e
libre
$(\mathsf{S}\hs,\mu_{\mathrm{sh}})$ dans l'alg\`{e}bre commutative gradu\'{e}e $\mathsf{S}\hs\otimes
\mathsf{S}\hs$ le long de l'homomorphisme $\Delta_{\mathrm{sh}}$. On calcule sans peine que
$(\mathsf{S}\big(\mathfrak{c}[1]\big),\Delta_{\mathrm{sh}},\epsilon,\kappa_\delta)$ est
toujours une cog\`{e}bre de Gerstenhaber. De plus, \'{e}tant donn\'{e} un morphisme
$\phi:(\mathfrak{c}_1,\delta_{\mathfrak{c}_1}\ra(\mathfrak{c}_2,\delta_{\mathfrak{c}_1})$
de cog\`{e}bres de Lie, il vient que l'unique morphisme de cog\`{e}bres coassociatives cocommutatives
gradu\'{e}es $\overline{\phi[1]}:\mathsf{S}(\mathfrak{c}_1[1])\ra \mathsf{S}(\mathfrak{c}_2[1])$
co-induit par $\phi[1]$ est un morphisme de cog\`{e}bres de Gerstenhaber.\\
Soit $\mathcal{CG}_{AN}$ la cat\'{e}gorie des cog\`{e}bres de Gerstenhaber telles
que
la cog\`{e}bre associative est dans $\mathcal{C}_{AN}$, le cocrochet d\'{e}cal\'{e} s'annule sur $1$ et
la cog\`{e}bre de Gerstenhaber quotient par $\korps 1$ est nilpotente dans le sens que pour
tout \'{e}l\'{e}ment il existe un entier positif $N$ tel que toutes les $N$\`{e}mes
it\'{e}rations de la comultiplication et du cocrochet d\'{e}cal\'{e} s'annulent. Soit $\hs$ un espace
vectoriel gradu\'{e}. En
prenant la cog\`{e}bre de Lie colibre $\mathfrak{c}=\underline{\hs[-1]^\otimes}$
on voit sans peine que la cog\`{e}bre de Gerstenhaber
$\mathsf{S}\big(\mathfrak{c}[1]\big)$, qui est \'{e}gale \`{a}
\begin{equation}
       \mathsf{G}\hs:=\mathsf{S}\Big(\big(\underline{\hs[-1]^\otimes}\big)[1]\Big)
\end{equation}
est une {\em cog\`{e}bre de Gerstenhaber colibre} dans la cat\'{e}gorie
$\mathcal{CG}_{AN}$: en fait, \'{e}tant donn\'{e}e une application lin\'{e}aire $\phi:C\ra\hs$ de degr\'{e} $0$
dans le diagramme (\ref{EqDiagrammeColibre}) on
utilise d'abord la structure de cog\`{e}bre de Lie de $C[-1]$ pour construire de
$\phi[-1]:=s\circ\phi\circ s^{-1}$ un
morphisme de cog\`{e}bres de Lie gradu\'{e}es $\psi:=\overline{\phi[-1]}:C[-1]\ra
\underline{\hs[-1]^\otimes}$. Ensuite, le morphisme de cog\`{e}bres $\overline{\psi[1]}:C\ra
\mathsf{S}\Big(\big(\underline{\hs[-1]^\otimes}\big)[1]\Big)$ co-induit par $\psi[1]$
pr\'{e}serve aussi les structures de cog\`{e}bres de Lie gradu\'{e}es sur $C$ et $\mathsf{G}\hs$.\\
De la m\^{e}me fa\c{c}on, on peut co-induire des applications lin\'{e}aires $d:C^+\ra \hs$
de degr\'{e} $k$ comme cod\'{e}rivations gradu\'{e}es $\overline{d}$ de degr\'{e} $k$ de
$C\ra\mathsf{G}\hs$ le long
d'un morphisme $\overline{\phi}:C\ra\mathsf{G}\hs$ de cog\`{e}bres de Gerstenhaber
gradu\'{e}es co-induit par $\phi:C^+\ra\hs$. Le cas particulier $C=\mathsf{G}\hs$ et
$\overline{\phi}=\id_{\mathsf{G}\hs}$ est int\'{e}ressant: puisque $\overline{d}$ est
uniquement d\'{e}termin\'{e}e par sa composante $\mathrm{pr}_\hs(\overline{d})=d$,
on introduit une multiplication $d_1\circ_T d_2:=d_1\circ \overline{d_2}$
sur l'espace $\Hom(\mathsf{G}^+\hs,\hs)$ d'apr\`{e}s le moule
(\ref{EqMultGerstenhaber}). Cette multiplication satisfait (\ref{EqIdGerstenhaber})
(avec $\circ_G$ remplac\'{e} par $\circ_T$), et le commutateur gradu\'{e} $[~,~]_T$
(l'antisym\'{e}tris\'{e} gradu\'{e} de $\circ_T$) est la structure d'une alg\`{e}bre de
Lie gradu\'{e}e sur $\Hom(\mathsf{G}^+\hs,\hs)$.

\subsection{Un aper\c{c}u des structures \`{a} homotopie pr\`{e}s}
 \label{AppSubSecFInfini}

Pour d\'{e}finir des {\em structures $\mathbf{F_\infty}$ (pour $\mathbf{F}=\mathbf{A,L,C,G}$)
sur un espace vectoriel
gradu\'{e} $\hs$}, on consid\`{e}re d'abord l'espace d\'{e}cal\'{e} $\hs[1]$:

$\mathbf{A_\infty}$: on regarde la cog\`{e}bre colibre gradu\'{e}e
$\mathbf{A_\infty}\hs:=\hs[1]^\otimes$.
 Une {\em structure $\mathbf{A_\infty}$ sur
$\hs$} est un \'{e}l\'{e}ment $d\in \Hom(\hs[1]^{\otimes +},\hs[1])$ de degr\'{e} $1$ avec
\begin{equation}\label{EqDefAinfini}
        d\circ_G d=0.
\end{equation}
De mani\`{e}re \'{e}quivalente, on parle d'une cod\'{e}rivation gradu\'{e}e $\overline{d}$
(co-induite par $d$) de degr\'{e}
$1$ et de carr\'{e} nul, i.e. $\overline{d}\circ\overline{d}=0$. Alors le
couple $(\hs[1]^\otimes,\overline{d})$ est appel\'{e} une {\em cog\`{e}bre codiff\'{e}rentielle
gradu\'{e}e}. Puisque l'espace vectoriel $\hs[1]^\otimes$ a une deuxi\`{e}me graduation
$\hs[1]^\otimes=\oplus_{r\in\mathbb N}
\hs[1]^{\otimes r}$ (on \'{e}crira parfois $\hs[1]^{\otimes r}=:\hs[1]^{[r]}$),
l'application $d$ se d\'{e}compose comme
$d=\sum_{r=1}^\infty d^r$ o\`{u} les $d^r$ se trouvent dans $\Hom(\hs[1]^{\otimes r},\hs[1])^1$
quel que soit $r\in\mathbb N$. Les d\'{e}cal\'{e}s
\begin{equation}\label{EqAinfiniDecales}
m^r:= d^r[-1]=s\circ d^r \circ (s^{\otimes r})^{-1} \in \Hom(\hs^{\otimes r},\hs)^{2-r}
\end{equation}
sont donc des morphismes de degr\'{e} $2-r$. En particulier, $m^1$ est de degr\'{e} $1$
et $m^2$ est de degr\'{e} $0$. L'identit\'{e} (\ref{EqDefAinfini}) pour $r=1$ est \'{e}quivalente \`{a} dire
que $m^1$ est une {\em codiff\'{e}rentielle} sur $\hs$ ($m^1m^1=0$). Pour $\alpha\in
\Hom(\hs^{\otimes k},\hs)$ soit $(m^1.\alpha)$ l'action usuelle d'une codiff\'{e}rentielle, i.e.
$(m^1.\alpha)=m^1\circ\alpha-(-1)^{|\alpha|}
\sum_{i=0}^{k-1}\alpha\circ\big(\id^{\otimes i}\otimes m^1 \otimes \id^{\otimes k-1-i}\big)$.
Au rang $r=2$ (\ref{EqDefAinfini}) exprime la compatibilit\'{e} de $m^1$ avec $m^2$, i.e.
$m^1.m^2=0$. Au rang $3$ de (\ref{EqDefAinfini}) on voit que
l'associateur de $m^2$,
$m^2\circ (m^2\otimes \id_\hs)-m^2\circ (\id_\hs\otimes m^2)$ est
un cobord, i.e. proportionnel \`{a} $m^1.m^3$ ce qui explique le nom {\em structure associative
\`{a} homotopie pr\`{e}s}. La structure d'une {\em alg\`{e}bre associative codiff\'{e}rentielle gradu\'{e}e} sur
$\hs$
est un exemple d'une structure $A_\infty$ avec $m^r=0$ pour $r\geq 3$.
Les signes pr\'{e}cis pour les identit\'{e}s des quantit\'{e}s d\'{e}cal\'{e}es
$m^r$ se calculent \`{a} l'aide de la r\`{e}gle de signe (\ref{EqRegleKoszul}),
(\ref{EqApplMultDecalees}) et co\"{\i}ncident
avec
ceux de Stasheff, 1963, voir \cite[p.294, Def.~2.1]{Sta63}. Soient $(\hs,d)$
et $(\hat{\hs},\hat{d})$ des structures $A_\infty$ sur $\hs$ et sur $\hat{\hs}$.
Un {\em morphisme d'alg\`{e}bres $A_\infty$} est un morphisme de cog\`{e}bres codiff\'{e}rentielles
gradu\'{e}es $\Phi:\hs[1]^\otimes \ra \hat{\hs}[1]^\otimes$, i.e. $\Phi$
pr\'{e}serve les comultiplications et les codiff\'{e}rentielles dans le sens
\begin{equation}\label{EqMorphAinfini}
 (\Phi\otimes \Phi)\circ \Delta = \hat{\Delta}\circ \Phi~~\mbox{  et  }~~
    \Phi\circ \overline{d}=\overline{\hat{d}}\circ\Phi.
\end{equation}

\noindent Toutes les autres structures $\mathbf{F_\infty}$ sont con\c{c}ues exactement
d'apr\`{e}s le m\^{e}me moule:

$\mathbf{L_\infty}$: on regarde la cog\`{e}bre colibre cocommutative gradu\'{e}e
$\mathbf{L_\infty}\hs:=\mathsf{S}\big(\hs[1]\big)$ qui est isomorphe
en tant qu'espace vectoriel \`{a} $\ve \hs$.
Une {\em structure $\mathbf{L_\infty}$ sur
$\hs$} est un \'{e}l\'{e}ment $d_L\in \Hom(\mathsf{S}^+\hs[1],\hs[1])$ de degr\'{e} $1$ avec
$d_L\circ_{NR} d_L=0$. A l'aide de la cod\'{e}rivation $\overline{d_L}$ de
$\mathsf{S}\big(\hs[1]\big)$
co-induite par $d$, le couple $\big(\mathsf{S}\big(\hs[1]\big),\overline{d_L}\big)$
est une cog\`{e}bre cocommutative codiff\'{e}rentielle gradu\'{e}e. En utilisant la deuxi\`{e}me graduation
de $\mathsf{S}(\hs[1])=\oplus_{r\in\mathbb N}\mathsf{S}^r(\hs[1])=:
\oplus_{r\in\mathbb N}\hs[1]^{[r]}$ on d\'{e}finit comme avant les applications $d_L^r$, et
les structures
d\'{e}cal\'{e}es $m_L^r:=d_L^r[-1]$ comme dans (\ref{EqAinfiniDecales}) sont des applications
de degr\'{e} $2-r$ dans l'espace $\Hom\big(\hs^{\otimes r},\hs\big)$ qui sont
automatiquement antisym\'{e}triques gradu\'{e}es. L'application $m_L^1$ est une codiff\'{e}rentielle,
on a
$m_L^1.m_L^2=0$, et le {\em jacobiateur} de $m_L^2$,
$m_L^2\circ (m_L^2\otimes \id)\circ
\big(\id^{\otimes 3}+\tau_{12}\tau_{23}+\tau_{23}\tau_{12}\big)$, est proportionnel \`{a}
$m_L^1.m_L^3$, d'o\`{u} le nom {\em structure d'alg\`{e}bre de Lie \`{a} homotopie pr\`{e}s}.
Les signes dans les identit\'{e}s des $m_L^r$ sont induits par ceux qu'on obtient \`{a}
l'aide de la r\`{e}gle de signe, voir \cite{LS93}.
Les morphismes $L_\infty$ sont des morphismes de cog\`{e}bres
codiff\'{e}rentielles gradu\'{e}es comme dans (\ref{EqMorphAinfini}). Un cas
particulier est la structure d'une {\em alg\'{e}bre de Lie codiff\'{e}rentielle gradu\'{e}e}
sur $\hs$, i.e. $m_L=m_L^1+m_L^2$.

$\mathbf{C_\infty}$: on regarde la cog\`{e}bre de Lie colibre gradu\'{e}e
$\mathbf{C_\infty}\hs:=\underline{\hs[1]^\otimes}$. Une {\em structure $\mathbf{C_\infty}$ sur
$\hs$} est la donn\'{e}e d'un \'{e}l\'{e}ment $d\in \Hom(\underline{\hs[1]^\otimes},\hs[1])$ de degr\'{e}
$1$ avec $d\circ_{H} d=0$. Tout ce qui a \'{e}t\'{e} dit pour le cas
$\mathbf{A_\infty}$ se traduit de mani\`{e}re analogue, par exemple la composante d\'{e}cal\'{e}e
$m^2:=d^2[-1]$ est une {\em multiplication associative commutative gradu\'{e}e
\`{a} homotopie pr\'{e}s}. Pour $r\in\mathbb N$, on \'{e}crira parfois la notation g\'{e}n\'{e}rale
$\hs[1]^{[r]}$ pour $\underline{\hs[1]^{\otimes r}}$.
Un cas
particulier est la structure d'une {\em alg\'{e}bre associative commutative codiff\'{e}rentielle
gradu\'{e}e} sur $\hs$, i.e. $m=m^1+m^2$.

$\mathbf{G_\infty}$: on regarde la cog\`{e}bre de Gerstenhaber colibre gradu\'{e}e
$\mathbf{G_\infty}\hs:=\mathsf{G}\big(\hs[1]\big)
=\mathsf{S}\big((\underline{\hs^\otimes})[1]\big)$
qui est isomorphe en tant qu'espace vectoriel \`{a} $\ve \underline{\hs^\otimes}$, notation
plus simple qu'on utilisera souvent.
Une {\em structure $\mathbf{G_\infty}$ sur
$\hs$} est la donn\'{e}e d'un \'{e}l\'{e}ment $d\in \Hom(\mathsf{G}^+(\hs[1]),\hs[1])$ de degr\'{e} $1$ avec
$d\circ_{T} d=0$. La famille des sous-espaces $\hs[1]^{[p_1,\ldots,p_r]}$ de
$\mathsf{G}\big(\hs[1]\big)$ pour $r,p_1,\ldots,p_r\in\mathbb N$ d\'{e}finis par
$(\underline{\hs^{\otimes p_1}})[1]
\bullet\cdots\bullet (\underline{\hs^{\otimes p_r}})[1]$
d\'{e}termine une d\'{e}composition de
$\mathsf{G}\big(\hs[1]\big)$, et les sommes
$\hs[1]^{[n]}\!\!\!:=\oplus_{r,p_1+\cdots +p_r=n}\hs[1]^{[p_1,\ldots,p_r]}$
donnent une deuxi\`{e}me graduation sur $\mathsf{G}\big(\hs[1]\big)$ dont la
filtration correspondante, $\hs^{[\leq n]}:=\oplus_{k=1}^n \hs[1]^{[k]}$,
est pr\'{e}serv\'{e}e par les structures alg\'{e}briques $\Delta$, $l$ et
$\overline{d}$. On note $d^{p_1,\ldots ,p_r}$ la restriction de $d$ \`{a}
$\hs[1]^{[p_1,\ldots,p_r]}$. Les composantes $d^{1}$, $d^{1,1}$ et $d^2$ sont particuli\`{e}rement
importantes: la d\'{e}cal\'{e}e $m^1:=d^1[-1]$ est une codiff\'{e}rentielle
sur $\hs$. De plus, la d\'{e}cal\'{e}e $m^{1,1}:=d^{1,1}[-1]$ se restreint \`{a} la structure
d'une alg\`{e}bre de Lie gradu\'{e}e \`{a} homotopie pr\`{e}s sur $\hs$. Finalement, la restriction
de la d\'{e}cal\'{e}e $m^2:=d^2[-1]$ \`{a} $\hs$ est \'{e}gale \`{a} la d\'{e}cal\'{e}e $m_{\hs[-1]}[1]$
d'une structure d'alg\`{e}bre commutative
gradu\'{e}e associative \`{a} homotopie pr\`{e}s sur $\hs[-1]$.\\
Si $(\hs,m_{\hs[-1]}[1],1,[~,~],\mathsf{b})$ est une alg\`{e}bre de Gerstenhaber codiff\'{e}rentielle
(i.e. la codiff\'{e}rentielle $\mathsf{b}$
pr\'{e}serve la multiplication $m_{\hs[-1]}$ et le crochet de Lie
gradu\'{e} $[~,~]$, il
s'ensuit que ceci induit une structure $G_\infty$ sur $\mathsf{G}(\hs[1])$
avec $d=d^1+d^{1,1}+d^2$.\\
Une g\'{e}n\'{e}ralisation importante de ce cas d'une structure $G_\infty$ s'obtient de la
fa\c{c}on suivante: soit
$(\mathfrak{b},[~,~]_\mathfrak{b},\delta_\mathfrak{b},\mathsf{b})$ une
{\em big\`{e}bre de Lie codiff\'{e}rentielle gradu\'{e}e}, i.e. $[~,~]_\mathfrak{b}$
est une structure d'alg\`{e}bre de Lie gradu\'{e}e, $\delta_\mathfrak{b}$ est une
structure de cog\`{e}bre de Lie gradu\'{e}e,
$\mathsf{b}:\mathfrak{b}\ra\mathfrak{b}$ est une application lin\'{e}aire de
degr\'{e} $1$ et de carr\'{e} nul telles que $\mathsf{b}$ pr\'{e}serve
$[~,~]_\mathfrak{b}$
et $\delta_\mathfrak{b}$, i.e.
\[
  \mathsf{b}\circ [~,~]_\mathfrak{b}=
    [~,~]_\mathfrak{b}\circ \big(\mathsf{b}\otimes \id +\id\otimes \mathsf{b}\big)
       \mathrm{~~,~~}
    \delta_\mathfrak{b}\circ\mathsf{b}
    =\big(\mathsf{b}\otimes \id +\id\otimes
    \mathsf{b}\big)\circ\delta_\mathfrak{b},
\]
et $[~,~]_\mathfrak{b}$ et $\delta_\mathfrak{b}$ sont compatibles, i.e.
\beas
    \delta_\mathfrak{b}\circ [~,~]_\mathfrak{b}
   & = & ([~,~]_\mathfrak{b}\otimes \id)
       \circ\big(\tau_{23}\circ(\delta_\mathfrak{b}\otimes \id)
                        +\id\otimes \delta_\mathfrak{b}\big)
       \\
   &   & + (\id\otimes [~,~]_\mathfrak{b})
       \circ \big(\delta_\mathfrak{b}\otimes \id
                   +\tau_{12}\circ(\id\otimes\delta_\mathfrak{b})\big).
\eeas
Dans le cas particulier o\`{u} $\mathfrak{b}$ est la cog\`{e}bre de Lie gradu\'{e}e
colibre $\underline{\hs^\otimes}$ munie du cocrochet $\delta$ canonique, il est
clair que la codiff\'{e}rentielle $\mathsf{b}$ est co-induite par sa
composante $\mathrm{pr}_\hs\circ\mathsf{b}=\sum_{p=1}^\infty \mathsf{b}^{p}$.
De plus, on peut montrer que tout crochet antisym\'{e}trique gradu\'{e} $[~,~]$ qui est compatible
avec $\delta$ est \'{e}galement d\'{e}termin\'{e} par sa composante
$\mathrm{pr}_\hs\circ [~,~]=\sum_{p_1,p_2=1}^\infty l^{p_1,p_2}$ o\`{u}
$l^{p_1,p_2}:\underline{\hs^{\otimes p_1}}\otimes\underline{\hs^{\otimes
p_2}}\ra\hs$.
Pour revenir au cas d'une big\`{e}bre de Lie gradu\'{e}e codiff\'{e}rentielle
g\'{e}n\'{e}rale,
on construit d'abord la cog\`{e}bre de Gerstenhaber
$\big(\mathsf{S}\big(\mathfrak{b}[1]\big),\Delta_{\mathrm{sh}},
\epsilon,l_{\delta_\mathfrak{b}}\big)$
\`{a} l'aide de la cog\`{e}bre de Lie gradu\'{e}e $(\mathfrak{b},\delta_\mathfrak{b})$,
comme d\'{e}j\`{a} mentionn\'{e} ci-dessus. La d\'{e}cal\'{e}e $d^1$ de $\mathsf{b}$ est une
application lin\'{e}aire de degr\'{e} $1$ de $\mathfrak{b}[1]$ dans
$\mathfrak{b}[1]$, et la d\'{e}cal\'{e}e $d^2$ de $[~,~]_\mathfrak{b}$ est une
application lin\'{e}aire de degr\'{e} $1$ de $S^2\big(\mathfrak{b}[1]\big)$ dans
$\mathfrak{b}[1]$. La cod\'{e}rivation
$\overline{d}$ de la cog\`{e}bre $\big(\mathsf{S}\big(\mathfrak{b}[1]\big),\Delta_{\mathrm{sh}},
\epsilon\big)$ co-induite par la somme $d:=d^1+d^2$ pr\'{e}serve
$l_{\delta_\mathfrak{b}}$ (gr\^{a}ce \`{a} la deuxi\`{e}me et troisi\`{e}me \'{e}quation de
compatibilit\'{e} ci-dessus) et est de carr\'{e} $0$ car $\mathsf{b}$ l'est,
$[~,~]_\mathfrak{b}$ satisfait l'identit\'{e} de Jacobi gradu\'{e}e et $\mathsf{b}$ pr\'{e}serve
$[~,~]_\mathfrak{b}$. Alors $\big(\mathsf{S}\big(\mathfrak{b}[1]\big),\Delta_{\mathrm{sh}},
\epsilon,l_{\delta_\mathfrak{b}}, \overline{d}\big)$ est une cog\`{e}bre de
Gerstenhaber codiff\'{e}rentielle gradu\'{e}e. Ceci implique la proposition
suivante montr\'{e}e dans \cite[p.43, Lemma 2.1]{GH03}:
\begin{prop}\label{Lemma 1.1}
Supposons donn\'{e}e une structure de big\`{e}bre de Lie
codiff\'{e}rentielle sur la cog\`{e}bre de
Lie gradu\'{e}e colibre $\underline{\hs^\otimes}$
dont la diff\'{e}rentielle et le crochet de Lie sont co-induits
respectivement par les applications $\mathsf{b}^p:\underline{\hs^{\otimes p}}\ra\hs$
et $l^{p_1,p_2}:\underline{\hs^{\otimes p_1}}\otimes \underline{\hs^{\otimes p_2}}\ra
\hs$.
Alors $\hs$ a une structure $G_\infty$ $d$ donn\'{e}e sur $\mathsf{G}\big(\hs[1]\big)$,
pour tous entiers $p,p_1,p_2,\ldots,p_r,\ldots\geq 1$, par les d\'{e}cal\'{e}es
\[
  d^p:=s^{-1}\circ\mathsf{b}^p\circ s_{\underline{\hs^{\otimes}}}, \qquad
  d^{p_1,p_2}:=s^{-1}\circ l^{p_1,p_2}\circ\big(s_{\underline{\hs^{\otimes}}}
                                                \otimes s_{\underline{\hs^{\otimes}}}\big)
\]
o\`{u} $d^{p_1,\dots, p_r}=0$ pour $r\geq 3$.
\end{prop}
Pour \'{e}noncer une {\em liaison entre structures $G_\infty$ et $L_\infty$},
on remarque que l'espace vectoriel gradu\'{e} $\hs$ munie du cocrochet trivial est une sous-cog\`{e}bre
de Lie de la cog\`{e}bre de Lie colibre $\underline{\hs^\otimes}$ via
l'injection canonique de $\hs$ sur $\underline{\hs^{\otimes 1}}$: en
effet, $\delta(x)=\Delta(x)-\tau\Delta(x)=0$ pour tout $x\in \hs$.
L'application co-induite $\overline{i}:\mathsf{S}(\hs[1])\ra
\mathsf{S}\big((\underline{\hs^\otimes})[1]\big)$ est donc un morphisme de
cog\`{e}bres de Gerstenhaber injectif.
\begin{prop} \label{PropGinfiniLinfini}
 Soient $\hs$ et $\Hs$ deux espaces vectoriels gradu\'{e}s.
 \begin{enumerate}
  \item Soit $d=\sum_{r=1}^\infty\sum_{p_1,\ldots,p_r=1}^\infty d^{p_1,\ldots,p_r}$
  une structure $G_\infty$ sur $\hs$. Alors la restriction $d_L$ de $d$ \`{a} la
  sous-cog\`{e}bre de Gerstenhaber $\mathsf{S}(\hs[1])$ de $\mathsf{G}(\hs[1])$
  d\'{e}finit une unique
  structure $L_\infty$ sur $\hs$ dont les composantes sont don\-n\'{e}es par
  $d_L^r:=d^{1,\ldots,1}$.
  \item Soient $d$ et $D$ des structures $G_\infty$ sur $\hs$ et $\Hs$, respectivement, et
  soit $\overline{\Phi}:(\mathsf{G}(\hs[1]),d)\ra(\mathsf{G}(\Hs[1]),D)$ un
  morphisme $G_\infty$ co-induit par l'application
  $\Phi:\mathsf{G}(\hs[1])\ra \Hs[1]$. Alors la restriction $\overline{\phi}$
  de $\overline{\Phi}$ \`{a} la sous-cog\`{e}bre $\mathsf{S}(\hs[1])$ prend ses
  valeurs dans $\mathsf{S}(\Hs[1])$ et d\'{e}finit un morphisme $L_\infty$
  entre $(\mathsf{S}(\hs[1]),d_L)$ et $(\mathsf{S}(\Hs[1]),D_L)$
 \end{enumerate}
\end{prop}
\begin{prooof}
  On \'{e}crit $\mathrm{pr}^G_\hs:\mathsf{G}(\hs[1])\ra \hs[1]$ pour la
  projection canonique.\\
 1. L'application $\overline{d}\circ \overline{i}$ est une cod\'{e}rivation
 de $\mathsf{S}(\hs[1])$ dans $\mathsf{G}(\hs[1])$ le long du morphisme
 $\overline{i}$. Elle est donc co-induite par sa projection
 $d_L:=\mathrm{pr}^G_{\hs[1]}\circ \overline{d}\circ \overline{i}$. Soit
 $\overline{d_L}$ la cod\'{e}rivation de la cog\`{e}bre $\mathsf{S}(\hs[1])$
 co-induite par $d_L$. Alors $\overline{i}\circ\overline{d_L}$ est \'{e}galement
 une cod\'{e}rivation de $\mathsf{S}(\hs[1])$ dans $\mathsf{G}(\hs[1])$ le long du morphisme
 $\overline{i}$. Puisque $\mathrm{pr}^G_{\hs[1]}\circ \overline{d}\circ \overline{i}$
 est \'{e}gal \`{a} $d_L=\mathrm{pr}_{\hs[1]}\circ \overline{d_L}
 =\mathrm{pr}^G_{\hs[1]}\circ \overline{i}\circ\overline{d_L}$ on peut
 conclure que $\overline{d}\circ
 \overline{i}=\overline{i}\circ\overline{d_L}$, et $\overline{i}$ est
 donc un morphisme de cog\`{e}bres de Gerstenhaber gradu\'{e}es codiff\'{e}rentielles.
 Gr\^{a}ce \`{a} cette \'{e}quation et l'injectivit\'{e} de $\overline{i}$, il vient
 que $\overline{d_L}\circ\overline{d_L}=0$, et $d_L$ est donc une
 structure $L_\infty$ sur $\hs$ dont les composantes sont des restrictions
 de $d$ \`{a} $\mathsf{S}(\hs[1])$, \`{a} savoir $d^{1,\ldots,1}$.\\
 2. D'apr\`{e}s ce qui pr\'{e}c\`{e}de, l'application $\overline{\Phi}\circ
 \overline{i}_{\mathsf{S}(\hs[1])}:\mathsf{S}(\hs[1])\ra\mathsf{G}(\Hs[1])$
 est un morphisme de cog\`{e}bres de Gerstenhaber codiff\'{e}rentielles, et donc
 co-induit par sa projection $\phi:=\mathrm{pr}^G_{\Hs[1]}\circ
 \overline{\Phi}\circ\overline{i}_{\mathsf{S}(\hs[1])}$. Soit
 $\overline{\phi}:\mathsf{S}(\hs[1])\ra \mathsf{S}(\Hs[1])$ le morphisme de cog\`{e}bres
 cocommutatives gradu\'{e}es co-induit par $\phi$. Il vient que les deux
 morphismes de cog\`{e}bres de Gerstenhaber $\overline{\Phi}\circ
 \overline{i}_{\mathsf{S}(\hs[1])}$ et
 $\overline{i}_{\mathsf{S}(\Hs[1])}\circ \overline{\phi}$ ont la m\^{e}me
 composante dans $\Hs[1]$, et sont donc \'{e}gaux. Ensuite,
 \beas
   \overline{i}_{\mathsf{S}(\Hs[1])}\circ \overline{\phi} \circ
   \overline{d_L}
   & = &\overline{\Phi}\circ\overline{i}_{\mathsf{S}(\hs[1])}\circ
    \overline{d_L}
   = \overline{\Phi}\circ\overline{d}\circ \overline{i}_{\mathsf{S}(\hs[1])}
   =\overline{D}\circ\overline{\Phi}\circ
   \overline{i}_{\mathsf{S}(\hs[1])}\\
   & = &\overline{D}\circ\overline{i}_{\mathsf{S}(\Hs[1])}\circ\overline{\phi}
   =\overline{i}_{\mathsf{S}(\Hs[1])}\circ \overline{D_L}\circ
   \overline{\phi},
 \eeas
 donc gr\^{a}ce \`{a} l'injectivit\'{e} de $\overline{i}_{\mathsf{S}(\Hs[1])}$,
 l'application $\overline{\phi}$ est un morphisme $L_\infty$.
\end{prooof}

\subsection{Accolades}
 \label{AppSubSecAccolades}

Soit
$\hs$ un espace vectoriel gradu\'{e}. Pour chaque entier $j$ on choisit un sous-espace $\Hs^j$
de l'espace $\Hom_\korps(\hs^{\otimes +},\hs)^j$ satisfaisant la condition que
pour tous $\phi\in\Hs^j,\psi\in\Hs^k$ et pour tout $i\in\mathbb N$ la
`$i$-\`{e}me composition'
\begin{equation}\label{EqImeComposition}
  \phi\circ_i\psi:=
        \sum_{j\in\mathbb N}\phi\circ\big(\mathrm{id}^{\otimes (i-1)}\otimes\psi\otimes
                        \mathrm{id}^{\otimes j}\big)
\end{equation}
soit un \'{e}l\'{e}ment de $\Hs^{j+k}$. Soit $\Hs$ l'espace gradu\'{e}
$\oplus_{j\in\mathbb Z}\Hs^j$. On consid\`{e}re l'application
\[
   r:\Hs^\otimes \otimes \hs^\otimes \ra \hs: \xi\otimes x\mapsto
   \mathrm{pr}_\Hs(\xi)\big(x\big)+\epsilon_\Hs (\xi)\mathrm{pr}_\hs(x).
\]
Puisque la cog\'{e}bre gradu\'{e}e produit tensoriel
$\Hs^\otimes\otimes\hs^\otimes$ est de classe $\mathcal{C}_{AN}$ et la
cog\`{e}bre $(\hs^\otimes,\Delta_{\hs^\otimes},\epsilon_{\hs^\otimes})$ est colibre
dans cette cat\'{e}gorie, il existe un
unique morphisme de cog\`{e}bres gradu\'{e}es $\rho':\hs^\otimes \leftarrow
\Hs^\otimes\otimes \hs^\otimes$ avec $\mathrm{pr}_\hs\circ \rho'=r$.
On \'{e}crira toujours $\rho'(\xi\otimes x)=\rho(\xi)(x)$ et l'on consid\`{e}re
$\rho$ comme application lin\'{e}aire
$\Hs^\otimes \ra \Hom(\hs^\otimes,\hs^\otimes)$. Comme $\rho'$ est un
morphisme de cog\`{e}bres, il vient que pour tous $\phi_1,\cdots,\phi_k\in\Hs$
\[
   \Delta_{\hs^\otimes}\circ\rho(\phi_1\cdots \phi_k)=
        \sum_{i=1}^{k+1}
        \big(\rho(\phi_1\cdots\phi_{i-1})\otimes\rho(\phi_i\cdots\phi_k)\big)\circ
        \Delta_{\hs^\otimes}.
\]
En particulier, $\rho(\phi_1)$ est une cod\'{e}rivation de
$(\hs^\otimes,\Delta_{\hs^\otimes})$, alors
$\rho(\phi_1)=\id_{\hs^\otimes}\star\phi_1\star\id_{\hs^\otimes}$ (o\`{u} l'on
utilise la convolution $\star$ par rapport \`{a} $\mu_{\hs^\otimes}$ et $\Delta_{\hs^\otimes}$
mentionn\'{e}e ci-dessus).
On montre par r\'{e}currence que
\beas
   \rho(\phi_1\cdots\phi_k)\!\!& = &
       \id_{\hs^\otimes} \star \phi_1 \star \id_{\hs^\otimes} \star\phi_2 \star
       \id_{\hs^\otimes} \star\cdots \star\id_{\hs^\otimes}
       \star\phi_k \star \id_{\hs^\otimes} \\
       &  = &\!\!\!\!
       \sum_{i_1\cdots i_{k+1}=0}^\infty
         \id_\hs^{\otimes i_1}\otimes \phi_1 \otimes
         \id_\hs^{\otimes i_2}\otimes \phi_2 \otimes \cdots \otimes
         \id_\hs^{\otimes i_k}\otimes \phi_k \otimes \id_\hs^{\otimes i_{k+1}}
\eeas
quels que soient $\phi_1,\cdots,\phi_k\in\Hs$,
et l'on d\'{e}finit les {\em accolades (braces)} (voir \cite{Kad02},
\cite{GV95a} \cite{GJ94}) par
\[
  \phi\{\phi_1,\cdots,\phi_k\}:=\phi\circ \rho(\phi_1\cdots\phi_k).
\]
Soit $m_K$ l'application
$\Hs^\otimes\otimes \Hs^\otimes \ra\Hs$ d\'{e}finie par $m_K(\xi\otimes \eta):=
\mathrm{pr}_\Hs(\xi)\circ
\rho(\eta)+\epsilon_{\Hs}(\xi)\mathrm{pr}_\Hs(\eta)$. L'unique morphisme
$\mu_K:\Hs^\otimes\otimes \Hs^\otimes\ra \Hs^\otimes$
de cog\`{e}bres gradu\'{e}es de classe $\mathcal{C}_{AN}$  co-induit par $m_K$ (i.e. tel que
$\mathrm{pr}_\Hs\circ \mu_K=m_K)$ satisfait
\beas
  \rho'\circ\big(\id_{\Hs^\otimes}\otimes\rho') & = &
                    \rho'\circ \big(\mu_K\otimes \id_{\hs^\otimes})\\
     \mu_K\circ\big(\mu_K\otimes \id_{\Hs^\otimes}\big)
      & = &\mu_K\circ\big(\id_{\Hs^\otimes}\otimes \mu_K\big)\\
      \mu_K(e\otimes \id_{\Hs^\otimes}) & = & \id_{\Hs^\otimes} =
      \mu_K(\id_{\Hs^\otimes}\otimes e).
\eeas
o\`{u} $e$ est l'\'{e}l\'{e}ment genre groupe dans $\Hs^\otimes$: puisque les membres
de droite et de gauches des \'{e}quations pr\'{e}c\'{e}dentes sont des morphismes de
cog\`{e}bres il suffit de v\'{e}rifier leurs projections sur $\hs$ et sur
$\Hs$.
Alors $(\Hs^\otimes,\mu_K,1,\Delta_{\Hs^\otimes},\epsilon_{\Hs^\otimes})$ est une {\em big\`{e}bre
gradu\'{e}e}, et $\rho'$ d\'{e}finit sur $(\hs^\otimes,\Delta_{\hs^\otimes})$ la structure
d'une {\em module-cog\`{e}bre gradu\'{e}e \`{a} gauche} par rapport \`{a} $\Hs^\otimes$.
On a la formule suivante pour $\mu_K$, par fois \'{e}crit $\bullet_K$:
\bea
  \lefteqn{\phi_1\cdots\phi_k\bullet_K\psi_1\cdots\psi_l =} \nonumber \\
   & & \sum_{0\leq s_1\leq\cdots \leq s_{2k}\leq l}
        \prod_{r=1}^k
        (-1)^{(|\phi_r|+\cdots+|\phi_k|)(|\psi_{s_{2r-3}+1}|+\cdots+|\psi_{s_{2r-1}}|)}
        \nonumber \\
   & & ~~~~~~~~~~~~~~~\psi_1\cdots\psi_{s_1}(\phi_1\{\psi_{s_1+1},\ldots,\psi_{s_2}\})
       \psi_{s_2+1}\cdots\psi_{s_3}\nonumber \\
   & & ~~~~~~~~~~~~~~~\cdots(\phi_k\{\psi_{s_{2k-1}+1},\ldots,\psi_{s_{2k}}\})
       \psi_{s_{2k}+1}\cdots\psi_{l}\label{EqDefMultKadeiAccol}
\eea
o\`{u} $s_{-1}:=0$. La multiplication de Gerstenhaber (gradu\'{e}e) (\ref{EqMultGerstenhaber})
$\phi\circ_G \psi$ est donc donn\'{e}e par $\phi\{\psi\}=\phi\circ\rho(\psi)$ pour
$\phi,\psi\in\Hs$.
Soit $m[1]\in\Hs^1$ la d\'{e}cal\'{e}e d'une multiplication associative gradu\'{e}e
$m:\hs[-1]\otimes \hs[-1]\ra\hs[-1]$. Il vient que le commutateur gradu\'{e}
$\mathsf{b}_K:=[m[1],~]_K$
par rapport \`{a} $\mu_K$ d\'{e}finit une d\'{e}rivation gradu\'{e}e de degr\'{e} $1$ de $\mu_K$
qui est de carr\'{e} $0$ gr\^{a}ce \`{a} l'associativit\'{e} de $m$. Puisque tous les \'{e}l\'{e}ments
de $\Hs$ sont des \'{e}l\'{e}ments primitifs par rapport \`{a} $\Delta_{\Hs^\otimes}$,
i.e. $\Delta_{\Hs^\otimes}(\phi)=\phi\otimes 1+ 1\otimes \phi$ il vient
que $\mathsf{b}_K$ est \'{e}galement une cod\'{e}rivation de la cog\`{e}bre
$\Hs^{\otimes}$. Pour $\phi_1,\ldots,\phi_k\in\Hs$ on calcule
\begin{eqnarray}
 \lefteqn{\mathsf{b}_K(\phi_1\cdots\phi_k)  = } \nonumber \\
        & &\sum_{0\leq s\leq k-1}(-1)^{|\phi_1|+\cdots+|\phi_{s}|}
                    \phi_1\cdots\phi_s\big(\mathsf{b}\phi_{s+1}\big)\phi_{s+2}
                     \cdots \phi_k \nonumber \\
         & &+ \sum_{0\leq s\leq k-2}(-1)^{|\phi_1|+\cdots+|\phi_{s+1}|}
              \phi_1\cdots\phi_s(\phi_{s+1}\cup\phi_{s+2})\phi_{s+3}\cdots\phi_k
                \label{EqDefCobordBKadei}
\end{eqnarray}
o\`{u} $\mathsf{b}=[m,~]_G$ est le cobord de Hochschild (\`{a} un signe global
pr\`{e}s) et $\cup$ est la multiplication {\em cup} d\'{e}finie avant.

La construction pr\'{e}c\'{e}dente se traduit litt\'{e}ralement au cas o\`{u} on remplace
$\hs^\otimes$ par la cog\`{e}bre colibre cocommutative gradu\'{e}e $\mathsf{S}\hs$:
ici $\Hs_S\subset \Hom(\mathsf{S}\hs,\hs)$,
$\rho'_S:\mathsf{S}\Hs_S\otimes \mathsf{S}\hs\ra \mathsf{S}\hs$ et
$\mu_S:\mathsf{S}\Hs_S\otimes \mathsf{S}\Hs_S\ra \mathsf{S}\Hs_S$
d\'{e}finies par des accolades sym\'{e}triques.

\subsection{Morphismes de formalit\'{e}}
 \label{AppSubSecMorphForm}

\noindent La proposition suivante est un exemple de la th\'{e}orie de
perturbation homologique et \'{e}tait montr\'{e}e pour le cas $G_\infty$ dans
\cite[p.45, Prop. 3.1]{GH03}.

\begin{prop}\label{Theorem 2.1}
 Soit $\Hs$ un espace vectoriel gradu\'{e}, soit $D=\sum_{r=1}^\infty
 D^{[r]}$ une structure $F_\infty$ sur $\Hs$ et soit $\hs$ la cohomologie
 de la codiff\'{e}rentielle $D_1$ d\'{e}cal\'{e}e. Soit $\psi^{[1]}:\hs\ra \Hs$ une
 application HKR, i.e. une injection lin\'{e}aire de degr\'{e} $0$ avec
 $Ker D^{[1]}=Im\psi^{[1]}\oplus Im D^{[1]}$.\\
 Alors pour tout entier positif $n\geq 2$ il existe une application lin\'{e}aire
 $\psi^{[n]}:\hs[1]^{[n]}\ra \Hs[1]$ de degr\'{e} $0$ et une application
 $d^{\prime [n]}:\hs[1]^{[n]}\ra \hs[1]$ de degr\'{e} $1$ telles que
 \begin{enumerate}
  \item $d':=\sum_{n=2}^\infty d^{\prime [n]}$ d\'{e}finit une structure $F_\infty$ sur
       $\hs$.
  \item L'application lin\'{e}aire $\overline{\psi}:\mathbf{F_\infty}\hs\ra\mathbf{F_\infty}\Hs$
  co-induite par $\psi:=\sum_{n=1}^\infty\psi^{[k]}$ en tant que morphisme de cog\`{e}bres
  gradu\'{e}es est un morphisme de
  structures $F_\infty$, i.e.
  $\overline{D}~\overline{\psi}=\overline{\psi}~\overline{d'}$.
 \end{enumerate}
 De plus $d'$ est uniquement d\'{e}termin\'{e} par le choix de $\psi$.
\end{prop}
\begin{prooof}
 Afin de simplifier la notation, on enl\`{e}vera le symbol $\circ$ pour la
 composition des applications lin\'{e}aires dans cette
 d\'{e}monstration et dans la suivante.
 Pour un choix arbitraire des $\psi^{[n]}$ $(n\geq 2)$ et des $d^{\prime
 [n]}$ ($n\geq 2$) on d\'{e}finit
 \[
   \overline{P}(\psi,d'):=\overline{D}~\overline{\psi}-\overline{\psi}~\overline{d'}
   ~\mbox{ et }~
   \overline{Q}(d'):= \overline{d'}~\overline{d'}.
 \]
Alors puisque $\overline{D}$ et $\overline{d'}$ sont des $F$-cod\'{e}rivations gradu\'{e}es et
$\overline{\psi}$ est un morphisme de $F$-cog\`{e}bres il vient que
$\overline{P}(\psi,d')$ est une $F$-cod\'{e}rivation le long de $\overline{\psi}$. De
plus, l'\'{e}quation $\overline{D}~\overline{D}=0$ entra\^{\i}ne l'identit\'{e}
suivante:
\[
   \overline{D}~\overline{P}(\psi,d')~+~\overline{P}(\psi,d')~\overline{d'}~+~
                   \overline{\psi}~\overline{Q}(d')=0.~~~~~~~~~~~~~~~(*)
\]
On va construire par r\'{e}currence $\psi$ et $d'$ telles que
$\overline{P}(\psi,d')$ et $\overline{Q}(d')$ s'annulent. D'abord on note
$\hs[1]^{[\leq n]}$ la filtration $\oplus_{k=0}^n\hs[1]^{[k]}$ de $\mathbf{F_\infty}\hs$ et
les applications tronqu\'{e}es $\psi^{[\leq n]}:=\psi^{[1]}+\cdots+\psi^{[n]}$
et $d^{\prime [\leq n]}:=d^{\prime [2]}+\cdots+d^{\prime [n]}$. De la construction
des morphismes et cod\'{e}rivations co-induites l'on
constate les faits g\'{e}n\'{e}raux suivants: 1. $\psi^{[k]}$ et $d^{\prime [k]}$
s'annulent sur $\hs[1]^{[l]}$ si $k>l$, donc $\overline{D}$,
$\overline{d'}$ et $\overline{\psi}$ pr\'{e}servent les filtrations.
2. Pour les restrictions il vient
\[
\left.\left(
   \overline{\psi} = \overline{\psi^{[\leq n]}}
= \left(\overline{\psi^{[\leq n-1]}}+\psi^{[n]}\right)\right)\right|_{\hs[1]^{[\leq n]}}
~~\mbox{ et }~~
\left.\left(\overline{d'} =
\overline{d^{\prime[\leq n]}}\right)\right|_{\hs[1]^{[\leq n]}},
\]
bien que la d\'{e}pendence de $\overline{\psi}$ de $\psi$ ne soit pas lin\'{e}aire
contrairement aux cod\'{e}rivations. 3. On a
$\overline{d^{\prime}}\big(\hs[1]^{[\leq n+1]}\big)\subset
\hs[1]^{[\leq n]}$ car $d^{\prime [1]}=0$. 4. Si une cod\'{e}rivation
$\overline{p}$:
$\mathbf{F_\infty}\hs\ra\mathbf{F_\infty}\Hs$ le long d'un morphisme s'annule sur
$\hs[1]^{[\leq n]}$, alors
$\overline{p}(\hs[1]^{[n+1]})=p^{[n+1]}(\hs[1]^{[n+1]})\subset
\Hs=\Hs^{[1]}$ car elle est d\'{e}termin\'{e}e par ses projections sur $\Hs$.\\
Pour $n=1$ on trouve
$\overline{P}(\psi,d')|_{\hs[1]^{[1]}}=D^{[1]}\psi^{[1]}=0$ par la
d\'{e}finition de $\psi^{[1]}$ et le fait que $d^{\prime [1]}=0$. Puisque tout $d^{\prime [2]}$
envoie $\hs[1]^{[\leq 2]}$ sur $\hs[1]^{[1]}=\hs[1]$, sur lequel tout $d'$ s'annule,
il vient que $\overline{Q}(d')$ s'annule sur $\hs[1]^{[\leq 2]}$.\\
 Supposons, par r\'{e}currence,
qu'il y ait $\psi^{[\leq n]}$ et $d^{\prime [\leq n]}$ tel que
$\overline{P}(\psi^{[\leq n]},d^{\prime [\leq n]})$ s'annule sur $\hs[1]^{[\leq n]}$
et $\overline{Q}(d^{\prime [\leq n]})$ s'annule sur $\hs[1]^{[\leq n+1]}$.
Alors pour tout choix $\psi^{[n+1]}:\hs[1]^{[n+1]}\ra \Hs[1]$ et
$d^{\prime [n+1]}:\hs[1]^{[n+1]}\ra \hs[1]$ il vient que (gr\^{a}ce aux faits
2. et 3. ci-dessus)
\beas
  \lefteqn{
  \overline{P}(\psi^{[\leq n+1]},d^{\prime [\leq n+1]})|_{\hs[1]^{[n+1]}}}\\
    & = &\overline{P}(\psi^{[\leq n]},d^{\prime [\leq n]})|_{\hs[1]^{[n+1]}}
      +D^{[1]}\psi^{[n+1]}-\psi^{[1]}d^{\prime [n+1]}.~~~~(**)
\eeas
D'un autre c\^{o}t\'{e}, regardons l'identit\'{e} $(*)$ pour $\psi=\psi^{[n]}$ et $d'=d^{\prime [n]}$
restreinte \`{a} $\hs[1]^{[n+1]}$: le terme
$\overline{\psi^{[\leq n]}}~\overline{Q}(d^{\prime [\leq n]})$ s'y annule par hypoth\`{e}se de
r\'{e}currence, ainsi que le terme
$\overline{P}(\psi^{[\leq n]},d^{\prime [\leq n]})\overline{d^{\prime [\leq n]}}$
gr\^{a}ce au fait 3. ci-dessus et l'hypoth\`{e}se de r\'{e}currence.
Finalement, gr\^{a}ce au fait 4. ci-dessus,
La restriction $\overline{P}(\psi^{[\leq n]},d^{\prime [\leq n]})|_{\hs[1]^{[n+1]}}$
prend ses valeurs dans $\Hs[1]$, et donc $(*)$ se r\'{e}duit \`{a}
\[
   0=\overline{D}~
\overline{P}(\psi^{[\leq n]},d^{\prime [\leq n]})|_{\hs[1]^{[n+1]}}
    =D^{[1]}
    \overline{P}(\psi^{[\leq n]},d^{\prime [\leq n]})|_{\hs[1]^{[n+1]}}.
\]
Par cons\'{e}quent, $\overline{P}(\psi^{[\leq n]},d^{\prime [\leq n]})|_{\hs[1]^{[n+1]}}$ prend ses valeurs
dans $Ker~D^{[1]}=Im~\psi^{[1]}\oplus Im~D^{[1]}$, et
on peut donc trouver un unique $d^{\prime [n+1]}$ et un $\psi^{[n+1]}$ tels que
le membre de gauche de l'\'{e}quation $(**)$ s'annule.\\
De plus, $\overline{Q}(d^{\prime [\leq n+1]})$ restreinte \`{a}
$\hs[1]^{[\leq n+1]}$ s'annule, car $\overline{Q}(d^{\prime [\leq n]})$ s'y
annule et $\overline{d^{\prime [\leq n]}}d^{\prime [n+1]}=0$
(car $d^{\prime [1]}=0$) et
$d^{\prime [n+1]}\overline{d^{\prime [\leq n]}}(\hs[1]^{[\leq n+1]})=0$
(gr\^{a}ce au fait 3.). Finalement, regardons l'identit\'{e} $(*)$ pour
$\psi=\psi^{[\leq n+1]}$ et $d'=d^{\prime [\leq n+1]}$ restreinte
\`{a} $\hs[1]^{[n+2]}$: d'apr\`{e}s ce qui pr\'{e}c\`{e}de
$\overline{P}(\psi^{[\leq n+1]},d^{\prime [\leq n+1]})$ s'annule sur
$\hs[1]^{[\leq n+1]}\supset \overline{d'}(\hs[1]^{[n+2]})$, donc le terme
$\overline{P}(\psi^{[\leq n+1]},d^{\prime [\leq n+1]})~\overline{d^{\prime [\leq n+1]}}$
s'annule sur $\hs[1]^{[n+2]}$. De plus, gr\^{a}ce au fait 4. les restrictions de
$\overline{P}(\psi^{[\leq n+1]},d^{\prime [\leq n+1]})$ et de
$\overline{\psi^{[\leq n+1]}}~\overline{Q}(d^{\prime [\leq n+1]})$ \`{a} $\hs[1]^{[n+2]}$
prennent leurs valeurs dans $\Hs$.
Alors, l'identit\'{e} $(*)$ a la forme
\[
  0=D^{[1]}\overline{P}(\psi^{[\leq n+1]},d^{\prime [\leq n+1]})|_{\hs[1]^{[n+2]}}
      +\psi^{[1]}\overline{Q}(d^{\prime [\leq n+1]})|_{\hs[1]^{[n+2]}},
\]
et puisque $Im~\psi^{[1]}\cap Im~D^{[1]}=\{0\}$ il vient que
$\overline{Q}(d^{\prime [\leq n+1]})$ s'annule sur $\hs[1]^{[n+2]}$, ce
qui termine la r\'{e}currence.
\end{prooof}

\noindent On remarque que la cod\'{e}rivation $\overline{d}$ de
$\mathbf{F_\infty}\hs$ co-induite par $d^{\prime [2]}=:d:\hs[1]^{[2]}\ra\hs[1]$
est de carr\'{e} $0$. D'apr\`{e}s le cas $n=2$ de la d\'{e}monstration pr\'{e}c\'{e}dente elle est donn\'{e}e par
la formule
\[
      \psi^{[1]}d=\psi^{[1]}d^{\prime [2]}=D^{[2]}\overline{\psi^{[1]}}.
\]
Evidemment, $d$ d\'{e}finit \'{e}galement une structure $F_\infty$ sur $\hs$.

Pour la construction des star-produits, il est souhaitable de transformer $d'$ r\'{e}sultant de la proposition \ref{Theorem 2.1}
en son premier terme $d$ :
on constate d'abord que l'espace $\Hom(\mathbf{F_\infty}\hs,\hs[1])$ est cofiltr\'{e} par la
filtration de $\mathbf{F_\infty}\hs$, i.e.
\[
\Hom(\mathbf{F_\infty}\hs,\hs[1])^{[\geq n+1]}
:=\big\{\phi\in \Hom(\mathbf{F_\infty}\hs,\hs[1])~|~\phi(\hs^{[\leq n]})=\{0\}\big\}
\]
A l'aide du crochet $[~,~]_F$ (o\`{u} $F=G,~NR,~H,~T$), on voit que
l'application $\phi\mapsto [\phi,d]_F$ est une codiff\'{e}rentielle
(i.e. de carr\'{e} $0$) sur tout espace
$\Hom(\mathbf{F_\infty}\hs,\hs[1])^{[\geq n+1]}$.
On en d\'{e}duit le crit\`{e}re suffisant pour transformer $d'$ en $d$
(montr\'{e} dans \cite{Tam98p} et \cite[p.48, Prop. 4.1]{GH03} pour le cas $G_\infty$):

\begin{prop}\label{Theorem 3.1} Soit $\hs$ un espace vectoriel gradu\'{e} et
$d'=\sum_{r=2}^\infty d^{\prime [r]}$ une structure $F_\infty$ sur $\hs$.
Si la cohomologie de $\big(\Hom(\mathbf{F_\infty}\hs,\hs[1])^{[\geq 2]},
[~,d]\big)$ s'annule, alors pour tout entier positif $n\geq 2$
il existe des applications
$\psi^{\prime [n]}:\hs[1]^{[n]}\ra\hs[1]$ de degr\'{e} $0$ telles que le morphisme
$\overline{\psi'}:\mathbf{F_\infty}\hs\ra\mathbf{F_\infty}\hs$
de $F$-cog\`{e}bres co-induit par $\psi':=\sum_{n=1}\psi^{\prime [n]}$
(avec $\psi^{\prime [1]}=\id_{\hs[1]}$)
est un morphisme $F_\infty$, i.e.
\[
   \overline{\psi'}~\overline{d}=\overline{d'}~\overline{\psi'}.
\]
\end{prop}
\begin{prooof}
 On utilisera les abbr\'{e}viations $\hs[1]^{[n]}$ et $\hs[1]^{[\leq n]}$ de la d\'{e}mon\-stra\-tion
 de la proposition \ref{Theorem 2.1} pr\'{e}c\'{e}dente et $\psi^{\prime [\leq n]}:=
 \psi^{\prime [1]}+\cdots +\psi^{\prime [n]}$. Pour un choix arbitraire
 des $\psi^{\prime [n]}$ on d\'{e}finit
 \[
    \overline{R}(\psi'):=\overline{\psi'}~\overline{d}-\overline{d'}~\overline{\psi'},
 \]
 qui est une $F$-cod\'{e}rivation gradu\'{e}e:
 $\mathbf{F_\infty}\hs\ra \mathbf{F_\infty}\hs$ de degr\'{e} $1$ le long de
 $\overline{\psi'}$. Puisque $\overline{d}~\overline{d}=0$ et
 $\overline{d'}~\overline{d'}=0$ on obtient l'identit\'{e}
 \[
    ~~~~~~~~~~~~~~~~~~~~~~~~~
    \overline{d'}~\overline{R}(\psi')~+\overline{R}(\psi')~\overline{d}=0.
    ~~~~~~~~~~~~~~~~~~~~~~~~~(*)
 \]
 Puisque $d^{\prime [1]}=0=d^{[1]}$ il vient que
 $\overline{R}(\psi')(\hs[1]^{[\leq n+1]})\subset \hs[1]^{[\leq n]}$ et
 donc $\overline{R}(\psi')|_{\hs[1]^{[\leq n+1]}}$ co\"{\i}ncide avec
 $\overline{R}(\psi^{\prime [\leq n]})|_{\hs[1]^{[\leq n+1]}}$.
 Le but est de construire les $\psi^{[n]}$ par r\'{e}currence telle que
 $\overline{R}(\psi^{\prime [\leq n]})|_{\hs[1]^{[\leq n+1]}}=0$.\\
 Puisque $d$ et $d^{\prime}$ s'annulent sur $\hs[1]$ et $d=d^{\prime [2]}$
 il vient que $\overline{R}(\psi^{\prime [\leq 1]})$ s'annule sur
 $\hs^{[\leq 2]}$.\\
 Supposons que l'on a choisi $\psi^{\prime [2]},\ldots,\psi^{\prime [n]}$
 telles que
 $\overline{R}(\psi^{\prime [\leq n]})|_{\hs[1]^{[\leq n+1]}}=0$. Il vient que
 \beas
   \lefteqn{\overline{R}(\psi^{\prime [\leq n+1]})|_{\hs[1]^{[n+2]}}
     =}\\
     & &\overline{R}(\psi^{\prime [\leq n]})|_{\hs[1]^{[n+2]}}
      +\psi^{\prime [n+1]}~\overline{d}|_{\hs[1]^{[n+2]}}
       - d\circ_F \psi^{\prime [n+1]}|_{\hs[1]^{[n+2]}}\\
     & & =\overline{R}(\psi^{\prime [\leq n]})|_{\hs[1]^{[n+2]}}
      +[\psi^{\prime [n+1]},d]_{F}|_{\hs[1]^{[n+2]}},
        ~~~~~~~~~~~~~~~~~~~~~(**)
 \eeas
 car les mon\^{o}mes en $\psi^{\prime [1]}\!\!\!,\ldots,\psi^{\prime [n+1]}$ de
 $\overline{\psi^{\prime [\leq n+1]}}|_{\hs[1]^{[n+2]}}$ qui contiennent
 $\psi^{\prime [n+1]}$ ne contiennent que $\psi^{\prime [n+1]}$ et
 $\psi^{\prime [1]}=\id$, d'o\`{u} l'apparition du terme
 $- d\circ_F \psi^{\prime [n+1]}|_{\hs[1]^{[n+2]}}$ dans $(**)$. Comme
 dans la d\'{e}monstration pr\'{e}c\'{e}dente, on cons\-tate que
 $\overline{R}(\psi^{\prime [\leq n]})|_{\hs[1]^{[n+2]}}$ prend ses
 valeurs dans $\hs[1]$. On regarde la restriction de l'identit\'{e} $(*)$
 pour $\psi^{\prime [\leq n]}$ \`{a} l'espace $\hs[1]^{[n+3]}$. Puisque $\overline{d'}$
 s'annule sur $\hs[1]$, la composition
 $\overline{d'}~\overline{R}(\psi^{\prime [\leq n]})|_{\hs[1]^{[n+3]}}$ ne d\'{e}pend que
 de la restriction de $\overline{R}(\psi^{\prime [\leq n]})$ \`{a}
 $\hs[1]^{[n+2]}$, et comme $\overline{d'}$ envoie
 $\hs[1]^{[n+3]}$ dans $\hs[1]^{[\leq n+2]}$ il vient
 \[
  0= d\circ_F
   \big(\overline{R}(\psi^{\prime [\leq n]})|_{\hs[1]^{[n+2]}}\big)
      +
   \big(\overline{R}(\psi^{\prime [\leq n]})|_{\hs[1]^{[n+2]}}\big)
        \circ_F d.
 \]
 Il vient que
 $\big(\overline{R}(\psi^{\prime [\leq n]})|_{\hs[1]^{[n+2]}}\big)$ se trouve
 dans $\Hom(\mathbf{F_\infty}\hs,\hs[1])^{[\geq 2]}$ et
 dans le noyau de la codiff\'{e}rentielle $[~,d]_F$, et puisque la cohomologie
 \'{e}tait nulle, on trouve $\psi^{\prime [n+1]}$ pour annuler le membre de
 gauche dans l'\'{e}quation $(**)$.
\end{prooof}

\end{appendix}

\end{document}